\documentclass[10pt]{amsart}
\pdfoutput=1
\usepackage{ mathrsfs }
\usepackage{amsmath, amsthm, amssymb,slashed,stmaryrd}
\usepackage{mathtools}
\usepackage{ifpdf}
\usepackage[pdftex]{graphicx}
\usepackage{tikz}
\usetikzlibrary{}
\usetikzlibrary{decorations.markings,matrix,arrows,decorations.pathmorphing,backgrounds,patterns,shapes.geometric,calc}

\usetikzlibrary{matrix,arrows,calc}
\usepackage[all]{xy}

\usepackage[pdftex,plainpages=false,hypertexnames=false,pdfpagelabels]{hyperref}
\usepackage{amsthm}

\theoremstyle{definition}
 \newtheorem{thm}{Theorem}[section]
    
 \newtheorem{cor}[thm]{Corollary}
 \newtheorem{lem}[thm]{Lemma}

 \newtheorem{prop}[thm]{Proposition}
 
 \newtheorem{defn}[thm]{Definition}
  
 \newtheorem{notation}[thm]{Notation}
 \newtheorem{ex}[thm]{Example}
  
 \newtheorem*{thm*}{Theorem}
 \theoremstyle{remark}
 \newtheorem{rmk}[thm]{Remark}

\def\beq{\begin{eqnarray}}
\def\eeq{\end{eqnarray}}
 \newcommand{\bp}{\begin{proof}[Proof]}
 \newcommand{\ep}{\end{proof}}

\def\susp{{{u}}}

\def\id{{{\rm id}}}
\def\Top{{\sf{Top}}}
\def\pr{{\rm{pr}}}
\def\Fr{{\rm{Fr}}}
\def\C{{\mathbb{C}}}
\def\X{{\mathfrak{X}}}
\def\V{{\mathcal{V}}}
\def\bX{{\mathbb{X}}}
\def\bY{{\mathbb{Y}}}
\def\bW{{\mathbb{W}}}
\def\Y{{\mathcal{Y}}}

\def\Z{{\mathbb{Z}}}

\def\End{{\rm End}}

\def\HH{{\rm H}}
\def\Rep{{\rm Rep}}

\def\cl{{\rm{cl}}}
\def\Pin{{\rm Pin}}
\def\Spin{{\rm Spin}}
\def\O{{\rm O}}
\def\Fred{{\rm Fred}}
\def\H{{\mathcal H}}

\def\Cl{{\rm Cl}}
\def\cCl{\mathbb{C}{\rm l}}
\def\U{{\rm U}}
\def\PU{{\rm PU}}
\def\PO{{\rm PO}}
\def\Pow{{\mathbb{P}}}
\def\cH{{\mathcal{H}}}
\def\cL{{\mathcal{L}}}

\def\R{{\mathbb R}}
\def\N{{\mathbb N}}
\def\Pow{{\mathbb P}}
\def\bS{{\mathbb S}}
\def\O{{\rm O}}
\def\SO{{\rm SO}}
\def\KO{{\rm KO}}
\def\KU{{\rm KU}}
\def\KR{{\rm KR}}
\def\Sym{{\rm Sym}}
\def\Pic{{\sf Pic}}

\def\K{{\rm K}}
\def\SW{\tau}
\def\Th{{\rm Th}}

\def\Tw{{\rm Tw}}

\def\Twist{{\rm Twist}}

\def\Gerbe{{\rm Gerbe}}

\def\res{{\rm res}}
\def\Stack{{\sf TopStack}}

\def\Grpd{{\sf TopGrpd}}
\def\pt{{*}}
\def\twocommute{\ensuremath{\rotatebox[origin=c]{30}{$\Rightarrow$}}}

\newcommand{\op}{{\sf{op}}}   
\newcommand{\sq}{/\!\!/}

\begin{document}
\title[Power operations preserve Thom classes in twisted K-theory]{Power operations preserve Thom classes in twisted equivariant Real K-theory}

\author{Daniel Berwick-Evans and Meng Guo}
\begin{abstract}
We construct power operations for twisted KR-theory of topological stacks. 
Standard algebraic properties of Clifford algebras imply that these power operations preserve universal Thom classes. As a consequence, we show that the twisted Atiyah--Bott--Shapiro orientation commutes with power operations. 

\end{abstract}


\maketitle 

\setcounter{tocdepth}{1}
\tableofcontents

\section{Introduction}
In this paper we construct universal Thom classes in twisted equivariant Real K-theory, 
\beq\label{eq:univThom}
&&\Th_{p,q}\in \KR^{\SW_{p,q}}(\R^{p,q}\sq\O_{p,q} ),\qquad p,q\in \N.
\eeq
Above, $\SW_{p,q}$ is the \emph{$(p,q)^{\rm th}$ Thom twist}, and $\R^{p,q}$ is the Real space (in Atiyah's sense~\cite{Atiyahreal}) given by $p$-copies of the trivial $C_2$-representation and $q$-copies of the sign representation. The standard action of the orthogonal group $\O_{p+q}$ on $\R^{p+q}$ determines an action of a Real Lie group $\O_{p,q}$ on $\R^{p,q}$, and  $\R^{p,q}\sq \O_{p,q}$ is the associated Real quotient stack, see~\S\ref{sec:Realgroups}. The classes~\eqref{eq:univThom} are \emph{universal} in the sense that a map $\X\to \pt\sq \O_{p,q}$ classifying a Real vector bundle on a Real stack $\X$ determines a pullback in $\KR$-theory that recovers the usual Thom class, e.g., in the flavors from \cite[\S12]{ABS}, \cite[Remark 3.2.22]{ST04}, \cite[\S3.6]{FHTI}; see~\S\ref{sec:univerthomintro}.

Next we construct twisted $\KR$-power operations 
\beq\label{eq:powermain}
\Pow_k \colon \KR^\sigma(\X) \to \KR^{\Pow_k\sigma }(\X^{\times k}\sq \Sigma_k),\qquad k\in \mathbb{N}
\eeq
for any Real topological stack~$\X$ and twist~$\sigma$, where $\X^{\times k}\sq \Sigma_k$ is the $k^{\rm th}$ symmetric power of the stack~$\X$. Concretely, the operations~\eqref{eq:powermain}
 come from tensor powers of Fredholm operators. 

Our main result is a compatibility between the constructions~\eqref{eq:univThom} and~\eqref{eq:powermain}.
 
\begin{thm} \label{thm:A}
Twisted power operations preserve universal Real Thom classes:
\beq
\begin{tikzpicture}[baseline=(basepoint)];
\node (A) at (0,0) {$\KR^{\SW_{p,q}}(\R^{p,q}\sq \O_{p,q})$};
\node (B) at (4.75,0) {$\KR^{\Pow_k\SW_{p,q}}((\R^{p,q})^{\times k}\sq (\O_{p,q}\wr \Sigma_k))$};
\node (C) at (10,0) {$\KR^{\SW_{kp,kq}}(\R^{kp,kq}\sq \O_{kp,kq})$};
\node (D) at (0,-.75) {$\Th_{p,q}$};
\node (E) at (5,-.75) {$\Pow_k(\Th_{p,q})=\res(\Th_{kp,kq})$};
\node (F) at (10,-.75) {$\Th_{kp,kq}.$};
\draw[->] (A) to node [above] {$\Pow_k$} (B);
\draw[->] (C) to node [above] {$\res$} (B);
\draw[|->] (D) to (E);
\draw[|->] (F) to (E);
\path (0,-.75) coordinate (basepoint);
\end{tikzpicture}\nonumber
\eeq
Above, $\res$ is restriction along the map of stacks $(\R^{p,q})^{\times k}\sq (\O_{p,q}\wr\Sigma_k)\hookrightarrow \R^{kp,kp}\sq \O_{kp,kq}$ determined by the  inclusion of the wreath product as a subgroup $\O_{p,q}\wr\Sigma_k\hookrightarrow\O_{kp,kq}$. 
\end{thm}

\begin{rmk}
Throughout, $\K$-theory groups will be compactly supported; see Remark~\ref{rmk:compactsupport}.
\end{rmk}

Theorem~\ref{thm:A} specializes to twisted (complex) KU-theory and twisted (real) KO-theory. 

\begin{cor} \label{cor:A}
Power operations preserve universal $\KU$-and $\KO$-Thom classes:
\beq
\begin{tikzpicture}[baseline=(basepoint)];
\node (A) at (0,0) {$\KU^{\SW_{n}}(\R^{n}\sq \O_{n})$};
\node (B) at (4.75,0) {$\KU^{\Pow_k\SW_{n}}((\R^{n})^{\times k}\sq (\O_{n}\wr \Sigma_k))$};
\node (C) at (10,0) {$\KR^{\SW_{kn}}(\R^{kn}\sq \O_{kn})$};
\node (AA) at (0,-.75) {$\KO^{\SW_{n}}(\R^{n}\sq \O_{n})$};
\node (BB) at (4.75,-.75) {$\KO^{\Pow_k\SW_{n}}((\R^{n})^{\times k}\sq (\O_{n}\wr \Sigma_k))$};
\node (CC) at (10,-.75) {$\KO^{\SW_{kn}}(\R^{kn}\sq \O_{kn})$};
\node (D) at (0,-1.5) {$\Th_{n}$};
\node (E) at (5,-1.5) {$\Pow_k(\Th_{n})=\res(\Th_{kn})$};
\node (F) at (10,-1.5) {$\Th_{kn}$};
\draw[->] (A) to node [above] {$\Pow_k$} (B);
\draw[->] (C) to node [above] {$\res$} (B);
\draw[->] (AA) to node [above] {$\Pow_k$} (BB);
\draw[->] (CC) to node [above] {$\res$} (BB);
\draw[|->] (D) to (E);
\draw[|->] (F) to (E);
\path (0,-.75) coordinate (basepoint);
\end{tikzpicture}\nonumber
\eeq
where now $\res$ is restriction along the map of stacks $(\R^{n})^{\times k}\sq (\O_{n}\wr\Sigma_k)\hookrightarrow \R^{kn}\sq \O_{kn}$ determined by the  inclusion of the wreath product as a subgroup $\O_{n}\wr\Sigma_k\hookrightarrow\O_{kn}$.
\end{cor}

In words, Theorem~\ref{thm:A} and Corollary~\ref{cor:A} state that a power operation applied to a Thom class can be recovered from the restriction of another Thom class. The proof boils down to the $\O_{p,q}\wr\Sigma_k$-equivariant isomorphism of Real super algebras, 
\beq\label{eq:maintensoriso}
\cCl_{p,q}^{\otimes k}\simeq \cCl_{kp,kq}.
\eeq
The action on the left side is the $\O_{p,q}$-action on each tensor factor $\cCl_{p,q}=\cCl(\R^{p,q})$ (from functoriality of the Clifford construction) together with the permutation of the tensor factors; the action on the right side of~\eqref{eq:maintensoriso} is through restriction along $\O_{p,q}\wr\Sigma_k=\O_{p,q}^{\times k}\rtimes \Sigma_k\hookrightarrow \O_{kp,kq}$ to block-diagonal orthogonal matrices and permutations of the blocks. 
Hence, Theorem~\ref{thm:A} and the isomorphism~\eqref{eq:maintensoriso} identify a completely algebraic origin for the compatibility between power operations and the (twisted) Real Atiyah--Bott--Shapiro orientation. 

We expand on this compatibility by computing the effect of power operations on the twisted topological pushforward determined by the Thom classes~\eqref{eq:univThom}. A simple version of this twisted pushforward takes as input a compact manifold $X$ with Real Spin$^c$-structure (see Definition~\ref{defn:spinstructure}), and a Real $G$-action by a compact Lie group that is not required to preserve the Real Spin$^c$- structure. This leads to a wrong-way map (see~\S\ref{sec:toppush}-\ref{sec:ABS})
\beq\label{Eq:twistedABS0}
&&\KR(X\sq G)\xrightarrow{\pi_!^{\rm top}} \KR^{\tau_{TX}^{-1}}(\pt\sq G),\qquad \pi\colon X\sq G\to \pt\sq G
\eeq
where the inverse Thom twisting of the tangent bundle $\tau_{TX}^{-1}$ is a pair of data: (i) a grading $G\to \Z/2$ that encodes when the $G$-action preserves or reverses the orientation on $X$ and (ii) a Real $\U(1)$-central extension~$\widehat{G}$ that lifts the $G$-action on $X$ to a Real $\widehat{G}$-action on the Clifford linear spinor bundle of~$X$. 



\begin{thm}[{Theorem~\ref{thm:ABS}}]\label{thm:3}
For all $k\in \N$, the twisted topological pushforward~\eqref{Eq:twistedABS0} is compatible with the~$k^{\rm th}$ total power operation, i.e., the diagram commutes
\beq\label{eq:ABSpow}
\begin{tikzpicture}[baseline=(basepoint)];
\node (A) at (0,0) {$\KR(X\sq G)$};
\node (B) at (5,0) {$\KR^{\tau_{TX}^{-1}}(\pt\sq G)$};
\node (C) at (0,-1.25) {$\KR(X^{\times k}\sq G\wr \Sigma_k)$};
\node (D) at (5,-1.25) {$\KR^{\tau^{-1}_{TX^{\times k}}}(\pt\sq G\wr \Sigma_k).$};
\draw[->] (A) to node [above] {$\pi_!^{\rm top}$} (B);
\draw[->] (B) to node [right] {$\Pow_k$} (D);
\draw[->] (A) to node [left] {$\Pow_k$} (C);
\draw[->] (C) to node [below] {$\pi_!^{\rm top}$} (D);
\path (0,-.75) coordinate (basepoint);
\end{tikzpicture}
\eeq
The analogous statement also holds in KU- and KO-theory. 
\end{thm}

Theorem~\ref{thm:ABS} below further verifies the families generalization of Theorem~\ref{thm:3}: power operations  are compatible with the twisted topological pushforward along a Real $G$-equivariant proper submersion $\pi\colon X\to Y$ (see Definition~\ref{defn:ABS}). This recovers Theorem~\ref{thm:3} when $Y=\pt$.

\begin{rmk}
We may identify $\KR^{\tau_{TX}^{-1}}(\pt\sq G)$ with $\Rep^{\tau_{TX}^{-1}}(G)$, the Grothendieck group of Real, projective, graded $G$-representations in Clifford modules for the grading and central extension specified by the twisting $\tau_{TX}^{-1}$. 
By the equivariant index theorem, the image of $1\in \KR(X\sq G)$ along~\eqref{Eq:twistedABS0} is the kernel of the Real Clifford linear Dirac operator on~$X$ viewed as a graded projective Real $G$-representation for the grading $G\to \Z/2$ and central extension $\widehat{G}$. A $G$-invariant spin structure (if it exists) fixes an isomorphism between the twisting $\tau_{TX}^{-1}$ and the degree twisting, resulting in the usual shift by the dimension of $X$ in~\eqref{Eq:twistedABS0}. 
\end{rmk}
After some contextualizing remarks below, the remainder of this introduction outlines the constructions leading to Theorems~\ref{thm:A} and~\ref{thm:3}.

\begin{rmk}
We recall that the K-theory of quotient stacks is constructed to coincide with (the non-completed) equivariant K-theory of $G$-spaces: by definition
\beq\label{eq:notBorel}
&&\KR^\tau_G(X)=\KR^\tau(X\sq G), \quad \KO^\tau_G(X)=\KO^\tau(X\sq G),\quad \KU^\tau_G(X)=\KU^\tau(X\sq G)
\eeq
for any compact Lie group $G$ acting on a topological space $X$ with twisting~$\tau$. 
\end{rmk}

\begin{rmk}\label{rmk:suspension}
The pullback along the quotient map $q\colon \R^{p,q}\to \R^{p,q}\sq \O_{p,q}$ 
\beq\label{eq:suspension}
\KR^{\SW_{p,q}}(\R^{p,q}\sq \O_{p,q})\xrightarrow{q^*} \KR^{q^*\SW_{p,q}}(\R^{p,q})\simeq \KR^{p,q}(\R^{p,q}),\qquad \Th_{p,q}\mapsto \susp_{p,q},
\eeq
sends the Thom class~\eqref{eq:univThom} to the $\KR$-suspension class $\susp_{p,q}$ of the $C_2$-representation sphere of~$\R^{p,q}$. Equivalently,~\eqref{eq:suspension} restricts the Thom class of the vector bundle $\R^{p,q}\sq \O_{p,q}\to \pt\sq \O_{p,q}$ to the fiber at the basepoint $\pt\to \pt\sq \O_{p,q}$. Ignoring twistings and $C_2$-equivariance,~\eqref{eq:suspension} recovers the definition of orientability from~\cite[\S12]{ABS}: a Thom class is defined as a class that restricts to the suspension class on each fiber. Turning~\eqref{eq:suspension} around, the Thom twisting measures the failure of the suspension class $\susp_{p,q}$ to descend to an $\O_{p,q}$-equivariant class. 
\end{rmk}

\begin{rmk} \label{rmk:introtwist}
In Theorem~\ref{thm:A} and Corollary~\ref{cor:A}, Thom twistings are essential to even parse a compatibility between power operations and Thom classes. Algebraically, this stems from the fact that the wreath products $\Spin_n\wr\Sigma_k$ are not natural subgroups of $\Spin(nk)$; this stands in contrast to the wreath products of orthogonal groups as subgroups of larger orthogonal groups. Geometrically, twistings arise from the fact that the action of the symmetric group~$\Sigma_k$ on the $k^{\rm th}$ power~$(\R^n)^{\times k}\simeq \R^{nk}$ need not preserve a spin structure, or even an orientation. Hence, power operations applied to Thom classes are unavoidably twisted. 
\end{rmk}

\begin{rmk}\label{rmk:introtwist2}
Building on the previous remark, twistings are also inevitable when applying power operations to classes in nonzero cohomological degree. To summarize, the $k^{\rm th}$ power operation on a $\cCl_n$-linear Fredholm operator yields a~$\cCl_n^{\otimes k}$-linear $\Sigma_k$-equivariant Fredholm operator. However, the $\Sigma_k$-equivariance involves the $\Sigma_k$-action on~$\cCl_n^{\otimes k}$ that permutes the tensor factors, which introduces a twisting. Indeed, this twisting is a specialization of the Thom twisting for the $\Sigma_k$-action on~$(\R^n)^{\times k}$. The twisting trivializes when the Clifford algebra is Morita trivial: in degrees $0$ mod $2$ for $\KU$ and degrees $0$ mod $8$ for $\KO$. 
\end{rmk}

\begin{rmk}
When $G$ is trivial, Theorems~\ref{thm:3} and~\ref{thm:ABS} generalize the classical compatibility between power operations and orientations of K-theory as maps of $H_\infty$-spectra~\cite[Proposition VII.7.2]{bmms}. These previous results restrict the Spin$^c$ orientation of $\KU$ to degrees 0 mod 2 and the spin orientation of $\KO$ to degrees 0 mod 8. This restriction is explained by the previous remark: power operations in other degrees necessitate twistings. 
\end{rmk}

\subsection{Representing stacks for twisted KR-theory} 
A reoccurring strategy in this paper is to formulate definitions and constructions using universal objects. This begins with the following stacky distillation of twisted Real K-theory that is somewhat implicit in the literature, e.g., see preprint versions of~\cite{FHTIII} and \cite[Example 1.3.19]{Sati2021EquivariantPI}. 

Let $\cH$ be a complex Hilbert space (typically infinite rank) with the structure of a graded, unitary, Real right $\cCl_{p,q}$-module where $\cCl_{p,q}$ is the $(p,q)^{\rm th}$ Real Clifford algebra; see~\S\ref{sec:RealCliff}. Let $\U^\pm_{\cCl_{p,q}}(\cH)$ denote the Real graded  topological group of unitary endomorphisms of~$\cH$ that commute with the Clifford action. We equip $\U^\pm_{\cCl_{p,q}}(\cH)$ with the compact open (or equivalently, strong~\cite{EspinozaUribe}) topology. Generalizing the classical notion due to Atiyah~\cite{Atiyahreal}, a \emph{Real stack} is a topological stack $\X$ with (weak) $C_2$-action, and a \emph{Real map} is a $C_2$-equivariant map of topological stacks. Consider the Real fibration
\beq\label{Eq:mainfibration}
\Fred_{\cCl_{p,q}}(\H)\sq \PU^\pm_{\cCl_{p,q}}(\H)\to \pt\sq \PU^\pm_{\cCl_{p,q}}(\cH)
\eeq
whose source is the quotient stack for the conjugation action of the projective unitary group on $\Fred_{\cCl_{p,q}}(\cH)$, the space of Real $\cCl_{p,q}$-linear Fredholm operators.

\begin{defn}[{Definitions~\ref{defn:twist} and~\ref{defn:twistedK} and Lemma~\ref{lem:KRcocycle}}]\label{defn:twiustedKthey}
For a Real local quotient stack~$\X$
the \emph{twisted $\KR$-theory} of $\X$ is the set of homotopy classes of Real maps
\beq\label{eq:mainKR}
\X\to \Fred_{\cCl_{p,q}}(\H)\sq \PU^\pm_{\cCl_{p,q}}(\H).
\eeq
Forgetting Real structures, the set of homotopy classes of maps~\eqref{eq:mainKR} is the \emph{twisted $\KU$-theory of $\X$}. For a stack $\X$ the \emph{twisted $\KO$-theory of $\X$} is the set of homotopy classes of Real maps~\eqref{eq:mainKR} for the trivial Real structure on $\X$ (i.e., trivial $C_2$-action). 
\end{defn}

In brief, the stacks in~\eqref{Eq:mainfibration} are representing objects for twisted equivariant KR-theory. To unpack this slightly, a map~\eqref{eq:mainKR} provides data 
\beq\label{eq:KRcocycle}
&&\begin{tikzpicture}[baseline=(basepoint)];
\node (A) at (0,0) {$\X$};
\node (B) at (4,0) {$ \Fred_{\cCl_{p,q}}(\H)\sq \PU^\pm_{\cCl_{p,q}}(\H)$};
\node (C) at (8,0) {$\pt\sq \PU^\pm_{\cCl_{p,q}}(\cH)$};
\draw[->] (A) to node [above] {$F$} (B);
\draw[->] (B) to (C);
\draw[->,bend right=10] (A) to node [below] {$\tau$} (C);
\path (0,0) coordinate (basepoint);
\end{tikzpicture}\iff [F]\in \KR^\tau(\X)
\eeq
yielding the indicated twisted $\KR$ class, where postcomposition with~\eqref{Eq:mainfibration} defines the \emph{twist}~$\tau$. The collection of twists on~$\X$, denoted $\Twist(\X)$, form a symmetric monoidal category. 

To connect with the Freed--Hopkins--Teleman definition of twisted K-theory, for topological groupoids $\bX$ and $\bY$ a map between their underlying stacks is equivalent to a zigzag
\beq\label{eq:zigzag}
\bX\xleftarrow{\sim} \bW\to \bY
\eeq
 of continuous functors where $\bW\xrightarrow{\sim} \bX$ is an essential equivalence~\cite{Pronk}; we review this language in~\S\ref{appen:stacks}. Hence, for a topological groupoid $\bX$ presenting the stack~$\X$, the map of stacks~\eqref{eq:mainKR} is equivalent to a zigzag~\eqref{eq:zigzag} for a continuous functor $\bW\to \Fred_{\cCl_{p,q}}(\H)\sq \PU^\pm_{\cCl_{p,q}}(\H)$. In~\S\ref{sec:Ktheorygrpod1} we show that an object of $\Twist(\X)$ recovers a twisted Hilbert bundle on~$\X$, and in~\S\ref{sec:Ktheorygrpod} we show that a map of stacks~\eqref{eq:mainKR} determines a family of twisted Fredholm operators in the framework of Freed--Hopkins--Teleman~\cite{FHTI}. Adding Real structures identifies Definition~\ref{defn:twiustedKthey} with Freed and Moore's Real twisted K-theory~\cite{FreedMoore,Gomi}.

 \begin{rmk}\label{rmk:compactsupport}
Unless stated otherwise, we will take $\KR^\tau(\X)$ to mean the \emph{compactly supported} $\tau$-twisted $\KR$-theory of $\X$, where compact supports are defined relative to the coarse moduli space of $\X$, see Definition~\ref{sec:coarse}. In concrete terms, the compact support condition demands that the family of Fredholm operators~\eqref{eq:mainKR} is invertible outside a closed substack of~$\X$ whose image in the coarse moduli space is a compact subspace. In particular, this condition is vacuous when the coarse moduli space is compact. Compactly supported $\KR$-theory can also be phrased as relative $\KR$-theory for a pair. For example, the Thom classes~\eqref{eq:univThom} are elements in the $\KR$-cohomology of the $\O_{p,q}$-equivariant pair $(D^{p,q},\partial D^{p,q})$ for the disk $D^{p,q}\subset \R^{p,q}$ and its boundary sphere $\partial D^{p,q}\simeq S^{p+q-1}$. The absolute $\KR$-groups defined below determine the $\KR$-groups of a pair via \cite[Definition 3.15]{FHTI}. 
 \end{rmk}

\subsection{Universal twisted power operations}
Turning to~\eqref{eq:powermain}, power operations are constructed by maps of Real topological stacks (Definition~\ref{defn:powerop})
\beq
\Sym^k(\Fred_{\cCl_{p,q}}(\H))\sq \PU^\pm_{\cCl_{p,q}}(\H))&\simeq & \Fred_{\cCl_{p,q}}(\H)^{\times k}\sq \PU^\pm_{\cCl_{p,q}}(\H)\wr \Sigma_k\nonumber \\
&\to& \Fred_{\cCl_{p,q}^{\otimes k}}(\H^{\otimes k})\sq \PU^\pm_{\cCl_{p,q}^{\otimes k}}(\H^{\otimes k})\label{eq:defofpowerop}
\eeq
where the source is the $k^{\rm th}$ symmetric power whose value on a stack $\X$ is $\Sym^k(\X):=\X^{\times k}\sq \Sigma_k$. The arrow~\eqref{eq:defofpowerop} is induced by the $k^{\rm th}$ tensor power of Fredholm operators. Functoriality of the $k^{\rm th}$ symmetric power then determines~\eqref{eq:powermain}, the $k^{\rm th}$ total power operation. Basic properties of symmetric powers give an essentially formal verification of the usual axioms for power operations~\cite[Chapter VIII,\S1]{bmms}, see Proposition~\ref{prop:powersproperties} below. Furthermore, for ordinary spaces $X$ with the action by a compact Lie group $G$ and trivial twisting, the operation $\Pow_k$ recovers Atiyah's power operations~\cite{AtiyahPow}, see Proposition~\ref{prop:Atiyah}. 

\begin{rmk}
For any $E_\infty$-ring spectrum $R$, general arguments give an action of $\Pic(R)$ on~$R$, see~\cite{ABGHR}.
This permits the construction of twisted power operations in $R$-theory using symmetric powers of an appropriate quotient of $\Pic(R)$ by $R$. The source of~\eqref{Eq:mainfibration} is an approximation to this quotient for $R=\KR$ (in $C_2$-equivariant spectra). However, the general setup only captures Borel $G$-equivariant KR-theory, whereas the language of topological stacks above encodes the correct (non-Borel) equivariant refinement as in~\eqref{eq:notBorel}. In other words, the equivariant power operations induced by~\eqref{eq:defofpowerop} land in genuine equivariant K-theory rather than the Atiyah--Segal completed version~\cite{Atiyah_Segal_complete}.  Through their connection with the representation theory of Lie groups, these genuine equivariant operations are both more powerful and more computable~\cite{AtiyahPow}. 
\end{rmk}

\subsection{Universal Thom classes}\label{sec:univerthomintro}
We refer to Example~\ref{ex:CnUn} for the complete definition of the Real stack $\R^{p,q}\sq \O_{p,q}$. 
In the following, $\bS_{p,q}$ is the spinor representation (viewed as a projective $\O_{p,q}$-representation) built from $\cCl_{p,q}$ regarded as a $\cCl_{p,q}$-$\cCl_{p,q}$ bimodule, see~\S\ref{sec:spinor}.

\begin{defn}[{Definitions~\ref{defn:Thomtwist} and~\ref{defn:Thomclass}}]
The universal Thom cocycles are maps of stacks
\beq\label{eq:Thomtwists}
\widehat{\Th}_{p,q}\colon  \R^n\sq \O_{p,q}\to \End_{\cCl_{p,q}}(\bS_{p,q})\sq \PU^\pm(\bS_{p,q})
\eeq
determined by the map $\R^{p,q}\subset \cCl_{p,q}\to \End_{\cCl_{p,q}}(\bS_{p,q})$ given by (left) Clifford multiplication and the homomorphism $\O_{p,q}\to \PU^\pm_{\cCl_{p,q}}(\bS_{p,q})$ given by the spinor representation. The \emph{universal Thom class}~\eqref{eq:univThom} is the twisted $\KR$-class underlying~\eqref{eq:Thomtwists}. 
\end{defn}

Theorem~\ref{thm:A} computes the image of~\eqref{eq:Thomtwists} under the tensor power operation~\eqref{eq:defofpowerop}, and the isomorphism~\eqref{eq:maintensoriso} provides the comparison with the desired higher rank Thom class.

The universal Thom isomorphism for~\eqref{eq:Thomtwists} is verified in Proposition~\ref{prop:Thomiso}, following the argument of Freed--Hopkins--Teleman~\cite[\S3.6]{FHTI}. To explain \emph{universality} in this context, observe that a Real map $V\colon X\to \pt\sq \O_{p,q}$ classifies a Real vector bundle $V_\C\to X$, where we extend the involution on $V$ to a $\C$-antilinear involution on the complexification~$V_\C$ (see Lemma~\ref{lem:RealVB}). Let $\tau_V$ denote the pullback of the Thom twisting $\SW_{p,q}$ along the classifying map for~$V$. The induced homomorphism on twisted $\KR$-groups
\beq\label{eq:univThom2}
&&\KR^{\SW_{p,q}}(\R^{p,q}\sq \O_{p,q})\to \KR^{\SW_V }(V),\qquad \Th_{p,q}\mapsto \Th_V,
\eeq
sends the universal class $\Th_{p,q}$ to a class~$\Th_V\in \KR^{\tau_V}(V)$ that implements a twisted Thom isomorphism for $V$ by pulling back the universal Thom isomorphism (Corollary~\ref{cor:Thomiso}). Forgetting Real structures, $\Th_V$ is the twisted Thom class of Freed--Hopkins--Teleman~\cite[\S3.6]{FHTI} and~\cite[1.92]{Vienna} and Donovan--Karoubi~\cite{DonovanKaroubi}. A Real Spin$^c$-structure $s$ on $V$ determines an untwisting,
$$
\KR^{\SW_V }(V) \stackrel{s}{\simeq} \KR^{p,q}(V)
$$
that sends $\Th_V$ to the untwisted Thom classes constructed by Stolz--Teichner~\cite[Remark~3.2.22]{ST04} and Atiyah--Bott--Shapiro~\cite{ABS}, see Proposition~\ref{prop:standardThom}. Finally, as any Thom class must restrict locally to the suspension class, Remark~\ref{rmk:suspension} demonstrates the~\eqref{eq:univThom} is the most general incarnation of the twisted Thom class.

\subsection{Twisted equivariant pushforwards and Atiyah--Bott--Shaprio maps}
The universal twisted Thom classes~\eqref{eq:Thomtwists} determine pushforwards along embeddings $f\colon \X\hookrightarrow \Y$ of proper smooth stacks following a standard recipe: compose a Thom isomorphism with the Pontryagin--Thom map (see Definition~\ref{defn:pushembedding})
\beq\label{eq:pushforwardalongembed0}
&&\KR^{f^*\sigma\otimes\tau_{\nu}^{-1}}(\X)\xrightarrow{\Th_{\nu}} \KR^{\sigma}(\nu)\xrightarrow{{\rm Pontryagin-Thom}} \KR^{\sigma}(\Y),\quad \sigma\in \Twist(\Y)
\eeq
for the Thom class $\Th_\nu$ and Thom twisting $\tau_\nu\in \Twist(\X)$ of the normal bundle $\nu$ to the embedding. 
With~\eqref{eq:pushforwardalongembed0} in hand, the pushforward along a $G$-equivariant submersion $\pi\colon X\to Y$ again follows the standard recipe: choose a Real $G$-equivariant embedding $X\hookrightarrow \R^{p,q}$ and compose a Thom isomorphism with the inverse to a (twisted equivariant) suspension isomorphism (see Definition~\ref{defn:ABS})
\beq\label{eq:intropushforward}
&&\pi_!^{\rm top}\colon \KR^{\pi^*\sigma\otimes \tau_\nu^{-1}}(X\sq G)\xrightarrow{\pi_!^{\rm top}}\KR^{\sigma}((\R^{p,q}\times Y)\sq G)
\xleftarrow[\smallsmile \Th_\rho]{\simeq}\KR^{\sigma\otimes \tau_\rho^{-1}}(Y\sq G)
\eeq
for $\nu$ the normal bundle to the embedding $X\sq G\hookrightarrow (\R^{p,q}\times Y)\sq G$ and $\tau_\rho$ the restriction of the Thom twisting along the $G$-representation $G\to \O_{p,q}$. 
This sketches the twisted equivariant Atiyah--Bott--Shapiro orientation, compare~\cite[\S4.1]{HebestreitJoachim},~\cite[\S9.2.2]{FreedHopkins}. Theorems~\ref{thm:3} and~\ref{thm:ABS} then follow quickly from compatibility of power operations with the Thom isomorphism (Proposition~\ref{prop:ThomisocommuteswithPow}), which itself is an observation using Theorem~\ref{thm:A}. 

\begin{rmk}
In string theory, D-brane charges are classified by K-theory, where the relevant flavor of K-theory depends on the type of string theory under consideration~\cite{Minasian,WittenDbrane,KapustinDbrane}. One variant takes an oriented submanifold $f\colon Q\hookrightarrow X$ (i.e., a D-brane) satisfying the anomaly cancellation condition $-W_3(TQ)=[f^*\tau]$ for $\tau\in \Twist(X)$ and a complex vector bundle $V\to Q$ with corresponding D-brane charge $f_!(V)\in \K^\tau(X)$ via the $\KU$-version of~\eqref{eq:pushforwardalongembed0}. Orientifold string theories enhance this construction to unoriented $Q$ and KR-valued D-brane charges~\cite[\S5.2]{WittenDbrane} \cite{GukovDbrane}. The construction~\eqref{eq:pushforwardalongembed0} offers a geometric cocycle counterpart to existing analytical tools, e.g.,~\cite{Bouwknegt,MMS,BMRS,Carey}
\end{rmk}



\subsection{Overview and conventions}

All tensor products below are $\Z/2$-graded. We use the notation $\ell^2$ for the $\Z/2$-graded version of the usual sequence space, so that $\ell^2$ is a $\Z/2$-graded complex Hilbert space with countably infinite-dimensional even and odd subspaces that are orthogonal with repsect to the inner product. We view topological stacks as a localization of topological groupoids at the essential equivalences, following~\cite{Pronk}. As such, we use the same notation for a topological groupoid and its underlying topological stack. When there is potential ambiguity, we indicate in the text when an arrow between topological groupoids is a continuous functor or a map of stacks, i.e., a zigzag~\eqref{eq:zigzag}. Unless otherwise stated, all topological stacks are assumed to be \emph{local quotient stacks}, see~\S\ref{sec:locquostacks}. 

\subsection{Acknowledgements} 
It is a pleasure to thank Tobi Barthel, Millie Deaton, Zach Halladay,  Andr\'e Henriques, Mike Hopkins, Yigal Kamel, Kiran Luecke, Charles Rezk, Nat Stapleton, Stephan Stolz, and Peter Teichner for enlightening conversations that helped to shape this work. We also thank Arun Debray, Yigal Kamel, Cameron Krulewski, Natalia Pacheco--Tallaj, and Luuk Stehouwer for feedback on an earlier draft. 
DBE is partially supported by the NSF through grants DMS-2205835 and DMS-2340239.

\section{The projective unitary group of a Real Clifford module}

\subsection{The category of Real graded  groups}\label{sec:Realgroups}

\begin{defn}[{\cite[\S1]{Atiyahreal}}] The category of \emph{Real spaces} has objects spaces $X$ with involution $(\overline{\phantom{g}})\colon X\to X$ and morphisms continuous maps $f\colon X\to Y$ compatible with the involutions, $f(\overline{x})=\overline{f(x)}$. 
\end{defn}

\begin{ex}\label{ex:RealCn}
The most basic examples of Real spaces are~$\C^n$ with involution given by the usual complex conjugation $z\mapsto \overline{z}$. 
\end{ex}

\begin{ex}\label{ex:Rpq}
Let $\R^{p,q}$ denote the Real space given by $\R^{p+q}=\R^p\oplus \R^q$ with  involution
\beq\label{eq:RpqRealstructure}
(x,y)\mapsto (x,-y), \quad x\in \R^p,\ y\in \R^q. 
\eeq
We note the isomorphisms of Real spaces $\R^{n,n}\simeq \C^n$ given by $(x,y)\mapsto x+iy=z$. 
\end{ex}

\begin{defn}[{\cite[\S5]{Atiyahbottp}}]
A \emph{Real group} is a group object in Real spaces, i.e., a group $G$ with an involution homomorphism $(\overline{\phantom{g}})\colon G\to G$. A \emph{Real homomorphism} $(G,(\overline{\phantom{g}}))\to (H,(\overline{\phantom{g}}))$ is a homomorphism $\varphi\colon G\to H$ compatible with the involutions, $\varphi(\overline{g})=\overline{\varphi(g)}$. A \emph{Real action} of a Real group $G$ on a Real space $X$ is an action compatible with the Real structure, $\overline{g\cdot x}=\overline{g}\cdot \overline{x}$ for $g\in G$ and $x\in X$. 
\end{defn}

\begin{ex}\label{ex:UnReal}
The Real structure on $\C^n$ promotes the unitary group $\U_n$ to a Real group with involution corresponding to complex conjugation of matrix entries. Equivalently,
\beq\label{eq:RealUn}
\overline{u}=(\overline{\phantom{g}})\circ u \circ (\overline{\phantom{g}}),\qquad u\in \U_n
\eeq
for pre- and post-composition by the involution on $\C^n$ from Example~\ref{ex:RealCn}. The standard~$\U_n$-action on $\C^n$ determines a Real action.
\end{ex}

\begin{ex}\label{ex:OpqReal}
Let $\O_{p,q}$ denote the Real group defined by the orthogonal group $\O_{p+q}$ with Real structure given by the same formula as~\eqref{eq:RealUn} but for conjugation by the involution~\eqref{eq:RpqRealstructure} on $\R^{p,q}$. The standard $\O_{p+q}$-action on $\R^{p+q}$ determines a Real action of $\O_{p,q}$ on $\R^{p,q}$. 
\end{ex}

\begin{ex}\label{ex:CnUn}
A Real action of a Real group $G$ on a Real space $X$ determines a Real stack $X\sq G$ in the sense of Definition~\ref{defn:Realstack}. In particular, the above examples give the Real stacks $\R^{p,q}\sq \O_{p,q}$ and $\C^n\sq \U_n$.  The standard isometric identification $\C^n\simeq \R^{2n}$ and associated injection $\U_n\hookrightarrow \O_{2n}$ determines a map of of Real stacks
$$
\C^n\sq \U_n\to \R^{n,n}\sq \O_{n,n}.
$$
\end{ex}

\begin{defn}\label{defn:gradedgroup}
A \emph{graded group} $(G,\epsilon)$ is a group $G$ with a homomorphism $\epsilon \colon G\to \Z/2$. A \emph{graded homomorphism} $(G,\epsilon)\to (H,\delta)$ is homomorphism $\varphi\colon G\to H$ that intertwines the gradings, i.e., $\delta \circ \varphi=\epsilon$. 
\end{defn}

\begin{defn}
A \emph{Real graded  group} $(G,\epsilon,(\overline{\phantom{g}}))$ is a group $G$ with a grading and Real structure satisfying the compatibility condition $\epsilon(\overline{g})=\overline{\epsilon(g)}.$ Real graded  groups are the objects of a category whose morphisms are Real, graded homomorphisms. 
\end{defn}

\begin{ex}\label{eq:Opqgraded}
The determinant homomorphism
$$
\epsilon=\det\colon \O_{p,q}\to \{\pm 1\}\simeq \Z/2
$$
endows $\O_{p,q}$ with the structure of a Real graded group. 
\end{ex}

\begin{rmk}
The involution $(\overline{\phantom{g}})\colon G\to G$ generates an action on $G$ by a group of order~2. This group plays a distinct role from the $\Z/2$ involved in the grading. To disambiguate these notationally, throughout the paper $C_2$ will denote the group whose action on $G$ is generated by the involution $(\overline{\phantom{g}})$, whereas $\Z/2$ is the target of the grading homomorphism~$\epsilon\colon G\to \Z/2$. 
\end{rmk}

\begin{rmk}
There are various ways to pass between combinations of Real and graded groups. For example, any group $G$ has the trivial grading $G\to \{e\}<\Z/2$ and the trivial Real structure given by the identity involution. As another example, for a Real graded  group $(G,\epsilon,(\overline{\phantom{g}}))$ the fixed point subgroup $G^{C_2}<G$ inherits a grading from the restriction of $\epsilon$, giving a functor from Real graded  groups to graded groups. Applied to Examples~\ref{ex:UnReal} and~\ref{ex:OpqReal},
$$
\O_p\times \O_q=\O_{p,q}^{C_2}\subset \O_{p,q} \qquad \O_n=\U_n^{C_2}\subset \U_n.
$$
Similarly the kernel ${\rm ker}(\epsilon)<G$ of the grading homomorphism gives a functor from Real graded  groups to Real groups. Applied to Example~\ref{eq:Opqgraded}, we get $\SO_{p,q}={\rm ker}(\det)\subset \O_{p,q}$. 
\end{rmk}

\begin{rmk}\label{rmk:C2quotient}
For a Real group, the $C_2$-action on $G$ determines the semidirect product $\tilde{G}=G\rtimes C_2$, and the kernel of the homomorphism $\phi\colon \tilde{G}\to C_2$ recovers $G$. This leads to an equivalence of categories between Real, graded groups in the sense above, and groups with a pair of homomorphisms $\epsilon \colon \tilde{G}\to \Z/2$ and $\phi\colon \tilde{G}\to C_2$ satisfying a compatibility condition; this is the perspective of \emph{$\phi$-twisted graded central extensions} taken in~\cite[\S7]{FreedMoore} and \cite{Gomi}. When considering (twisted) power operations, it turns out to be slightly more convenient to follow Atiyah, taking a Real structure as an involution~\cite[\S5]{Atiyahbottp}. We compare these points of view in more detail in~\S\ref{sec:Realstacks}. 
\end{rmk}

\subsection{The unitary group of a Real graded  Hilbert space as a Real graded  group}

\begin{defn}
A \emph{$\Z/2$-graded Hilbert space} is a graded inner product space $\H=\H_+\oplus \H_-$ over $\C$ with the property that the even and odd subspaces~$\H_+$ and~$\H_-$ are orthogonal and there exists and isomorphism $\H_+\simeq \H_-$. We allow $\H$ to be either finite or (countably) infinite dimensional. 
\end{defn}
\begin{rmk}
We emphasize that the isomorphism $\H_+\simeq \H_-$ is property, not data. This property guarantees the existence of odd isometries, and hence that the unitary groups below admit a nontrivial grading in the sense of Definition~\ref{defn:gradedgroup}. 
\end{rmk}

\begin{rmk}
Often we will use that a $\Z/2$-graded Hilbert space $\H$ is equivalent to a Hilbert space equipped with an involution~$(-1)^F$ satisfying properties. One recovers the grading from eigenspaces, $(-1)^F|_{\H_\pm}=\pm \id_{\H_\pm}$. From this perspective an even (respectively, odd) endomorphism $\H\to \H$ commutes (respectively, anticommutes) with the involution~$(-1)^F$. 
\end{rmk}

\begin{defn}A \emph{Real structure} on a $\Z/2$-graded Hilbert space is a grading preserving (i.e., even) map $r\colon \cH\to \overline{\cH}$ satisfying $\overline{r}\circ r=\id_\cH$. 
\end{defn}

Exercising care with some signs, the category of $\Z/2$-graded Hilbert spaces and isometries has a symmetric monoidal structure given by the $\Z/2$-graded Hilbert tensor product, e.g., see~\cite[pages~89-91]{strings1}. 

\begin{defn}
The \emph{unitary group} of a Real $\Z/2$-Real graded  Hilbert space is 
\beq
 \U^\pm(\H)&:=&\{u\in \U(\H)\mid (-1)^F u(-1)^F=+ u \ {\rm or} \ (-1)^F u(-1)^F=-u \}, \nonumber\\
 &=& \left\{u\in \U(\H)\mid u=\left[\begin{array}{cc} * & 0 \\ 0 & *\end{array}\right] \ {\rm or} \ u=\left[\begin{array}{cc} 0 & * \\ * & 0\end{array}\right]\right\}\label{eq:gradedunitary}
\eeq
where $\U(\H)$ is the unitary group of the Hilbert space $\cH$ (ignoring the $\Z/2$-grading), and the block diagonal form in~\eqref{eq:gradedunitary} is relative to the decomposition~$\H=\H^+\oplus \H^-$ into even and odd subspaces. The \emph{grading homomorphism} $\epsilon \colon  \U^\pm(\H)\to \Z/2$ records whether an element of $\U^\pm(\cH)$ preserves or reverses the grading on~$\cH$. The \emph{Real structure} on $\U^\pm(\cH)$ is $u\mapsto \overline{r}\circ u\circ r$ for the Real structure $r\colon \cH\to \overline{\cH}$ on~$\cH$. 
\end{defn}

The $C_2$-fixed points $\U^\pm(\cH)^{C_2}\subset \U^\pm(\cH)$ for the Real structure are the $\R$-linear unitary maps $\cH\to \cH$, and hence restrict to maps on the real $\Z/2$-graded Hilbert space $\cH_\R=\cH^{C_2}\subset \cH$. This subgroup can therefore be identified with the graded projective orthogonal group $\O^\pm(\cH_\R)\simeq \U^\pm(\cH)^{C_2}$, defined in analogy to~\eqref{eq:gradedunitary} but for the real Hilbert space $\cH_\R$. 

\begin{defn}
The quotient by $\U(1)$ scalars defines the \emph{projective unitary group}
\beq\label{eq:projectiveexactseq}
1\to \U(1)\to \U^\pm(\cH)\to \PU^\pm(\cH)\to 1.
\eeq
The group $\PU^\pm(\cH)$ inherits a Real graded  structure from $\U^\pm(\cH)$ characterized by~\eqref{eq:projectiveexactseq} being an exact sequence in Real graded groups where $\U(1)$ has the trivial grading and Real structure from complex conjugation $\U(1)\subset \C$. 
\end{defn}



Functoriality gives a map of exact sequences whose vertical arrows are inclusions of $C_2$-fixed points
\beq\label{eq:POH}
\begin{array}{ccccccccc} 1 &\to& \{\pm 1\}&  \to& \O^\pm(\cH_\R)&\to&  \PO^\pm(\cH_\R)&\to& 1\\
&&\downarrow && \downarrow &&\downarrow \\
1&\to& \U(1)&  \to& \U^\pm(\cH)&\to &  \PU^\pm(\cH)&\to& 1, \end{array}
\eeq
and so we obtain an analogous graded projective orthogonal group.

When $\cH$ is finite-dimensional, $\U^\pm(\cH)$ is a Real graded  Lie group, and in particular a topological group. When $\cH$ is infinite-dimensional, we endow $ \U^\pm(\H)$ with the strong operator topology, which is equivalent~\cite{EspinozaUribe,Schottenloher} to the weak operator topology or the compact open topology used in~\cite{AtiyahSegaltwists}.\footnote{The continuous map from $\U^\pm(\H)$ in the (more common) norm topology to $\U^\pm(\H)$ in the strong topology is not a homeomorphism~\cite{EspinozaUribe,Schottenloher}.}  This topology lifts~\eqref{eq:POH} to a diagram in topological groups, and the long exact sequence for the fibrations computes the nonzero homotopy groups
\beq\label{eq:homotopygrps}
\begin{array}{ll}  \pi_0\PU^\pm(\H)\simeq \Z/2,& \pi_2\PU^\pm(\H)\simeq \Z\\
\pi_0\PO^\pm(\cH_\R)\simeq \Z/2,& \pi_1\PO^\pm(\cH_\R)\simeq \Z/2\end{array}
\eeq
using contractibility of the connected component of the identity element $\U^+(\cH)=\U(\cH^+)\times \U(\cH^-)<\U^\pm(\cH)$ and $\O^+(\cH)=\O(\cH^+)\times \O(\cH^-)<\O(\cH)$~\cite[Proposition A2.1]{AtiyahSegaltwists}.

\subsection{Real Clifford algebras and the Real, graded $\Pin^c$-groups}\label{sec:RealCliff} The following is standard e.g.,~\cite[\S4]{Atiyahreal} or~\cite[\S{I.10}]{LM}, though there are differences in signs inherited from choices of convention. We establish our conventions below. 

For $\R^{p,q}$ as in Example~\ref{ex:Rpq}, define the Real Clifford algebra $\cCl_{p,q}=\cCl(\R^{p,q})$ to be the Clifford algebra (over $\C$) of $\R^{p+q}$ with antilinear algebra involution $\cCl_{p,q}\to \cCl_{p,q}$ determined by the involution of $\R^{p,q}$. For $a\in \cCl_{p,q}$, let $a\mapsto \overline{a}$ denote this involution generating a $C_2$-action on $\cCl_{p,q}$. In terms of generators and relations, 
\beq\label{eq:Cliffdefn}
&&\cCl_{p,q}=\langle e_1,\dots, e_{p+q} \mid e_je_k+e_ke_j=2\delta_{jk}\rangle,\qquad \overline{e_j}=\left\{\begin{array}{cc} e_j & j\le p\\ -e_j & j >p. \end{array}\right.
\eeq
The (real) fixed point subalgebra $\Cl_{p,q}=\cCl_{p,q}^{C_2}\subset \cCl_{p,q}$ is generated by
\beq\label{eq:realsubalgebra}
&&\{e_1,\dots, e_p, ie_{p+1},\dots , ie_{p+q}\} \subset \cCl_{p,q},\quad e_j^2=+1, \ j\le p, \ \ f_j^2=-1, \ f_j:=ie_j, \ j>p
\eeq
and hence is identified with real Clifford algebra $\Cl_{p,q}\subset \cCl_{p,q}$ associated to the quadratic form of mixed signature $(p,q)$ on $\R^{p+q}$. The graded tensor product gives an isomorphism of Real Clifford algebras,
\beq\label{eq:tensorofCliffordalgebras}
\cCl_{p,q}\otimes \cCl_{p',q'}\simeq \cCl_{p+p',q+q'}. 
\eeq

The \emph{Pin$^c$ groups}, $\Pin^c_{p,q}$, are defined as the group generated by vectors of unit norm $\R^{p,q}\subset \cCl_{p,q}$ together with unit norm scalars $\U(1)\subset \C\subset \cCl_{p,q}$. The grading and Real structure on $\cCl_{p,q}$ restrict to $\Pin^c_{p,q}$, yielding a Real graded  central extension 
\beq\label{eq:pinpinc}
1\to \U(1)  \to \Pin^c_{p,q}\xrightarrow{{\rm Ad}}   \O_{p,q}\to 1,\qquad   \O_{p,q}\xrightarrow{\det} \Z/2
\eeq
for the trivial Real structure on $\O_{p,q}$, the Real structure on $\U(1)$ from complex conjugation, and the grading on $\O_{p,q}$ given by the determinant. The \emph{Pin groups} are defined as the graded subgroups $\Pin_{p,q}:=(\Pin_{p,q}^c)^{C_2}\subset \Pin_{p,q}^c$ of $C_2$-fixed points. The \emph{Spin$^c$ groups} are the Real subgroups $\Spin^c_{p,q}:=\ker(\det\circ {\rm Ad})\subset \Pin^c_{p,q}$. Finally, the \emph{Spin groups} are the further $C_2$-fixed points, $\Spin_{p,q}\subset (\Spin^c_{p,q})^{C_2}\subset \Spin^c_{p,q}$

\begin{rmk}
The homomorphism ${\rm Ad}$ in~\eqref{eq:pinpinc} is the \emph{twisted adjoint representation} that sends a unit norm generator $x\in \R^{p,q}\subset \cCl_{p,q}$ to the reflection in $\O_{p,q}$ through the plane orthogonal to~$x$, e.g., see~\cite[page 14]{LM}. Hence the grading induced from the Clifford algebra $\Pin^c_{p,q}\subset \cCl_{p,q}$ is compatible with the grading from the sign of the determinant, 
\beq\label{eq:Pingrading}
(-1)^{|\omega|}=\det({\rm Ad}(\omega)),\qquad \omega\in \Pin^c_{p,q}\subset \cCl_{p,q}, \ |\omega|\in \{0,1\}.
\eeq
\end{rmk}

\subsection{Real Clifford modules and $\Pin^c$-representations}\label{sec:spinor}
A \emph{Real} right (respectively, left) $\cCl_{p,q}$-module is a a complex $\Z/2$-graded $\cCl_{p,q}$-module $\cH=\cH^0\oplus \cH^1$ with an even antilinear involution $v\mapsto \overline{v}$ satisfying 
$$
\overline{v\cdot a}=\overline{v}\cdot \overline{a}, \ \ ({\rm respectively,} \ \ \overline{a\cdot v}=\overline{a}\cdot \overline{v}),\qquad a\in \cCl_{p,q}. 
$$
Right modules and left module can be exchanged by taking opposite algebras. For example, left modules for $\cCl_{p,q}$ are right modules over $\cCl_{p,q}^\op$, the opposite Real algebra with multiplication
\beq\label{eq:superopposite}
a\cdot_{\op} b=(-1)^{|a||b|}b\cdot a, \qquad a,b\in \cCl_{p,q}
\eeq
and the Real structure is determined by the same $C_2$-action on generators $\R^{p,q}\subset \cCl_{p,q}^\op$. 

A right (respectively, left) Clifford module is \emph{unitary} if the endomorphisms associated generators $e\in \R^{p+q}\subset \cCl_{p,q}$ are skew-adjoint (respectively, self-adjoint). 
The following lemma justifies this definition. 

\begin{lem}\label{lem:pinorrep}
    A left (or right) Real $\cCl_{p,q}$-module $V$ canonically determines a Real, graded, unitary $\Pin^c_{p,q}$-representation on $V$. 
    In particular, any left Clifford module~$V$ uniquely determines a Real, graded projective representation
\beq\label{eq:projrepOpq}
\O_{p,q}\to \PU^\pm(V).
\eeq
\end{lem} 
\bp
Given a left $\cCl_{p,q}$-module~$V$, a Real, graded $\Pin^c_{p,q}$-representation comes from the composition,
\beq\label{eq:spinorrep}
\Pin^c_{p,q}\subset \cCl_{p,q}\xrightarrow{\gamma} \End(V), 
\eeq
where the algebra map $\gamma$ is the module structure. 
To see that the above factors through the graded unitary group, $\Pin^c_{p,q}\to \U^\pm(V)\subset \End(V)$, we compute the action of the generators,
\beq\label{eq:unitaryrep}
&&\langle \gamma(e)\cdot v,\gamma(e)\cdot w\rangle=\langle v,(\gamma(e)^*\gamma(e))\cdot w\rangle =\langle v,\gamma(e^2)\cdot w\rangle=\langle v,w\rangle, \quad e\in \R^{p,q}\subset\cCl_{p,q}
\eeq
using that $\gamma(e)^*=\gamma(e)$ and $e^2=1$. 
Since $\U(1) \subset \Pin^n_{p,q}$ acts by scalar multiplication in~\eqref{eq:spinorrep}, the representation~\eqref{eq:spinorrep} is equivalent to a projective representation of $\O_{p,q}$ for which the $\U(1)$-extension of $\O_{p,q}$ associated to~\eqref{eq:projrepOpq} is the $\Pin^c$-extension of $\O_{p,q}$; see~\eqref{eq:pinpinc} below. 

The argument for right Clifford modules is similar: identify a right Clifford module $V$ with a left module for the opposite algebra. Then the calculation~\eqref{eq:unitaryrep} carries an additional sign from the opposite multiplication~\eqref{eq:superopposite}, and the skew-adjointness of the generators implies that the associated $\Pin^c$-representation is unitary. In this representation, the subgroup $\U(1)\subset\Pin^c_{p,q}$ acts with weight $-1$ on $V$, and hence one obtains a projective representation of $\O_{p,q}$ for the \emph{dual} of the Pin$^c$ extension. 
\ep

The direct sum and graded tensor product of Real unitary Clifford modules is again a Real unitary Clifford module \cite[Remark~I.10.10]{LM}, e.g., for right modules
$$
(V_{\cCl_{p,q}})\otimes (V'_{\cCl_{p',q'}} )\simeq (V\otimes V')_{\cCl_{p+p',q+q'}},\qquad  V_{\cCl_{p,q}}\oplus  V'_{\cCl_{p,q}}\simeq (V\oplus V')_{\cCl_{p,q}} .
$$

\begin{rmk}
Taking fixed points for the involution on a Real Clifford module $V_\R\subset V$ produces a module over the real subalgebra $\Cl_{p,q}$ from~\eqref{eq:realsubalgebra}. If the Real left Clifford module~$V$ is unitary, the real $\Cl_{p,q}$-module~$V_\R$ has $e_j$ acting by self-adjoint maps for $j\le p$ and skew-adjoint maps for $f_j=(ie_j)$ with $j>p$. 
\end{rmk}




\subsection{The projective unitary group of a Real Clifford module}

\begin{defn}\label{defn:UCl}
Let $\cH$ be a Real, unitary right $\cCl_{p,q}$-module. Define the subgroup
\beq\label{eq:restrictedunitaryCl}
 &&\U_{\cCl_{p,q}}^\pm(\H):=\{u\in \U^\pm (\H)\mid u(v\cdot e) = u(v)\cdot e, \quad \forall e\in \R^{p+q}\subset \cCl_{p,q},\  v \in \H \}, 
\eeq
of \emph{Clifford linear unitary maps} for $\U^\pm(\H)$ the Real graded  group~\eqref{eq:gradedunitary}. The grading and Real structure on $\U^\pm(\cH)$ restrict to the subgroup $\U_{\cCl_{p,q}}^\pm(\cH)<\U^\pm(\cH)$, which is a Real graded  topological group for the subspace topology. 
\end{defn}



The quotient by scalars $\U(1)< \U_{\cCl_{p,q}}^\pm(\H)$ gives the exact sequence of Real graded  topological groups
\beq\label{eq:PUHgradedgerbe}
&&1\to \U(1)\to \U_{\cCl_{p,q}}^\pm(\H)\to \PU_{\cCl_{p,q}}^\pm(\H)\to 1
\eeq
where $\U(1)$ has Real structure from complex conjugation and is trivially graded. The Real structure on $\PU_{\cCl_{p,q}}^\pm(\H)$ comes from the Real structure on $\cH$, and grading $\epsilon=\pi_0\colon \PU_{\cCl_{p,q}}^\pm(\H)\to \Z/2$ is determined by the grading on $\U^\pm_{\cCl_{p,q}}(\cH)$ (using continuity and the fact that $\U(1)$ is connected).


Next consider the diagrams (compare \cite[page 10]{AtiyahSegaltwists}) of Real graded  groups
\beq\label{eq:tensoronU}
\begin{tikzpicture}[baseline=(basepoint)];
\node (A) at (0,0) {$\U(1)\times \U(1)$};
\node (B) at (5,0) {$\U^\pm_{\cCl_{p,q}}(\H)\times \U^\pm_{\cCl_{p',q'}}(\H')$};
\node (C) at (0,-1.25) {$\U(1)$};
\node (D) at (5,-1.25) {$\U^\pm_{\cCl_{p,q}\otimes \cCl_{p',q'}}(\H\otimes \H').$};
\draw[->] (A) to  (B);
\draw[->] (B) to node [right] {$\otimes$} (D);
\draw[->] (A) to node [left] {$\times$} (C);
\draw[->] (C) to (D);
\path (0,-.75) coordinate (basepoint);
\end{tikzpicture}
\eeq
The maps $\otimes$ are induced by the graded tensor product of $\Z/2$-graded Hilbert spaces; the diagram commutes for left vertical arrow given by multiplication in $\U(1)$. Hence, the diagrams~\eqref{eq:tensoronU} determine homomorphisms
\beq\label{eq:tensoronPU1}
&&\PU_{\cCl_{p,q}}^\pm(\H)\times \PU_{\cCl_{p',q'}}^\pm(\H')\to \PU_{\cCl_{p,q}\otimes \cCl_{p',q'}}^\pm(\H\otimes \H').
\eeq
Removing one copy of the scalars in~\eqref{eq:tensoronU} gives the homomorphisms
\beq\label{eq:UactiononPU}
&&\U^\pm_{\cCl_{p,q}}(\H)\times \PU^\pm_{\cCl_{p',q'}}(\H')\to \PU^\pm_{\cCl_{p,q}\otimes \cCl_{p',q'}}(\H\otimes \H')
\eeq
that tensor an isometry and projective transformation to get a projective transformation. An isometric embedding $\H\subset \H'$ of unitary Real ${\cCl_{p,q}}$-modules induces homomorphisms
\beq\label{eq:embedding1}
\PU_{\cCl_{p,q}}^\pm(\H)\hookrightarrow \PU_{\cCl_{p,q}}^\pm(\H').
\eeq


\begin{rmk}\label{rmk:RealMorita}
There are isomorphisms of topological groups
$$
\PU_{\cCl_{p+1,q+1}}^\pm(\H\otimes_\C M_{1,1})\simeq \PU_{\cCl_{p,q}}^\pm(\H),
$$
$$
\PU_{\cCl_{p,q+8}}^\pm(\H\otimes_\R M_{0,8})\simeq \PU_{\cCl_{p,q}}^\pm(\H),
$$
$$
\PU_{\cCl_{p+8,q}}^\pm(\H\otimes_\R M_{8,0})\simeq \PU_{\cCl_{p,q}}^\pm(\H),
$$
where $M_{1,1}\simeq \C^{1|1}$ is the Real irreducible graded $\C$-$\cCl_{1,1}$ bimodule, i.e., the bimodule implementing the $(1,1)$-periodicity in the Atiyah--Bott--Shapiro construction~\cite[I.10]{LM}. 
Similarly $M_{8,0}\simeq \R^{8|8}$ is the irreducible $\R$-$\Cl_8$ bimodule, and $M_{0,8}$ is the irreducible $\R$-$\Cl_{-8}$-bimodule. 
\end{rmk} 

\begin{rmk}
Analogous to~\eqref{eq:POH}, $C_2$-fixed subgroups determine the orthogonal groups
$$
\O_{\Cl_{p,q}}^\pm(\cH_\R)=\U_{\cCl_{p,q}}^\pm(\cH)^{C_2}<\U_{\cCl_{p,q}}^\pm(\cH).
$$
By functoriality,~\eqref{eq:PUHgradedgerbe} gives a map of exact sequences whose vertical arrows are inclusions of $C_2$-fixed points
\beq\nonumber
\begin{array}{ccccccccc} 1 &\to& \{\pm 1\}&  \to& \O^\pm_{\Cl_{p,q}}(\cH_\R)&\to&  \PO^\pm_{\Cl_{p,q}}(\cH_\R)&\to& 1\\
&&\downarrow && \downarrow &&\downarrow \\
1&\to& \U(1)&  \to& \U^\pm_{\cCl_{p,q}}(\cH)&\to &  \PU^\pm_{\cCl_{p,q}}(\cH)&\to& 1, \end{array}
\eeq
and the homomorphisms~\eqref{eq:tensoronPU1},~\eqref{eq:UactiononPU} and~\eqref{eq:embedding1} restrict to the fixed point subgroups. 
\end{rmk}

\section{Projective Hilbert bundles and Real graded  gerbes}\label{sec:ProjHilbgerbes}

The central extension of topological groups
\beq\label{eq:itsaline}
\U(1)\to \U^\pm_{\cCl_{p,q}}(\cH)\to \PU^\pm_{\cCl_{p,q}}(\cH)
\eeq 
can equivalently be viewed as a multiplicative principal $\U(1)$-bundle: multiplication on~$\U^\pm_{\cCl_{p,q}}(\cH)$ provides an isomorphism of $\U(1)$-bundles over $\PU^\pm_{\cCl_{p,q}}(\cH)$, 
\beq\label{eq:multiplicativeline}
&&\begin{array}{c} p_1^*\U^\pm_{\cCl_{p,q}}(\cH)\otimes p_2^*\U^\pm_{\cCl_{p,q}}(\cH)\simeq m^*\U^\pm_{\cCl_{p,q}}(\cH), \\ 
p_1,p_2,m\colon \PU^\pm_{\cCl_{p,q}}(\cH)\times \PU^\pm_{\cCl_{p,q}}(\cH) \to \PU^\pm_{\cCl_{p,q}}(\cH)\end{array}
\eeq
for pullbacks along the projections and multiplication; the isomorphism~\eqref{eq:multiplicativeline} of $\U(1)$-bundles furthermore satisfies an associativity condition on $\PU^\pm_{\cCl_{p,q}}(\cH)^{\times 3}$. 

The Real graded  group $\U^\pm_{\cCl_{p,q}}(\cH)$ endows the multiplicative principal bundle~\eqref{eq:itsaline} with graded and Real structures. To spell this out, the \emph{grading} regards the $\U(1)$-bundle over the grading-preserving component $\U^+_{\cCl_{p,q}}(\cH)\to \PU^+_{\cCl_{p,q}}(\cH)$ as even whereas the grading-reversing component $\U^-_{\cCl_{p,q}}(\cH)\to \PU^-_{\cCl_{p,q}}(\cH)$ is an odd bundle. The multiplication on $\U^\pm_{\cCl_{p,q}}(\cH)$ is compatible with the graded tensor product of these $\U(1)$-bundles. Similarly, the \emph{Real structure} on $\U^\pm_{\cCl_{p,q}}(\cH)\to \PU^\pm_{\cCl_{p,q}}(\cH)$ is the $C_2$-action on $\U^\pm_{\cCl_{p,q}}(\cH)$ and $\PU^\pm_{\cCl_{p,q}}(\cH)$ determined by the Real structure on~$\cH$. Hence,~\eqref{eq:itsaline} is a Real bundle in the sense of Atiyah~\cite{Atiyahreal}, where the $C_2$-action on the $\U(1)$-fibers is complex conjugation. As $\U^\pm_{\cCl_{p,q}}(\cH)$ is a Real topological group,~\eqref{eq:multiplicativeline} is an isomorphism of Real bundles. 

In this section, we analyze these structures for the classifying stacks associated to~\eqref{eq:itsaline}. 


\subsection{Projective unitary groups and Real graded  gerbes}

We refer to~\S\ref{sec:appenRealgerbe} for a review of Real graded gerbes. 

\begin{lem}\label{lem:gradedRealgerbe}
The central extension~\eqref{eq:itsaline} determines a Real graded  gerbe over~$\pt\sq \PU^\pm_{\cCl_{p,q}}(\cH)$. 
\end{lem}
\bp The $\U(1)$-central extension of groups~\eqref{eq:itsaline} determines a $\pt\sq \U(1)$-extension of groupoids, and hence a gerbe over the stack $\pt\sq \PU^\pm_{\cCl_{p,q}}(\cH)$, (see Lemma~\ref{lem:BX})
\beq\label{eq:fibrationmain}
\pt\sq \U(1) \to \pt\sq\U^\pm_{\cCl_{p,q}}(\cH)\to \pt\sq \PU^\pm_{\cCl_{p,q}}(\cH).
\eeq 
The grading on $\PU^\pm_{\cCl_{p,q}}(\cH)$ determines a map of topological groupoids $\epsilon \colon \pt\sq \PU^\pm_{\cCl_{p,q}}(\cH)\to \pt\sq(\Z/2)$, which then promotes~\eqref{eq:fibrationmain} to a graded gerbe over the stack $\pt\sq \PU^\pm_{\cCl_{p,q}}(\cH)$. Finally, the Real structure on $\PU^\pm_{\cCl_{p,q}}(\cH)$ promotes~\eqref{eq:fibrationmain} to a fibration in Real stacks, coming from a fibration of topological groupoids with (strict) $C_2$-action. The $C_2$-action on the $\pt\sq \U(1)$ fibers is complex conjugation on $\U(1)$, and this is compatible with the grading. Hence~\eqref{eq:fibrationmain} has the structure of a Real graded  gerbe over $\pt\sq \PU^\pm_{\cCl_{p,q}}(\cH)$. 
\ep

\begin{cor}\label{cor:gradedRealgerbe}
Let $\X$ be a Real stack. A Real map $\tau\colon \X\to \pt\sq \PU^\pm_{\cCl_{p,q}}(\cH)$ determines a Real graded  gerbe on $\X$. An isomorphism between maps $\tau\Rightarrow \tau'$ determines an isomorphism between Real graded  gerbes on $\X$. 
\end{cor}
\bp
This follows from Lemma~\ref{lem:gradedRealgerbe} by taking the pullback in Real stacks, 
\beq\label{eq:twisttogerbe}
&&
\begin{tikzpicture}[baseline=(basepoint)];
\node (AA) at (0,1) {$\pt\sq \U(1)$};
\node (A) at (0,0) {$\X^\tau$};
\node (BB) at (4,1) {$\pt\sq \U(1)$};
\node (B) at (4,0) {$\pt\sq \U^\pm_{\cCl_{p,q}}(\H)$};
\node (C) at (0,-1) {$\X$};
\node (D) at (4,-1) {$\pt\sq \PU^\pm_{\cCl_{p,q}}(\H)$};
\node (DD) at (7,-1) {$\pt\sq(\Z/2)$.};
\node (P) at (.5,-.3) {\scalebox{1.5}{$\lrcorner$}};
\draw[double equal sign distance] (AA) to  (BB);
\draw[->] (A) to  (B);
\draw[->] (AA) to  (A);
\draw[->] (BB) to  (B);
\draw[->] (B) to  (D);
\draw[->] (A) to  (C);
\draw[->] (C) to node [above] {$\tau$} (D);
\draw[->] (D) to node [above] {$\epsilon$} (DD);
\draw[->,bend right=15] (C) to node [below] {$\epsilon_\tau$} (DD);
\path (0,0) coordinate (basepoint);
\end{tikzpicture}
\eeq
Forgetting the Real and graded structures, $\X^\tau\to \X$ defines a gerbe on $\X$~\cite{TuStacks}, and this gerbe acquires a grading by the indicated map $\epsilon_\tau$. Furthermore, the $C_2$-action on $\X^\tau$ must restrict on each $\pt\sq \U(1)$-fiber to the same action as on $\pt\sq \U^\pm_{\cCl_{p,q}}(\cH)$, i.e., the action coming from the Real structure on $\cH$. This gives the action generated by complex conjugation on $\pt\sq \U(1)$. Furthermore, the map $\epsilon_\tau$ is $C_2$-invariant, as it is the composition of the $C_2$-equivariant map $\tau$ with the $C_2$-invariant map $\epsilon$. This gives $\X^\tau\to \X$ the structure of a Real graded  gerbe. 
\ep

\begin{lem}\label{lem:tensortwist}
The tensor product homomorphism~\eqref{eq:UactiononPU} determines a map of stacks that classifies the tensor product of Real graded  gerbes in Lemma~\ref{lem:gradedRealgerbe}. 
\end{lem}

\bp
We start by spelling out the statement in greater detail. Consider the diagram in Real stacks
\beq\label{eq:tensortwist}
&&
\begin{tikzpicture}[baseline=(basepoint)];
\node (A) at (0,0) {$\pt\sq \U^\pm_{\cCl_{p,q}}(\H)\boxtimes \pt\sq \U^\pm_{\cCl_{p',q'}}(\H')$};
\node (C) at (0,-1.5) {$\pt\sq \PU^\pm_{\cCl_{p,q}}(\H)\times \pt\sq \PU^\pm_{\cCl_{p',q'}}(\H')$};
\node (D) at (5.5,-1.5) {$\pt\sq \PU^\pm_{\cCl_{p,q}\otimes \cCl_{p',q'}}(\H\otimes \H')$};
\node (B) at (5.5,0) {$\pt\sq \U^\pm(\H\otimes \H')$};
\node (E) at (9.5,-1.5) {$\pt\sq(\Z/2).$};
\node (P) at (.7,-.4) {\scalebox{1.5}{$\lrcorner$}};
\draw[->] (A) to  (B);
\draw[->] (B) to  (D);
\draw[->] (A) to  (C);
\draw[->] (D) to node [above] {$\epsilon_{\cH\otimes \cH'}$} (E);
\draw[->] (C) to node [above] {$\otimes$} (D);
\draw[->,bend right=10] (C) to node [below] {$\epsilon_\cH \cdot \epsilon_{\cH'}$} (E);
\path (0,-.75) coordinate (basepoint);
\end{tikzpicture}
\eeq
The arrow labeled by $\otimes$ is determined by the homomorphism~\eqref{eq:UactiononPU}, and by Corollary~\ref{cor:gradedRealgerbe} it classifies a Real graded  gerbe over $\pt\sq \PU^\pm_{\cCl_{p,q}}(\H)\times \pt\sq \PU^\pm_{\cCl_{p',q'}}(\H')$. The statement of the lemma is that this gerbe is isomorphic to the (external) tensor product of Real graded  gerbes in the upper left corner of~\eqref{eq:tensortwist}. As the right vertical arrow is a fibration of topological groupoids, one can compute this pullback in topological groupoids, which in turn comes down to a computation about the topological groups in~\eqref{eq:tensortwist}. Indeed, the tensor product on unitary operators $\U^\pm_{\cCl_{p,q}}(\cH)\times \U^\pm_{\cCl_{p,q}}(\cH')\to \U^\pm_{\cCl_{p,q}}(\cH\otimes \cH')$ multiplies the $\U(1)$-scalars in each fiber over $\PU^\pm_{\cCl_{p,q}}(\cH)\times \PU^\pm_{\cCl_{p,q}}(\cH')$ as in~\eqref{eq:tensoronU}, and hence corresponds to the tensor product of $\U(1)$-bundles in~\eqref{eq:itsaline}. The compatibility with gradings and Real structures follows directly from~\eqref{eq:UactiononPU} being a homomorphism of Real graded  groups.
\ep

\begin{cor}\label{cor:monoidalgerbe}
Let $\X$ and $\X'$ be Real stacks. Under Lemma~\ref{lem:gradedRealgerbe}, the composition 
\beq\label{eq:tensorproductoftwists}
&&\X\times \X'\xrightarrow{\tau\times \tau'} \pt\sq \PU^\pm_{\cCl_{p,q}}(\cH)\times \pt\sq \PU^\pm_{\cCl_{p',q'}}(\cH)\xrightarrow{\otimes} \pt\sq \PU^\pm_{\cCl_{p,q}\otimes \cCl_{p',q'}}(\cH)
\eeq
classifies the external tensor product of the Real graded  gerbes classified by $\tau$ and $\tau'$, and 
$$
\X\xrightarrow{\Delta} \X\times \X\xrightarrow{\tau\times \tau'} \pt\sq \PU^\pm_{\cCl_{p,q}}(\cH)\times \pt\sq \PU^\pm_{\cCl_{p',q'}}(\cH)\xrightarrow{\otimes} \pt\sq \PU^\pm_{\cCl_{p,q}\otimes \cCl_{p',q'}}(\cH)
$$
classifies the (internal) tensor product of Real graded  gerbes on $\X$, where $\Delta$ is the diagonal map. 
\end{cor}
\bp
The second statement follows from the first, for which it suffices to verify the universal example; this is the content of Lemma~\ref{lem:tensortwist}. 
\ep

\begin{rmk} \label{rmk:orientifold}
The $C_2$-action on $\U^\pm_{\cCl_{p,q}}(\cH)$ can also be used to form the semidirect product $\U^\pm_{\cCl_{p,q}}(\cH)\rtimes C_2$ and a central extension analogous to~\eqref{eq:itsaline}. Passing to stacks, we obtain
\beq\label{eq:Realgerbes2ways}
&&\begin{tikzpicture}[baseline=(basepoint)];
\node (Am) at (-1.2,0) {$C_2\circlearrowright$};
\node (A) at (0,0) {$\pt\sq \U(1)$};
\node (B) at (3,0) {$\pt\sq \U_{\cCl_{p,q}}^\pm(\cH)$};
\node (C) at (6,0) {$\pt\sq \PU^\pm_{\cCl_{p,q}}(\cH)$};
\node (D) at (9,0) {$\pt$};
\node (AA) at (0,-1.5) {$\pt\sq \U(1)$};
\node (BB) at (3,-1.5) {$\pt\sq \U^\pm_{\cCl_{p,q}}(\cH)\rtimes C_2$};
\node (CC) at (6,-1.5) {$\pt\sq \PU^\pm_{\cCl_{p,q}}(\cH)\rtimes C_2$};
\node (DD) at (9,-1.5) {$\pt\sq C_2$};
\draw[->] (A) to  (B);
\draw[->] (B) to (C);
\draw[->] (C) to (D);
\draw[->] (AA) to  (BB);
\draw[->] (BB) to (CC);
\draw[->] (CC) to (DD);
\draw[double equal sign distance] (A) to  (AA);
\draw[->] (B) to (BB);
\draw[->] (C) to (CC);
\draw[->] (D) to (DD);
\path (0,-.75) coordinate (basepoint);
\end{tikzpicture}
\eeq
where the top row is the fibration in $C_2$-stacks associated to~\eqref{eq:itsaline}, with $C_2$-quotient the lower row. Conversely, the maps to $\pt\sq C_2$ in the lower row classify principal $C_2$-bundles with total stack the upper row. The top row defines a Real graded  gerbe, whereas the lower row encodes the same data in stacks over $\pt\sq C_2$, i.e., \emph{orientifolds}. In terms of groupoids, the lower row of~\eqref{eq:Realgerbes2ways} is a \emph{$\phi$-twisted graded central extension} \cite[Definition 7.23.iii]{FreedMoore} of $\pt\sq \PU^\pm_{\cCl_{p,q}}(\cH)\rtimes C_2$ where $\phi\colon \pt\sq \PU^\pm_{\cCl_{p,q}}(\cH)\rtimes C_2\to \pt\sq C_2$ is determined by the evident homomorphism. Taking $C_2$-quotients in Corollary~\ref{cor:gradedRealgerbe} leads to a variant in which Real stacks are replaced by stacks over $\pt\sq C_2$. We refer to~\S\ref{sec:appenRealgerbe} for further discussion. 
\end{rmk}

\begin{rmk}\label{rmk:cohomology}
Isomorphism classes of gerbes are classified by certain cohomology groups. Such invariants pull back along maps of stacks, and the Real graded  gerbes in Corollary~\ref{cor:gradedRealgerbe} are classified by the pullback of universal cohomological invariants on the Real stack $\pt\sq \PU_{\cCl_{p,q}}^\pm(\cH)$. The cohomology of a global quotient agrees with the cohomology of the Borel construction, and the relevant invariants are elements of the set
\beq\label{eq:setofinvariants}
\HH^3(B\PU_{\cCl_{p,q}}^\pm(\cH);\underline{\Z})\times \HH^1(B\PU_{\cCl_{p,q}}^\pm(\cH);\Z/2)\simeq \Z\times \Z/2
\eeq
where $\underline{\Z}$ is the $C_2$-equivariant local system on $B\PU_{\cCl_{p,q}}^\pm(\cH)$ for the involution $n\mapsto -n$ on~$\Z$. This local system comes from $\U(1)$ with its complex conjugation action via the $C_2$-equivariant exponential sequence $\Z\to \R\to \U(1)$. The monoidal structure on twists is compatible with the group structure on the set~\eqref{eq:setofinvariants} for the extension classified by the cocycle
\beq\label{eq:groupstruc}
&&\begin{array}{cccc} \HH^1(B\PU_{\cCl_{p,q}}^\pm(\cH);\Z/2)\times \HH^1(B\PU_{\cCl_{p,q}}^\pm(\cH);\Z/2) &\xrightarrow{\bigcup}&\HH^2(B\PU_{\cCl_{p,q}}^\pm(\cH);\Z/2)\\ 
&\xrightarrow{\beta} &\HH^3(B\PU_{\cCl_{p,q}}^\pm(\cH);\underline{\Z})
\end{array}
\eeq
given by the composition of the cup product and the Bockstein homomorphism; see~\cite[\S2.4]{FHTI}, \cite[Remark 7.28]{FreedMoore} or \cite[Proposition 2.13]{Gomi} for details. 
\end{rmk}

\subsection{The category of Real twistings}

\begin{defn} \label{defn:twist}
For a Real local quotient stack $\X$, the groupoid $\Tw_\X$ of \emph{twists} is the category whose objects are maps of stacks
\beq\label{eq:twistingdefn}
 \tau\colon \X\to \pt\sq \PU^\pm_{\cCl_{p,q}}(\H)
\eeq
for a Real Clifford module $\cH$ and Real map $\tau$. Morphisms are 2-commuting triangles
\beq\label{eq:twistmorphisms}
\begin{tikzpicture}[baseline=(basepoint)];
\node (A) at (0,-.75) {$\X$};
\node (B) at (2.5,0) {$\pt \sq \PU^\pm_{\cCl_{p,q}}(\H)$};
\node (C) at (2.5,-1.25) {$\pt \sq \PU^\pm_{\cCl_{p,q}}(\H')$};
\draw[->] (A) to  node [above] {$\tau$} (B);
\draw[->] (A) to node [below] {$\tau'$} (C);
\draw[->] (B) to (C);
\path (0,-.75) coordinate (basepoint);
\end{tikzpicture}
\eeq
with vertical arrow induced by an isometric embedding $\H\subset \H'$ of Real Clifford modules. The assignment $\X\mapsto \Tw_\X$ determines a presheaf of categories on Real stacks. 
\end{defn}

Let $\Twist_\X^\times$ denote the subcategory of $\Twist_\X$ whose morphisms are isomorphisms, i.e., the vertical arrow in~\eqref{eq:twistmorphisms} is induced by an isomorphism of Clifford modules.

Define functors $\boxtimes \colon \Tw_\X\times \Tw_{\X'}\to \Tw_{\X\times \X'}$ that on objects send twists $\tau$ and $\tau'$ to the \emph{external tensor product}, denoted $\tau\boxtimes \tau'$, and given by the composition~\eqref{eq:tensorproductoftwists}. The value of $\boxtimes$ on morphisms is defined analogously using the graded tensor product of operators. Similarly, the composition
\beq\label{eq:tensoroftwists}
\otimes \colon \Tw_\X\times \Tw_{\X}\xrightarrow{\boxtimes} \Tw_{\X\times \X}\xrightarrow{\Delta^*} \Tw_\X
\eeq
of the external tensor product and pullback along the diagonal $\Delta\colon \X\to \X\times \X$ endows $\Tw_\X$ with a monoidal structure. The symmetry of the tensor product of $\Z/2$-Real graded  Hilbert spaces further endows $\Tw_\X$ with the structure of a symmetric monoidal category.

\begin{defn} \label{defn:trivtwist}
A \emph{trivialization} of a twist $\tau$ is a 2-commuting diagram of Real stacks
\beq
\begin{tikzpicture}[baseline=(basepoint)];
\node (B) at (4,0) {$\pt\sq \U^+_{\cCl_{p,q}}(\cH)$};
\node (D) at (4,-1) {$\pt\sq \PU^\pm_{\cCl_{p,q}}(\H)$};
\node (C) at (0,-1) {$\X$};
\draw[->,dashed] (C) to (B);
\draw[->] (B) to (D);
\draw[->] (C) to node [below] {$\tau$} (D);
\path (0,-.75) coordinate (basepoint);
\end{tikzpicture}\nonumber
\eeq
where the lift is to the classifying stack for the subgroup $\U^+_{\cCl_{p,q}}(\cH)<\U^\pm_{\cCl_{p,q}}(\cH)$ of even, Clifford linear unitary operators on $\cH$. 
\end{defn}

\begin{rmk}\label{eq:unitarybundles}
The topological groups $\U^+_{\cCl_{p,q}}(\cH)$ are contractible~\cite[Proposition A2.1]{AtiyahSegaltwists}. Hence, when $X$ is a space every map to $\pt\sq \U^+_{\cCl_{p,q}}(\cH)$ is null-homotopic. However, this is typically false for maps $\X\to \pt\sq \U^+_{\cCl_{p,q}}(\cH)$ out of a topological stack $\X$. For example, non-isomorphic unitary $G$-representations determine non-homotopic maps $\pt\sq G\to \pt\sq \U^+(\cH)$. 
\end{rmk}

\begin{lem}
There is a weak symmetric monoidal 2-functor from $\Twist_\X^\times$ to the 2-category of Real graded  gerbes over $\X$ that sends trivializable twists to trivializable Real gerbes. 
\end{lem}

\bp
This follows from and Corollaries~\ref{cor:gradedRealgerbe} and~\ref{cor:monoidalgerbe}: pullbacks of gerbes are unique up to 1-isomorphism which is itself unique up to further 2-isomorphism, and hence the pullback of Real graded  gerbes constructs a (weak) functor from the 1-groupoid of maps $\X\to \pt\sq \PU^\pm(\cH)$ to the 2-groupoid of gerbes. Compatibility data of the monoidal structure also follows from the universal property of pullbacks. Again by the universal property of pullbacks, a trivialization in the sense of Definition~\ref{defn:trivtwist} leads to a trivializing section of a the graded Real gerbe pulled back along~$\tau$. 
\ep

\begin{rmk}
We emphasize that under the above 2-functor non-isomorphic maps $\X\to \pt\sq \PU^\pm(\cH)$ can lead to 1-isomorphic Real graded  gerbes owing to the existence of non-isomorphic maps to $\pt\sq \U^+(\cH)$, see Remark~\ref{eq:unitarybundles}. However, for each Real graded  gerbe over a local quotient stack there is an essentially unique twist; see Corollary~\ref{cor:absorbing} below. 
\end{rmk}

\subsection{Twisted Real Clifford bundles}
\begin{defn}
Let $\X$ be a Real stack. For a map $\tau\colon \X\to \pt\sq \PU^\pm_{\cCl_{p,q}}(\cH)$ the associated \emph{twisted Real Clifford bundle} is the pullback in stacks
\beq\label{eq:twistHilb}
&&
\begin{tikzpicture}[baseline=(basepoint)];
\node (A) at (0,0) {$\cH^\tau$};
\node (B) at (3,0) {$\cH\sq \U^\pm_{\cCl_{p,q}}(\H)$};
\node (C) at (0,-1.5) {$\X$};
\node (D) at (3,-1.5) {$\pt\sq \PU^\pm_{\cCl_{p,q}}(\H).$};
\node (P) at (.7,-.4) {\scalebox{1.5}{$\lrcorner$}};
\draw[->] (A) to  (B);
\draw[->] (B) to  (D);
\draw[->] (A) to  (C);
\draw[->] (C) to node [below] {$\tau$} (D);
\path (0,-.75) coordinate (basepoint);
\end{tikzpicture}
\eeq
When $\cCl_{p,q}=\C$, a twisted Real Clifford bundle is a \emph{twisted Real Hilbert bundle.} 
\end{defn}

To unpack this, a twisted Real Hilbert bundle $\cH^\tau\to \X$ is a Real Hilbert bundle $\cH^\tau\to \X^\tau$ classified by $\X^\tau\to \pt\sq \U^\pm(\cH)$ for $\X^\tau\to \X$ the Real graded  gerbe classified by $\tau$ under Lemma~\ref{lem:gradedRealgerbe}. Furthermore, for any point $\pt\to \X$, the fiber of $\cH^\tau$ at $\pt\sq \U(1)\hookrightarrow \X^\tau$ is the Real Clifford module $\cH$ with $\U(1)$ acting by scalars, and the Real structure on $\cH$ is compatible with the complex conjugation action on $\U(1)$. We refer to Lemma~\ref{lem:groupoidHilbbundle} an explicit groupoid description. 

\begin{rmk}
Twisted Real Clifford bundles and twisted Real Hilbert bundles are the objects of a bicategory whose 1-morphisms are maps fiberwise equivalences of stacks $f\colon \cH^\tau\to \cH^{\tau'}$ covering maps of Real graded  gerbes $\X^\tau\to \X^{\tau'}$.
\end{rmk}

\begin{cor}
The external tensor product of twisted Real Clifford bundles $\cH^\tau\boxtimes \cH'^{\tau'}$ is a twisted Real Clifford bundle over the external tensor product of graded gerbes $\X^\tau\boxtimes \X'^{\tau'}$. Similarly, the internal tensor product of $\cH^\tau\otimes \cH'^{\tau'}$ is a twisted Real Clifford bundle over the internal tensor product of graded gerbes $\X^\tau\otimes \X'^{\tau'}$.
\end{cor}
\bp 
We again appeal to the universal example, where the first claim follows from the pullback in Real stacks
\beq\nonumber
&&
\begin{tikzpicture}[baseline=(basepoint)];
\node (A) at (0,0) {$(\cH\sq \U^\pm_{\cCl_{p,q}}(\H))\boxtimes (\cH'\sq \U^\pm_{\cCl_{p',q'}}(\H'))$};
\node (C) at (0,-1.5) {$\pt\sq \PU^\pm_{\cCl_{p,q}}(\H)\times \pt\sq \PU^\pm_{\cCl_{p',q'}}(\H')$};
\node (B) at (5.5,0) {$(\cH\otimes \cH')\sq \U^\pm(\H\otimes \H')$};
\node (D) at (5.5,-1.5) {$\pt\sq \PU^\pm_{\cCl_{p,q}\otimes \cCl_{p',q'}}(\H\otimes \H').$};
\node (P) at (.7,-.4) {\scalebox{1.5}{$\lrcorner$}};
\draw[->] (A) to  (B);
\draw[->] (B) to  (D);
\draw[->] (A) to  (C);
\draw[->] (C) to node [above] {$\otimes$} (D);
\path (0,-.75) coordinate (basepoint);
\end{tikzpicture}
\eeq
The second claim follows after further pulling back along the diagonal $\pt\sq \PU^\pm_{\cCl_{p,q}}(\H)\xrightarrow{\Delta} \pt\sq \PU^\pm_{\cCl_{p,q}}(\H)\times \pt\sq \PU^\pm_{\cCl_{p,q}}(\H)$. 
\ep

Non-isomorphic twisted Real Clifford bundles can be defined over the same underlying Real graded  gerbe. Indeed, the homomorphisms~\eqref{eq:UactiononPU} determine maps of stacks
$$
\pt\sq (\U^\pm(\H')\times \PU^\pm(\H))\to \pt\sq \PU^\pm(\H'\otimes \H),
$$
that tensor a Hilbert bundle and a projective Hilbert bundle. The pullback gerbe in~\eqref{eq:twisttogerbe} is unchanged (up to isomorphism) by this tensor product with a Hilbert bundle, but in general a map $\tau\colon \X\to \pt\sq \PU^\pm(\cH)$ will not be isomorphic to its tensor product with a classifying map for a (non-projective) Hilbert bundle. The following identifies the subcategory with this special property, compare~\cite[Definition 6.8]{Heinloth}. 

\begin{defn} A Real map $\tau\colon \X\to \pt\sq \PU^\pm_{\cCl_{p,q}}(\H)$ is \emph{absorbing} if for any Real Hilbert bundle classified by $\rho\colon \X\to \pt\sq \U^+(\cH')$ there is a 2-commuting diagram in Real stacks
\beq\nonumber
\begin{tikzpicture}[baseline=(basepoint)];
\node (A) at (0,0) {$\X$};
\node (B) at (4,0) {$\pt\sq (\PU^\pm_{\cCl_{p,q}}(\H)\times \U^+(\H'))$};
\node (C) at (0,-1.5) {$\pt\sq \PU^\pm_{\cCl_{p,q}}(\H)$};
\node (D) at (4,-1.5) {$\pt\sq \PU^\pm_{\cCl_{p,q}}(\H\otimes \H')$};
\draw[->] (A) to node [above] {$\tau\times \rho$} (B);
\draw[->] (B) to  (D);
\draw[->] (A) to node [left] {$\tau$} (C);
\draw[->] (D) to node [below] {$\simeq$} (C);
\path (0,-.75) coordinate (basepoint);
\end{tikzpicture}\eeq
where the lower horizontal arrow is induced by a choice of isomorphism $\H\otimes \H'\simeq \H$ (which requires $\H$ to be infinite rank). A Real map $\tau\colon \X\to \pt\sq \PU^\pm_{\cCl_{p,q}}(\H)$ is \emph{locally absorbing} if this property holds after restricting $\tau$ to any closed substack of $\X$. 
\end{defn}

In other words, a twisted Real Clifford bundle $\cH^\tau \to \X$ is absorbing if for all Hilbert bundles $\cH'\to \X$ there exists an isomorphism $\cH^\tau\otimes \cH'\simeq \cH^\tau$ of twisted Real Clifford bundles over $\X$. The following crucially uses that $\X$ is a local quotient stack.  

\begin{prop}\label{prop:universal}
A Real graded  gerbe on a Real local quotient stack $\X$ determines a locally absorbing Real Clifford bundle classified by a map $\tau\colon \X\to \pt\sq \PU^\pm_{\cCl_{p,q}}(\cH)$ that is unique up to isomorphism. 
\end{prop}
\bp
The following argument is essentially \cite[Lemma~3.11]{FHTI}, see also \cite[Proposition 6.12]{Heinloth}. Given a gerbe on $\X$, Lemma~\ref{lem:BX} constructs an atlas $\bX_0\to \X$ on which the gerbe is equivalent to a Real graded  central extension $\widehat{\bX}\to \bX$ for the presentation $\bX=\{\bX_1\rightrightarrows \bX_0\}$ of~$\X$ with strict $C_2$-action on~$\bX$~\cite[Proposition 1.5]{Romagny}. Then~\cite[Corollary A.33]{FHTI} constructs a locally absorbing Hilbert bundle $\cH\to \widehat{\bX}$
that is unique up to unitary equivalence. Since the grading on $\widehat{\bX}$ is determined by a map $\epsilon_1\colon \widehat{\bX}_1\to \Z/2$ (and in particular, is constant on $\widehat{\bX}_0$), by \cite[Proposition A2.1]{AtiyahSegaltwists} we may choose a trivialization of $\cH$ on $\widehat{\bX}_0=\bX_0$. With such a trivialization fixed, $\cH$ is equivalent data to a continuous, $C_2$-equivariant map $\widehat{\bX}_1\to \U^\pm_{\cCl_{p,q}}(\H)$ compatible with the grading, and satisfying a cocycle condition on $\widehat{\bX}_2=\widehat{\bX}_1\times_{\bX_0}\widehat{\bX}_1$.
 
By construction, $\cH$ has a fiberwise Real action by~$\U(1)$. Let $\cH(1)\subset \cH$ denote the subbundle where this fiberwise action is through scalars, i.e., the weight~1 subspace. Then $\cH(1)$ defines a $\U(1)$-equivariant map $\widehat{\tau}_1\colon \widehat{\bX}_1\to \U^\pm_{\cCl_{p,q}}(\H)$. Taking quotients by this action gives a continuous $C_2$-equivariant map $\tau_1\colon \bX_1\to \PU^\pm_{\cCl_{p,q}}(\H)$. The cocycle condition for $\widehat{\tau}_1$ implies the cocycle condition for $\tau_1$ on ${\bX}_2={\bX}_1\times_{\bX_0}{\bX}_1$, and hence we have constructed a continuous $C_2$-equivariant functor of topological groupoids  $\bX\to \pt\sq \PU^\pm_{\cCl_{p,q}}(\H)$. This determines maps of stacks $\tau\colon \X\to \pt\sq \PU^\pm_{\cCl_{p,q}}(\H)$, which is locally absorbing because the Hilbert bundle~$\cH$ was locally absorbing. 
\ep

\begin{cor}\label{cor:absorbing}
For every $p,q\in \mathbb{N}$ and every Real graded  gerbe on a local quotient stack~$\X$, there is a Real twisting $\tau\colon \X\to \pt\sq \PU^\pm_{\cCl_{p,q}}(\cH)$ that is unique up to isomorphism. 
\end{cor}
\bp
For a given Real graded  gerbe, take the locally absorbing twist. 
\ep

\subsection{Groupoid presentations for Real graded  gerbes and twisted Hilbert bundles}\label{sec:Ktheorygrpod1}

\begin{lem}\label{lem:grpd1}
Given a Real map $\tau\colon \X\to \pt\sq \PU^\pm_{\cCl_{p,q}(\cH)}$, there exists a groupoid presentation~$\bX$ of~$\X$ for which $\tau$ determines a Real graded  central extension of $\bX$. Furthermore, an isomorphism between Real maps $\tau\colon \X\to \pt\sq \PU^\pm_{\cCl_{p,q}(\cH)}$ determines an isomorphism of Real graded  central extensions of~$\bX$. 
\end{lem}
\bp
This follows from pulling back the Real graded gerbe in Lemma~\ref{lem:gradedRealgerbe} along $\tau$ and applying Lemma~\ref{lem:BX}.
\ep

\begin{rmk}\label{rmk:lem:grpd1}
 Lemma~\ref{lem:BX} can made totally explicit in the setting of Lemma~\ref{lem:grpd1}, affording a description of the Real graded central extension that can be helpful in examples. To explain this, first let $\tau_1\colon \bX_1\to \PU^\pm_{\cCl_{p,q}}(\cH)$ be a continuous map determining a (strictly) $C_2$-equivariant continuous functor $\bX\to \pt\sq \PU^\pm_{\cCl_{p,q}}(\cH)$ between Real topological groupoids whose underlying map of Real stacks is~$\tau\colon \X\to \pt\sq \PU^\pm_{\cCl_{p,q}(\cH)}$; the existence of this functor follows from~\cite[Proposition 1.5]{Romagny}. Then consider the diagrams in $C_2$-spaces
\beq\label{eq:gerbepullbacksFHT}
&&\begin{tikzpicture}[baseline=(basepoint)];
\node (A) at (0,0) {$\bX_1^\tau$};
\node (B) at (3,0) {$\U^\pm_{\cCl_{p,q}}(\H)$};
\node (C) at (0,-1.5) {$\bX_1$};
\node (D) at (3,-1.5) {$ \PU^\pm_{\cCl_{p,q}}(\H)$};
\node (DD) at (4.8,-1.5) {$ \Z/2$};
\node (P) at (.7,-.4) {\scalebox{1.5}{$\lrcorner$}};
\draw[->] (A) to node [above] {$\widehat{\tau}_1$} (B);
\draw[->] (B) to  (D);
\draw[->] (A) to  (C);
\draw[->] (C) to node [above] {$\tau_1$} (D);
\draw[->,bend right=25] (C) to node [below] {$\epsilon_1$} (DD);
\draw[->] (D) to node [above] {$\pi_0$} (DD);
\path (0,-.75) coordinate (basepoint);
\end{tikzpicture}\hspace{-.2in} \begin{tikzpicture}[baseline=(basepoint)];
\node (A) at (0,0) {$\bX_2=\bX_1\times_{\bX_0} \bX_1$};
\node (B) at (4.5,0) {$\PU^\pm_{\cCl_{p,q}}(\H)\times \PU^\pm_{\cCl_{p,q}}(\H)$};
\node (C) at (0,-1.5) {$\bX_1$};
\node (D) at (4.5,-1.5) {$\PU^\pm_{\cCl_{p,q}}(\H).$};
\draw[->] (A) to node [above] {$\tau_1\times \tau_1$} (B);
\draw[->] (B) to node [right] {$m$} (D);
\draw[->] (A) to node [left] {$c$} (C);
\draw[->] (C) to node [below] {$\tau_1$} (D);
\path (0,-.75) coordinate (basepoint);
\end{tikzpicture}
\eeq
 The pullback square on the left is a $C_2$-equivariant graded $\U(1)$-bundle $\bX^\tau_1\to \bX_1$ where the $C_2$-action on fibers is by complex conjugation. The square on the right produces the data of a $C_2$-equivariant isomorphism of graded $\U(1)$-bundles covering composition in~$\bX$, which completes the data of a Real graded central extension of $\bX$ representing the Real graded gerbe classified by~$\tau$. 
\end{rmk}

\begin{lem}\label{lem:groupoidHilbbundle}
In the groupoid presentation from the previous lemma, the twisted Real Clifford bundle $\cH^\tau\to \X$ determines a Real Hilbert bundle over the Real graded  central extension~$\bX^\tau$ of $\bX$, where the value of the grading $\epsilon_\tau\colon \bX^\tau\to \Z/2$ and determines whether a morphism acts by a grading preserving or grading reversing isometry of $\cH$. 
\end{lem}

\bp
The map $\widehat{\tau}_1$ in the left square of~\eqref{eq:gerbepullbacksFHT} gives the claimed Real Clifford bundle over~$\bX$. The grading on $\bX^\tau$ comes from the composite $\bX_1^\tau\to \U^\pm_{\cCl_{p,q}}(\cH)\xrightarrow{\epsilon} \Z/2$, verifying the grading claim. 
\ep

\subsection{Comparison with other categories of K-theory twists}

\begin{lem}
The lower row in~\eqref{eq:Realgerbes2ways} is a $\phi$-twisted graded central extension of $\pt\sq \PU^\pm_{\cCl_{p,q}}(\cH)\rtimes C_2$ in the sense of ~\cite[Definition 7.23.iii]{FreedMoore} and~\cite[Definition~2.3]{Gomi}. The upper row in~\eqref{eq:Realgerbes2ways} is a Real graded  $S^1$-groupoid in the sense of~\cite[Definition 3.1]{Moutuou}.
\end{lem}

\bp
After unwinding the definitions in these references, this follows immediately from viewing~\eqref{eq:itsaline} as a Real graded  $\U(1)$-bundle with multiplicative structure~\eqref{eq:multiplicativeline}.
%
\ep

\begin{cor}\label{cor:twistcompare}
The pullback along a $C_2$-equivariant map of stacks $\X\to \pt\sq \PU^\pm(\cH)$ determines a Real twisting in the sense \cite{Moutuou}. Passing to $C_2$-quotients, the map $\X\sq C_2\to \pt\sq \PU^\pm(\cH)\rtimes C_2$ determines a twisting in the sense of \cite[Definition 7.23]{FreedMoore}, \cite[Definition 2.3]{Gomi}. 
\end{cor}
\bp
This follows formally from naturality: in all definitions, Real graded  central extensions of groupoids pull back along functors~\cite[Remark~7.31]{FreedMoore} and \cite[\S2.7]{MoutuouThesis}. One can also unpack this concretely as in~\eqref{eq:gerbepullbacksFHT}, where the line bundles $\bY^\tau_1\to \bY_1$ acquire a Real structure for $C_2$-action on $\bY_1$ coming from a strictification of the $C_2$-action on $\X$ (using~\cite{Romagny}). 
\ep


\begin{lem}
The pullback in Real stacks~\eqref{eq:twistHilb} determines an Rg module over a groupoid presentation of $\X$ in the sense of \cite[Definition 8.1.2]{MoutuouThesis}. The map on $C_2$-quotients $\cH^\tau\sq C_2\to \X\sq C_2$ gives a $\phi$-twisted graded Hilbert bundle in the sense of \cite[Definition 7.23.iv]{FreedMoore} and \cite[\S2.5]{Gomi}. Hence, these pullbacks determine twists for KR-theory in the sense \cite{MoutuouThesis,FreedMoore,Gomi}. 
\end{lem}
\bp
By construction, we have the factorization $\cH^\tau\to \X^\tau\to \X$ where $\X^\tau \to \X$ is a Real graded  gerbe and $\cH^\tau\to \X^\tau$ is a Real graded  Hilbert bundle classified by the $\pt\sq \U(1)$-equivariant map $\X^\tau\to \pt\sq \U^\pm(\cH)$. Choosing a groupoid presentation for this classifying map, we obtain a $\U(1)$-equivariant map of spaces $\bX_1^\tau\to \U^\pm(\cH)$, i.e., the $\U(1)$ in the Real graded  central extension of $\bX$ acts by scalars on the Hilbert bundle $\cH^\tau\to \bX$. Ignoring Real structures, the statement now follows immediately from \cite[Remark 3.7]{FHTI}. Adding Real structures, one compares with the more involved definition, but this is also straightforward, e.g., see \cite[Definition 2.8]{Kubota}. 
\ep

\section{Twisted $\KR$-theory of topological stacks}\label{sec:KR}

\subsection{Real Fredholm operators}\label{sec:Fred}
We follow conventions from~\cite[Remark 4.12]{Gomi} below. 
For a Real, unitary right $\cCl_{p,q}$-module $\cH$ with conventions as in~\S\ref{sec:RealCliff}, let $\Fred_{\cCl_{p,q}}(\cH)$ denote the set of odd, self-adjoint Fredholm operators $F\colon \cH\to \cH$, that commute with the right $\cCl_{p,q}$-action, $F^2-\id$ compact, and the spectrum of $F$ is in $[-1,1]$. Take the subspace topology on $\Fred_{\cCl_{p,q}}(\cH)$ for the inclusion
$$
\Fred_{\cCl_{p,q}}(\cH)\to B(\cH)\times K(\cH),\qquad F\mapsto (F,F^2-\id)
$$
where $B(\cH)$ is the space of bounded operators in the compact open topology and $K(\cH)$ is the space of compact operators in the norm topology. The unitary group $\U_{\cCl_{p,q}}^\pm(\cH)$ acts on~$\Fred_{\cCl_{p,q}}(\cH)$ by conjugation,
\beq\label{eq:UactiononFred}
\U_{\cCl_{p,q}}^\pm(\cH)\times \Fred_{\cCl_{p,q}}(\cH)\to \Fred_{\cCl_{p,q}}(\cH),\qquad (u,F)\mapsto uFu^{-1}
\eeq
and this action is continuous for $\U_{\cCl_{p,q}}^\pm(\cH)$ in the compact open (or equivalently, strong) topology \cite[\S{A}.2]{AtiyahSegaltwists}. Furthermore,~\eqref{eq:UactiononFred} is $C_2$-equivariant for the action on operators inherited from the Real structure on $\cH$, and factors through the quotient by $\U(1)$-scalars to give the top $C_2$-equivariant map 
\beq\label{eq:PactiononFred}
\begin{tikzpicture}[baseline=(basepoint)];
\node (A) at (0,0) {$\PU_{\cCl_{p,q}}^\pm(\cH)\times \Fred_{\cCl_{p,q}}(\cH)$};
\node (B) at (4,0) {$\Fred_{\cCl_{p,q}}(\cH)$};
\node (C) at (0,-1.5) {$\PO_{\Cl_{p,q}}^\pm(\cH_\R)\times \Fred_{\Cl_{p,q}}(\cH_\R)$};
\node (D) at (4,-1.5) {$\Fred_{\Cl_{p,q}}(\cH_\R)$};
\draw[->] (A) to (B);
\draw[->] (D) to  (B);
\draw[->] (C) to (A);
\draw[->] (C) to (D);
\path (0,-.75) coordinate (basepoint);
\end{tikzpicture}\qquad (u,F)\mapsto uFu^{-1}\eeq
where the lower row is the $C_2$-fixed points, giving the action of the projective orthogonal group on real Fredholm operators on the $\Cl_{p,q}$-module $\cH_\R:=\cH^{C_2}$. 

\begin{defn}
For $\cH$ a Real $\cCl_{p,q}$-module, let $\Fred_{\cCl_{p,q}}(\cH)\sq \PU^\pm_{\cCl_{p,q}}(\cH)$ denote the Real stack associated with the $C_2$-equivariant map~\eqref{eq:PactiononFred}. 
\end{defn}

\begin{rmk}
The above is the Real generalization of a Fredholm operator construction of twisted complex K-theory from \cite[\S{A}.2]{AtiyahSegaltwists} and \cite[\S{A}.5]{FHTI}. A comprehensive reference for this Real generalization is \cite[\S{B}.2 and Remark 4.12]{Gomi}. An early reference for Fredholm operators as a representing space for KR-theory is~\cite[Theorem~II]{Matumoto}. Similar ideas using representing stacks for $\KR$-theory are in \cite[Example 1.3.19]{Sati2021EquivariantPI}. 
\end{rmk}

For an isometric inclusion $\H\hookrightarrow \H'$ of right $\cCl_{p,q}$-modules, there is a 2-commutative square of Real stacks
\beq\nonumber
\begin{tikzpicture}[baseline=(basepoint)];
\node (A) at (0,0) {$\Fred_{\cCl_{p,q}}(\cH)\sq \PU_{\cCl_{p,q}}^\pm(\cH)$};
\node (B) at (5,0) {$\Fred_{\cCl_{p,q}}(\cH')\sq \PU_{\cCl_{p,q}}^\pm(\cH')$};
\node (C) at (0,-1.25) {$\pt \sq \PU_{\cCl_{p,q}}^\pm(\cH)$};
\node (D) at (5,-1.25) {$\pt \sq \PU_{\cCl_{p,q}}^\pm(\cH')$};
\draw[->] (A) to  (B);
\draw[->] (B) to (D);
\draw[->] (A) to (C);
\draw[->] (C) to (D);
\path (0,-.75) coordinate (basepoint);
\end{tikzpicture}
\eeq
determined by strictly commuting squares of topological groupoids. The lower horizontal arrows are determined by the homomorphism~\eqref{eq:embedding1}, while the upper arrows are given by
\beq\label{eq:stabilizeFred}
&&\Fred_{\cCl_{p,q}}(\H)\to \Fred_{\cCl_{p,q}}(\H'),\qquad F\mapsto F\oplus \epsilon, \quad \H\oplus \H^\perp\simeq \H'
\eeq
where $\H^\perp\subset \H$ is the orthogonal complements to the embedding and $\epsilon\in \Fred_{\cCl_{p,q}}(\H^\perp)$ is a basepoint, i.e., an invertible Fredholm operator. The space of invertible Fredholm operators is contractible, so the map~\eqref{eq:stabilizeFred} is unique up to a contractible space of choices.

The tensor product of Clifford linear Fredholm operators 
\beq\label{eq:Fredtensor}
\Fred_{\cCl_{p,q}}(\H)\times \Fred_{\cCl_{p,q}}(\H')&\to& \Fred_{\cCl_{p,q}}(\H\otimes \H'),\\
(F,F') &\mapsto& F\otimes \id+\id\otimes F'\nonumber
\eeq
is equivariant relative to the homomorphism~\eqref{eq:tensoronPU1}, giving 2-commutative diagrams in stacks 
\beq\nonumber
\begin{tikzpicture}[baseline=(basepoint)];
\node (A) at (0,0) {$\Fred_{\cCl_{p,q}}(\H)\times \Fred_{\cCl_{p,q}}(\H')\sq (\PU^\pm(\H)\times \PU^\pm(\H'))$};
\node (B) at (7,0) {$\Fred_{\cCl_{p,q}}(\H\times \H')\sq \PU^\pm(\H\otimes \H')$};
\node (C) at (0,-1.25) {$\pt \sq (\PU^\pm(\H)\times \PU^\pm(\H'))$};
\node (D) at (7,-1.25) {$\pt \sq \PU^\pm(\H\otimes\H').$};
\draw[->] (A) to  (B);
\draw[->] (B) to (D);
\draw[->] (A) to (C);
\draw[->] (C) to (D);
\path (0,-.75) coordinate (basepoint);
\end{tikzpicture}
\eeq


\begin{rmk}\label{rmk:vectorbundle}
A finite rank Real (super) vector bundle $V\to X$ determines a map $X\to \Fred_{\cCl_{p,q}}(\H)$ that is unique up to a contractible space of choices. Indeed, there is a contractible space of Real embeddings $V\hookrightarrow X\times \H$ over $X$, realizing $V$ as a summand in the trivial Real Hilbert bundle over~$X$, so $X\times \H\simeq V\oplus V^\perp$. With respect to this decomposition, define a Real Fredholm operator $F=0_V\oplus \id_{V^\perp}$ whose kernel is $V$. 
\end{rmk}


\subsection{Twisted KR-theory}

\begin{defn}[Twisted KR, KO, and KU]\label{defn:twistedK} 
For a Real local quotient stack $\X$ and a Real map $\tau\colon \X\to \pt\sq \PU^\pm_{\cCl_{p,q}}(\cH)$, define the \emph{$\tau$-twisted KR-theory of $\X$} as the homotopy classes of Real (i.e., $C_2$-invariant) sections
\beq\label{eq:KRdefn}
\KR^{\tau}(\X):=\pi_0(\Gamma(\X;\tau^*\Fred_{\cCl_{p,q}}(\H))^{C_2}).
\eeq
The \emph{$\tau$-twisted $\KO$-theory} of a local quotient stack $\X$ is the $\KR$-theory of $\X$ for the trivial Real structure, i.e., trivial $C_2$-action, 
$$
\KO^\tau(\X):=\KR^{\tau}(\X)=\pi_0(\Gamma(\X;\tau^*\Fred_{\cCl_{p,q}}(\H))^{C_2}).
$$
For a local quotient stack $\X$ and a map $\tau\colon \X\to \pt\sq \PU^\pm_{\cCl_{p,q}}(\cH)$, define the \emph{$\tau$-twisted KU-theory of $\X$} as the homotopy classes of sections
$$
\KU^{\tau}(\X):=\pi_0(\Gamma(\X;\tau^*\Fred_{\cCl_{p,q}}(\H))),
$$
i.e., the same as~\eqref{eq:KRdefn} but without taking Real ($C_2$-invariant) sections. 
\end{defn}

\begin{rmk}
We remind that spaces of sections as in~\eqref{eq:KRdefn} can be computed relative to any groupoid presentation of $\X$, see~\eqref{eq:spaceofsections} below. 
\end{rmk}

\begin{lem}\label{lem:KRcocycle}
A class in the $\tau$-twisted KR-theory of $\X$ is represented by a Real map,
\beq\label{eq:twistKasmapofstacks}
F\colon \X\to \Fred_{\cCl_{p,q}}(\H)\sq \PU^\pm_{\cCl_{p,q}}(\H)
\eeq
where the twist $\tau$ is the post-composition with the projection to $\pt\sq \PU^\pm_{\cCl_{p,q}}(\H)$. For a stack $\X$ (without specified Real structure), a class in the $\tau$-twisted KU-theory of~$\X$ is a map~\eqref{eq:twistKasmapofstacks} (without $C_2$-equivariance data). Similarly, a class in the $\tau$-twisted KO-theory of~$\X$ is represented by a Real map~\eqref{eq:twistKasmapofstacks} for the trivial $C_2$-action on~$\X$, which is equivalent to a map of stacks
\beq\label{eq:twistKOasmapofstacks}
\X\to \Fred_{\Cl_{p,q}}(\cH_\R)\sq \PO^\pm_{\Cl_{p,q}}(\cH_\R).
\eeq
for $\Fred_{\Cl_{p,q}}(\cH_\R)$ the set of Fredholm operators on the real $\Cl_{p,q}$-module $\cH_\R:=\cH^{C_2}\subset \cH$. 
\end{lem}
\bp
The claims that identify a section of a bundle of Fredholm operators with a map of stacks follow directly from the universal property of pullbacks. It remains to verify that~\eqref{eq:twistKasmapofstacks} reduces to \eqref{eq:twistKOasmapofstacks} for the trivial $C_2$-action on $\X$. Consider the diagram
\beq\label{eq:KOmapdiagraom}
\begin{tikzpicture}[baseline=(basepoint)];
\node (A) at (0,0) {$\X^{C_2}$};
\node (B) at (4,0) {$(\Fred_{\cCl_{p,q}}(\cH)\sq \PU_{\cCl_{p,q}}^\pm(\cH))^{C_2}$};
\node (C) at (0,-1.25) {$\X$};
\node (D) at (4,-1.25) {$\Fred_{\cCl_{p,q}}(\cH)\sq \PU_{\cCl_{p,q}}^\pm(\cH)$};
\draw[->] (A) to  (B);
\draw[->,right hook-latex] (B) to (D);
\draw[->,right hook-latex] (A) to node [left] {$\simeq$} (C);
\draw[->] (C) to (D);
\path (0,-.75) coordinate (basepoint);
\end{tikzpicture}
\eeq
that includes the $C_2$-fixed points. The vertical equivalence on the left follows from the $C_2$-action on $\X$ being trivial. The claim then follows after identifying the (homotopy) $C_2$-fixed points in the upper right of~\eqref{eq:KOmapdiagraom} with the stack underlying the (strict) $C_2$-fixed points in the presenting groupoid \cite[\S2]{Romagny}
$$
(\Fred_{\cCl_{p,q}}(\cH)\sq \PU_{\cCl_{p,q}}^\pm(\cH))^{C_2}\simeq \Fred_{\cCl_{p,q}}(\cH)^{C_2}\sq \PU_{\cCl_{p,q}}^\pm(\cH)^{C_2}\simeq \Fred_{\cCl_{p,q}}(\cH_\R)\sq \PO_{\Cl_{p,q}}^\pm(\cH_\R).
$$
By the direct computation above we obtain the claimed stack: the $C_2$-invariant subalgebra of $\cCl_{p,q}$ is the real Clifford algebra $\Cl_{p,q}$, Fredholm operators commuting with the Real structure can be identified with Fredholm operators on the real Hilbert space $\cH_\R$, and unitary operators commuting with the Real structure can be identified with orthogonal operators on the real Hilbert space $\cH_\R$. 
\ep


\begin{defn}\label{eq:defnproduct}
The \emph{external product} of twisted K-groups 
\beq
\KR^{\tau}(\X)\times \KR^{\tau'}(\X')&\xrightarrow{\boxtimes} & \KR^{\tau\boxtimes \tau'}(\X\times \X')\nonumber
\eeq
is induced by the tensor product of Fredholm operators~\eqref{eq:Fredtensor}. Similarly, the composition
\beq
\KR^{\tau}(\X)\times \KR^{\tau'}(\X)\xrightarrow{\boxtimes}  \KR^{\tau\boxtimes \tau'}(\X\times \X)\xrightarrow{\Delta^*}\KR^{\tau\otimes \tau'}(\X) \nonumber
\eeq
determines the (internal) product on twisted KR-theory for $\Delta\colon \X\to \X\times \X$ the diagonal. 
\end{defn} 
\begin{rmk}\label{rmk:Hilbertspacedata}
The twisted K-groups  in Definition~\ref{defn:twistedK} are defined relative to a specific twisted Hilbert bundle specified by~$\tau$. This differs from the formalism from~\cite{FHTI} where the twist is labeled by a graded gerbe. These approaches are equivalent if one always chooses a locally absorbing twisted Hilbert bundle from Proposition~\ref{prop:universal}, see also \cite[Remark~3.17]{FHTI}. 
\end{rmk}

\begin{notation}
For a Real stack $\X$, let $\KR^{p,q}(\X)$ denote the twisted KR-theory of $\X$ for the locally absorbing twist classified by the map $\tau\colon \X\to  \pt\sq \U^+_{\cCl_{p,q}}(\cH)\to \PU^\pm_{\cCl_{p,q}}(\cH)$ associated with the trivial Real graded  gerbe on $\X$. 
\end{notation}

\begin{rmk}
In the existing literature, it is perhaps more common to use additive notation for the tensor product of twists~\eqref{eq:tensoroftwists} in Definition~\ref{eq:defnproduct}. This is consistent with the binary operation on the group $\HH^3(\X;\Z)$ that classifies gerbes and their corresponding twists. For our purposes, we prefer the tensor product notation, as it emphasizes that the product of twists is inherited from the tensor product of Hilbert bundles and Fredholm operators. Alas, our choice leads to the somewhat awkward formula $(p,q)\otimes (p',q')\simeq (p+p',q+q')$ for the degree twistings above inherited from the isomorphism~\eqref{eq:tensorofCliffordalgebras}. 
\end{rmk}


\subsection{Computing twisted K-theory relative to a groupoid presentation}\label{sec:Ktheorygrpod}


\begin{prop}
The space of sections~\eqref{eq:KRdefn} is homeomorphic to the space of sections of the Fredholm family in \cite[Remark 7.37]{FreedMoore}, \cite[Definition~3.1]{Gomi}, and also the space defining twisted KR-theory in \cite[Definition 7.3.2]{MoutuouThesis} and \cite[\S6.2]{HMSV}. Finally, these maps are multiplicative: the products from Definition~\ref{eq:defnproduct} are compatible with the previous constructions of twisted KR-theory. 
\end{prop}
\bp
By general descent arguments (e.g., see \cite[Lemma 4.3]{Heinloth} or \cite[\S{A.3}]{FHTI}), one can compute the spaces of sections $\Gamma(\X;\tau^*\Fred_{\cCl_{p,q}}(\H)^{C_2})$ in Definition~\ref{defn:twistedK} using any groupoid presentation of~$\X$. Using Lemma~\ref{lem:grpd1} and Remark~\ref{rmk:lem:grpd1}, choose a continuous $C_2$-equivariant functor $\tau\colon \bX\to \pt\sq \PU^\pm(\H)$ presenting the twist for a Real atlas $\bX_0\to \X$. This provides a Real graded  central extension~$\bX^\tau$ and a $\tau$-twisted Clifford bundle on $\bX$ via~\eqref{eq:gerbepullbacksFHT}. The space of sections computed using this atlas is
\beq\label{eq:spaceofsections}
&&\Gamma(\X;\tau^*\Fred_{\cCl_{p,q}}(\H))^{C_2}\simeq \{F\colon \bX_0\to \Fred_{\cCl_{p,q}}(\H)\mid s^*F=\tau_1\cdot t^*F \cdot \tau_1^{-1} \ {\rm over} \ \bX_1\}^{C_2}
\eeq
for source and target maps $s,t\colon \bX_1\to \bX_0$ and where $\tau_1\colon \bX_1\to \PU^\pm(\H)$ acts by conjugation on Fredholm operators as indicated. This action descends from a $\U^\pm(\H)$-action over the graded central extension $\bX^\tau\to \bX$ from the pullback~\eqref{eq:twistHilb} for $\bX^\tau_1$ defined as the pullback~\eqref{eq:gerbepullbacksFHT}. Using that the quotient by the free $\U(1)$-action on $\bX_1^\tau$ is $\bX_1$, this gives the alternative description
\beq
\Gamma(\X;\tau^*\Fred_{\cCl_{p,q}}(\H))&\simeq& \{F\colon \bX_0\to \Fred_{\cCl_{p,q}}(\H)\mid \widehat{s}^*F=\widehat{\tau}_1\cdot \widehat{t}^*F \cdot \widehat{\tau}_1^{-1} \ {\rm over} \ \widehat{\bX}_1\}^{\U(1)\rtimes C_2}\nonumber\\
&\simeq&\Gamma(\widehat{\bX};\widehat{\tau}^*\Fred_{\cCl_{p,q}}(\H))^{\U(1)\rtimes C_2}\label{eq:FHTmodel}
\eeq
where $\widehat{s},\widehat{t}\colon \widehat{\bX}_1\to \widehat{\bX}_0=\bX_0$ are the source and target maps for a graded central extension and $\widehat{\tau}_1\colon \widehat{\bX}_1\to \U^\pm(\H)$  is the $\U(1)$-principal bundle with $C_2$-action determined by the pullback~\eqref{eq:gerbepullbacksFHT}. When $\X=M$ is a manifold, this recovers the definition of twisted KR-theory from \cite[\S6.2]{HMSV}, and when $\X$ is a Lie groupoid it recovers the definition \cite[Definition 7.3.2]{MoutuouThesis}. Forgetting the $C_2$-action recovers the Freed--Hopkins--Teleman model for twisted K-theory of $\bX$~\cite[\S3.4]{FHTI}.  


To translate to the context of \cite{FreedMoore,Gomi}, we must translate from Real stacks to stacks over $\pt\sq C_2$ by taking $C_2$-quotients. However, there is a canonical identification between $C_2$-invariant sections and sections over the $C_2$-quotient,
\beq\label{eq:twistKasmapofstacksC2}
&&F\in \Gamma(\X;\tau^*\Fred_{\cCl_{p,q}}(\H))^{C_2}=\Gamma(\X\sq C_2;\tau^*\Fred_{\cCl_{p,q}}(\H)).\nonumber 
\eeq
One again chooses a groupoid presentation $\bY$ of $\X\sq C_2$ in which $F$ is represented by a functor between topological groupoids. 
%
Using Remark~\ref{rmk:orientifold}, we obtain a $\phi$-twisted graded central extension $\widehat{\bY}$ of $\bY$ where the grading and Real structure are encoded by maps $\epsilon \colon \bY_1\to \Z/2$ and $\phi\colon \bY_1\to C_2$, respectively, see~\eqref{eq:Realgerbes2ways}. 
Over this $\phi$-twisted graded central extension $\widehat{\bY}$, we obtain a Clifford bundle determined by a map $\rho\colon \bY_1^\tau\to \U^\pm_{\cCl_{p,q}}(\cH)\rtimes C_2$, where $\phi$ agrees with the postcomposition with the projection to $C_2$, acting by the Real structure on $\cH$. The compact open topology on $\U^\pm_{\cCl_{p,q}}(\cH)$ allows us to identify $\rho$ with a morphism of (trivial) bundles over $\bY_1^\tau$ with fiber $\cH$, which in turns gives a $\U(1)$-invariant map $\rho\colon \bY_1^\tau\times \cH\to \cH$ for the diagonal $\U(1)$-action on $\bY_1^\tau$ and the action on $\cH$ by scalars. Taking the quotient by this (free) $\U(1)$-action gives a map 
\beq\label{eq:Gomicocycle}
{}^\phi \cL\otimes_\C \cH \simeq \bY_1^\tau\times_{\U(1)} \cH\to \cH
\eeq
where $\cL=\bY_1^\tau\times_{\U(1)}\C$ is the complex line bundle associated with the $\U(1)$-bundle $\bY_1^\tau\to \bY_1$. The value of $\phi\colon \bY_1\to C_2$ determines whether~\eqref{eq:Gomicocycle} is $\C$-linear or $\C$-antilinear, and the value $\epsilon\colon \bY_1\to \Z/2$ determines whether the line $\cL$ is even or odd. The isomorphism~\eqref{eq:Gomicocycle} determines a homeomorphism on spaces of Fredholm operators $s^*\Fred(\cH)\to t^*\Fred(\cH)$ over $\bY_1$ determined by the assignment $F\mapsto \id_{\cL}\otimes F$. This homemorphism is the cocycle defining Gomi's Fredholm family associated to the twisting~$\tau$ \cite[Definition~3.1]{Gomi}, see also the sketch in \cite[Remark 7.37]{FreedMoore}.
Hence, the space of $C_2$-invariant sections of the fiber bundle $\tau^*\Fred(\cH)$ is homemorphic to the space of sections computed in the groupoid presentation for $\bY$, and we have shown that the groupoid presentation recovers the space of sections from \cite[Definition~3.1]{Gomi}. 

Finally, we turn to the compatibility with products. This follows from the fact that products in all models are induced by the tensor product of twisted Clifford bundles and the tensor product of Fredholm operators, e.g., see \cite[\S3.5.3]{FHTI}.  
\ep


\begin{rmk}
The previous identification allows us to immediately deduce all the formal properties of twisted KR-theory using the work by previous authors, e.g., homotopy, excision, and exactness axioms \cite[Theorem~3.11]{Gomi} and Bott periodicity \cite[\S3.2]{Gomi}. We also obtain various comparison with other forms of twisted K-theory, see \cite[\S3.4]{Gomi}. 
\end{rmk}

\section{Universal twisted Real Thom classes}

\subsection{The universal twisted Real Thom class}

Endow $\cCl_{p,q}$ with the hermitian inner product for which monomials in the generators $e_i$ are orthonormal. For this Hilbert space structure, $\cCl_{p,q}$ is a unitary $\cCl_{p,q}$-$\cCl_{p,q}$ bimodule. 

\begin{defn} \label{defn:spinorrep}
For a fixed $(p,q)\in \N\times \N$, the \emph{spinor representation}, denoted $\bS_{p,q}$, is the Real, graded, unitary $\Pin^c_{p,q}$-representation determined by applying Lemma~\ref{lem:pinorrep} to $\cCl_{p,q}$ as a Real unitary left module over itself. 
\end{defn}

Since the left and right $\cCl_{p,q}$-actions on $\cCl_{p,q}$ commute, $\bS_{p,q}$ has a right $\cCl_{p,q}$-action commuting with the projective $\O_{p,q}$-representation. Hence, in this case the homomorphism~\eqref{eq:projrepOpq} factors through the Clifford-linear unitary group
\beq\label{eq:Thomtwist}
\O_{p,q}\to \PU^\pm_{\cCl_{p,q}}(\bS_{p,q}).
\eeq

\begin{defn}\label{defn:Thomtwist}
The \emph{Thom twist} is the map of stacks
\beq\label{eq:theThomtwist}
\pt\sq \O_{p,q}\to \pt\sq \PU^\pm_{\cCl_{p,q}}(\bS_{p,q})
\eeq
determined by the Real, graded homomorphism~\eqref{eq:Thomtwist}.
\end{defn}

\begin{rmk}\label{rmk:stackybS}
We provide an alternative (equivalent) description of the Thom twist following \cite[Remark 3.2.22]{ST04}. Let $\cCl(\R^{p,q})$ denote the Real Clifford algebra $\cCl_{p,q}$ with its natural $\O_{p,q}$-action by algebra automorphisms, and let $\cCl_{p,q}$ denote the Clifford algebra with the trivial $\O_{p,q}$-action. We regard these algebras with $\O_{p,q}$-actions as algebra bundles over the stack~$\pt\sq \O_{p,q}$. One can ask for a bundle of invertible $\cCl(\R^{p,q})$-$\cCl_{p,q}$-bimodules over $\pt\sq \O_{p,q}$. However, such a bimodule exists only after pulling the algebra bundles back along~$\pt\sq \Pin_{p,q}^c\to \pt\sq \O_{p,q}$, where the $\Pin_{p,q}^c$-representation $\bS_{p,q}$ from Definition~\ref{defn:spinorrep} furnishes the desired invertible bimodule bundle on $\pt\sq \Pin_{p,q}^c$. More formally, $\pt\sq \Pin_{p,q}^c$ is the terminal object in the bicategory of stacks over $\pt\sq \O_{p,q}$ on which the algebra bundles $\cCl(\R^{p,q})$ and $\cCl_{p,q}$ on $\pt\sq \O_{p,q}$ are (Morita) equivalent. 
\end{rmk}


For $\Pin^c_{p,q}$ acting on $\R^{p,q}$ through the homomorphism $\Pin^c_{p,q}\to \O_{p,q}$, there is a $\Pin_{p,q}^c$-equivariant map valued in $\cCl_{p,q}$-linear endomorphisms
\beq\label{eq:Thomfamofops}
\R^{p,q}\subset \cCl_{p,q} \to \End_{\cCl_{p,q}}(\bS_{p,q})\qquad x\mapsto \cl_x
\eeq
sending a vector $x\in\R^{p,q}$ to the endomorphism of $\bS_{p,q}$ given by the left Clifford action by~$x$. By construction, the family of operators $\cl_x$ commutes with the right Clifford action. The claimed $\Pin_{p,q}^c$-equivariance of~\eqref{eq:Thomfamofops} follows from the fact that $\Pin^c_{p,q}$ acts on $\cl_x$ through the twisted adjoint representation~\eqref{eq:pinpinc},
\beq\label{eq:twistedadjoint}
\tilde{g}\cdot {\rm cl}_x=\det(g)\varphi_{\tilde{g}}\cl_x\varphi_{\tilde{g}}^{-1}=\cl_{gx},\qquad \tilde{g}\in \Pin^c_{p,q}, \ x\in \R^{p,q}
\eeq
where $g\in \O_{p,q}$ is the image of $\tilde{g}\in \Pin^c_{p,q}$ under ${\rm Ad}\colon \Pin^c_{p,q}\to \O_{p,q}$, and the sign $\det(g)\in \{\pm 1\}$ tracks the grading of $\tilde{g}$, as discussed before~\eqref{eq:Pingrading}.

\begin{lem}
The maps~\eqref{eq:Thomtwist} and~\eqref{eq:Thomfamofops} give a map of Real graded  gerbes
\beq\label{eq:preThom}
&&
\begin{tikzpicture}[baseline=(basepoint)];
\node (A) at (0,0) {$\R^{p,q}\sq \Pin^c_{p,q}$};
\node (B) at (5,0) {$\Fred_{\cCl_{p,q}}(\bS_{p,q})\sq \U_{\cCl_{p,q}}^\pm(\bS_{p,q})$};
\node (C) at (0,-1.5) {$\R^{p,q}\sq \O_{p,q}$};
\node (D) at (5,-1.5) {$\pt \sq \PU_{\cCl_{p,q}}^\pm(\bS_{p,q})$};
\draw[->] (B) to  (D);
\draw[->] (C) to (D);
\draw[->] (A) to (B);
\draw[->] (A) to (C);
\path (0,-.75) coordinate (basepoint);
\end{tikzpicture}
\eeq
covering the Real map of stacks $\R^{p,q}\sq \O_{p,q}\to \pt \sq \PU_{\cCl_{p,q}}^\pm(\bS_{p,q})$ given by the composition of the projection $\R^{p,q}\sq \O_{p,q}\twoheadrightarrow \pt\sq \O_{p,q}$ and the Thom twist~\eqref{eq:theThomtwist}.
\end{lem}
\bp
First we note that~\eqref{eq:Thomfamofops} lands in the subspace of odd self-adjoint operators, and (since $\bS_{p,q}$ is finite rank) this is identified with the space of Fredholm operators $\Fred_{\cCl_{p,q}}(\bS_{p,q})\subset \End_{\cCl_{p,q}}(\bS_{p,q})$. By~\eqref{eq:twistedadjoint}, the map~\eqref{eq:Thomfamofops} is equivariant relative to the homomorphism~\eqref{eq:Thomtwist} and hence defines the upper horizontal arrow in~\eqref{eq:preThom} as a map of stacks. 
Compatibility with Real and graded structures is clear. To see that this Real map of stacks is a map of Real graded gerbes, observe that the homomoprhism $\Pin^c_{p,q}\to \U^\pm(\bS_{p,q})$ is $\U(1)$-equivariant for the action on the center of the source and scalars on the target. We conclude that the upper arrow in~\eqref{eq:preThom} is $\pt\sq \U(1)$-equivariant and therefore a map of Real graded gerbes. 
\ep

Let $S_n=\bS_{n,0}^{C_2}\subset \bS_{n,0}$ denote the real subspace of the spinor representation. By the same argument as in Lemma~\ref{lem:pinorrep},  $S_n$ is naturally of $\Pin_n$-representation in right $\Cl_n$-modules. The following is straightforward. 

\begin{cor}
There are maps of graded groupoids,
$$
\R^n\sq \Pin^c_n\to \Fred_{\cCl_n}(\bS_n)\sq \U^\pm_{\cCl_n}(\bS_n)\qquad \R^n\sq \Pin_n\to \Fred_{\Cl_n}(S_n)\sq \O^\pm_{\Cl_n}(S_n),
$$
where the first map corresponds to the underlying map of stacks~\eqref{eq:preThom} (forgetting the Real structures, so $n=p+q$) and the second comes from applying~\eqref{eq:preThom} to $\R^n$ with the trivial Real structure (so $(p,q)=(n,0)$). 
\end{cor}

\begin{defn}\label{defn:Thomclass}
Define the $(p,q)^{\rm th}$ \emph{universal $\KR$-Thom cocycle} as the Real map 
\beq\label{eq:Thom}
\widehat{\Th}_{p,q}\colon \R^{p,q}\sq \O_{p,q}\to \Fred_{\cCl_{p,q}}(\bS_{p,q})\sq \PU_{\cCl_{p,q}}^\pm(\bS_{p,q})
\eeq
associated to~\eqref{eq:preThom}, with underlying \emph{Thom class} $\Th_{p,q}\in \KR^{\SW_{p,q}}(\R^{p,q}\sq \O_{p,q})$. 

Forgetting the $C_2$-actions, the map~\eqref{eq:Thom} defines the $(p+q)^{\rm th}$ \emph{universal $\KU$-Thom cocycle} with underlying class $\Th_{p+q}\in \KU^{\SW_{p,q}}(\R^{p+q}\sq \O_{p+q})$. Restricting to case of a trivial $C_2$-action on $\R^n$, 
\beq\label{eq:ThomKO}
\widehat{\Th}_{n}\colon \R^{n}\sq \O_{n}\to \Fred_{\Cl_{n}}(S_n)\sq \PO^\pm_{\Cl_n}(S_n)
\eeq
defines the \emph{universal $\KO$-Thom cocycle} with underlying class $\Th_n\in \KO^{\SW_n}(\R^n\sq \O_n)$. In~\eqref{eq:ThomKO}, $S_n=\bS_{n}^{C_2}=\Cl_n$ is the Pin representation built from the real Clifford algebra as a bimodule over itself. 
\end{defn}


\begin{rmk}
Following Remark~\ref{rmk:compactsupport}, we observe that the map~\eqref{eq:Thom} satisfies the required compact support condition: Clifford multiplication by $x\in \R^{p,q}$ determines an invertible operator for $x\ne 0$. 
\end{rmk}

\begin{rmk}\label{eq:universalThomHilb}
The twisted Hilbert bundle associated to the twist~\eqref{eq:Thomtwist} is not absorbing. We describe the absorbing twist on $\pt\sq \O_n$ for the Real graded gerbe classified by~\eqref{eq:Thomtwist} presently. Consider $ L^2(\O_{p,q})\otimes \ell^2\otimes \bS_{p,q}$ as a representation of $\Pin^c_{p,q}$ for the action on $\bS_{p,q}$ from Definition~\ref{defn:spinorrep} and action on $L^2(\O_{p,q})$ through the homomorphism $\Pin^c_{p,q}\to \O_{p,q}$. By the Peter--Weyl Theorem, we have an embedding $ L^2(\O_{p,q})\otimes \ell^2\otimes \bS_{p,q}\subset L^2(\Pin_{p,q})\otimes \ell^2$. Indeed, this embedding is as the weight~1 subspace $\cH(1)\subset L^2(\Pin^c_{p,q})\otimes \ell^2$ for the $\U(1)$-action by the subgroup $\U(1)<\Pin^c_{p,q}$: $\U(1)$ acts on $\bS_{p,q}$ by scalars and acts trivially on $L^2(\O_{p,q})$. By \cite[Lemma~3.11]{FHTI}, this gives the absorbing twisted Hilbert bundle for the Real graded gerbe classified by~\eqref{eq:Thomtwist}. 
\end{rmk}


\subsection{The Real twisted Thom class of a vector bundle}

Following Atiyah~\cite{Atiyahreal}, a \emph{Real vector bundle} is a Real map of stacks $\V\to \X$ for a complex vector bundle $\V\to \X$ where the involution defining the Real structure acts complex antilinearly on $\V$. Maps between Real vector bundles are defined in the obvious way. We construct Thom isomorphisms for the following subcategory of Real vector bundles.

\begin{lem}\label{lem:RealVB}
A Real map $\X\to \pt\sq \O_{p,q}$ determines a Real vector bundle over $\X$. A Real map $\X\to \Y$ over $\pt\sq \O_{p,q}$ determines a morphism between the corresponding Real vector bundles on~$\X$ and $\Y$. 
\end{lem}

\bp
The map of Real stacks $\R^{p,q}\sq \O_{p,q}\to \pt\sq \O_{p,q}$ has the structure of a $C_2$-equivariant real vector bundle. Taking the complexification $\C^{p,q}:=\R^{p,q}\otimes \C$ and extending the involution on $\R^{p,q}$ complex antilinearly to $\C^{p,q}$ defines a Real vector bundle $\C^{p,q} \sq \O_{p,q}\to \pt\sq \O_{p.q}$, which pulls back to a Real vector bundle over $\X$ along a Real map $\X\to \pt\sq \O_{p,q}$. The remainder of the statement follows from naturality of pullbacks. 
\ep

Hereafter, we will assume that all Real vector bundle are of the form in Lemma~\ref{lem:RealVB}. 

\begin{rmk}
A (twisted) Thom isomorphism must restrict in a local trivialization of the given vector bundle to the suspension isomorphism. By Remark~\ref{rmk:suspension} the Real vector bundles in Lemma~\ref{lem:RealVB} comprise all bundles for which it makes sense to ask for a (twisted) $\KR$-Thom isomorphism. Note that this includes the class of Real vector bundles that come from the complexification of a $C_2$-equivariant real vector bundle, e.g., $\R$-vector bundles in Real stacks. 
\end{rmk}

\begin{rmk}
Another class of Real vector bundles included in Lemma~\ref{lem:RealVB} are those that pull back $\C^n\sq \U_n\to \pt\sq \U_n$ for the Real structure on $\C^n\sq \U_n$ from Example~\ref{ex:CnUn}. Indeed, these are a specialization of bundles classified by $\X\to \pt\sq \O_{n,n}$, with ensuing Thom isomorphism related to the complex (or MU) orientation of K-theory. We refer to~\cite{ZachYigal} for a detailed comparison between this complex orientation and the Real Spin$^c$ orientation.  
\end{rmk}

\begin{defn}\label{defn:spinstructure}
Given a Real vector bundle $\V$ classified by $\X\to \pt\sq \O_{p,q}$, the category of \emph{Real Spin$^c$ structures} is the category of Real lifts
\beq
\begin{tikzpicture}[baseline=(basepoint)];
\node (B) at (3,0) {$\pt\sq \Spin^c_{p,q}$};
\node (C) at (0,-1.5) {$\X$};
\node (D) at (3,-1.5) {$\pt\sq \O_{p,q}.$};
\draw[->] (B) to  (D);
\draw[->] (C) to node [below] {$\V$} (D);
\draw[->,dashed] (C) to (B);
\path (0,-.75) coordinate (basepoint);
\end{tikzpicture}
\eeq
The category of \emph{Spin$^c$ structures} on $\V$ is  the category of lifts without Real structures. For a real vector bundle $\V$ over a stack $\X$, the category of spin structures is the category of Real lifts for the trivial $C_2$-action on $\X$. 

\end{defn}
\begin{rmk}
When $\X$ has the trivial Real structure, a Real vector bundle $\V\to \X$ is the same data as a real vector bundle $\V^{C_2}\to \X$ given by a map $\X\to \pt\sq \O_n=\pt\sq \O_{n,0}$, and a spin structure is a lift to $\pt\sq \Spin_n=(\pt\sq \Spin^c_n)^{C_2}$. 
\end{rmk}

\begin{defn}
Given a Real vector bundle $\V$ classified by $\X\to \pt\sq \O_{p,q}$,
the \emph{Stiefel--Whitney gerbe} of $\V$ is the pullback $\X^{\SW_\V}$ in Real graded  gerbes
\beq\label{eq:pinpincstru}
&&\begin{tikzpicture}[baseline=(basepoint)];
\node (A) at (0,0) {$\X^{\SW_V}$};
\node (AA) at (0,1.5) {$\pt\sq \U(1)$};
\node (BB) at (3,1.5) {$\pt\sq \U(1)$};
\node (B) at (3,0) {$\pt\sq \Pin^c_{p,q}$};
\node (C) at (0,-1.5) {$\X$};
\node (D) at (3,-1.5) {$\pt\sq \O_{p,q}$};
\node (DD) at (6,-1.5) {$\pt\sq(\Z/2).$};
\node (P) at (.7,-.4) {\scalebox{1.5}{$\lrcorner$}};
\draw[double equal sign distance] (AA) to  (BB);
\draw[->] (AA) to  (A);
\draw[->] (BB) to  (B);
\draw[->] (A) to  (B);
\draw[->] (B) to  (D);
\draw[->] (A) to  (C);
\draw[->] (C) to node [below] {$\V$} (D);
\draw[->] (D) to node [below] {$\det$} (DD);
\path (0,0) coordinate (basepoint);
\end{tikzpicture}
\eeq
\end{defn}

\begin{rmk}
Following Remark~\ref{rmk:cohomology}, the Stiefel--Whitney gerbe is classified by 
\beq\label{eq:DDclass}
(W_3(\V),w_1(\V))\in \HH^3(\X;\underline{\Z})\times \HH^1(\X;\Z/2)
\eeq
where $\underline{\Z}$ is the $C_2$-equivariant local system on~$\X$ for the Real structure. 
\end{rmk}

We refer to Definition~\ref{defn:trivializationgerbe} for the notion of a trivialization of a Real graded gerbe.

\begin{lem}\label{lem:SWgerbe}
The category of Real Spin$^c$ structures is equivalent to the category of trivializations of the Stiefel--Whitney gerbe.
\end{lem}
\begin{proof}
The first piece of data in a trivialization of the Stiefel--Whitney gerbe is a trivialization of the grading $\X\to \pt\sq \O_{p,q}\to \pt\sq (\Z/2)$, i.e., an orientation on~$\V$. With an orientation fixed, we obtain a factorization of the classifying map through $\pt\sq \SO_{p,q}$. We have the pullback square in Real Lie groups, 
\beq\nonumber
&&
\begin{tikzpicture}[baseline=(basepoint)];
\node (A) at (0,0) {$\Spin^c_{p,q}$};
\node (B) at (3,0) {$\Pin^c_{p,q}$};
\node (C) at (0,-1.5) {$\SO_{p,q}$};
\node (D) at (3,-1.5) {$\O_{p,q}.$};
\node (P) at (.7,-.4) {\scalebox{1.5}{$\lrcorner$}};
\draw[->] (B) to  (D);
\draw[->] (C) to (D);
\draw[->] (A) to (B);
\draw[->] (A) to (C);
\path (0,-.75) coordinate (basepoint);
\end{tikzpicture}
\eeq
This implies that the Stiefel--Whitney gerbe for an oriented vector bundle is isomorphic to the pullback of the Real gerbe $\pt\sq \Spin^c_{p,q}\to \pt\sq \SO_{p,q}$ with the trivial grading. We conclude that a trivialization of the Stiefel--Whitney gerbe determines a Real Spin$^c$ structure on~$\V$. It is straightforward to check that an isomorphism between trivializations of the Stiefel--Whitney gerbe corresponds to isomorphic lifts to $\pt\sq \Spin^c_{p,q}$. 
\ep

Consider the pullback squares 
\beq
\begin{tikzpicture}[baseline=(basepoint)];
\node (A) at (0,0) {$\V$};
\node (B) at (3,0) {$\R^{p,q}\sq \O_{p,q}$};
\node (C) at (0,-1.5) {$\X$};
\node (D) at (3,-1.5) {$\pt\sq \O_{p,q}$};
\node (P) at (.7,-.4) {\scalebox{1.5}{$\lrcorner$}};
\draw[->] (B) to  (D);
\draw[->] (C) to (D);
\draw[->] (A) to (B);
\draw[->] (A) to (C);
\path (0,-.75) coordinate (basepoint);
\end{tikzpicture}\qquad\qquad
\begin{tikzpicture}[baseline=(basepoint)];
\node (A) at (0,0) {$\V^{\SW_\V}$};
\node (B) at (3,0) {$\X^{\SW_V}$};
\node (C) at (0,-1.5) {$\V$};
\node (D) at (3,-1.5) {$\X.$};
\node (P) at (.7,-.4) {\scalebox{1.5}{$\lrcorner$}};
\draw[->] (B) to  (D);
\draw[->] (C) to (D);
\draw[->] (A) to (B);
\draw[->] (A) to (C);
\draw[->, bend left, dashed] (C) to (A);
\draw[->, bend right, dashed] (D) to (B);
\path (0,-.75) coordinate (basepoint);
\end{tikzpicture}
\eeq
The first square gives a map of stacks $\V\to \R^{p,q}\sq \O_{p,q}$, and the second shows that a trivializing section of the Stiefel--Whtiney gerbe $\X^{\SW_\V}$ determines a trivializing section of the Real graded gerbe~$\V^{\SW_\V}\to \V$. 

\begin{defn} \label{defn:ThomofV}
Given a Real vector bundle $\V\to \X$, the \emph{Real twisted $\KR$-Thom class of $\V$} is the image of the universal Thom class under the pullback along $\V\to \R^{p,q}\sq \O_{p,q}$,  
\beq
&&\
\KR^{\SW_{p,q}}(\R^{p,q}\sq \O_{p,q})\to \KR^{\SW_{p,q}}(\V),\qquad \Th_{p,q}\mapsto \Th_\V.\nonumber
\eeq
Given a vector bundle $\V\to \X$ with a chosen Real Spin$^c$ structure, the (untwisted) \emph{Real $\KR$-Thom class of $\V$} is the image of the universal Thom class under the composition 
\beq
&&\KR^{\SW_{p,q}}(\R^{p,q}\sq \O_{p,q})\to \KR^{\SW_{p,q}}(\V)\simeq \KR^{p,q}(\V)\nonumber
\eeq
where the isomorphism is the pullback along the map of stacks $\V\to \V^{\SW_\V}$ gotten from a trivializing section of the Stiefel--Whitney gerbe. 

When $\X$ is not given a Real structure, respectively, is equipped with the trivial Real structure, the \emph{twisted $\KU$-}, respectively, \emph{$\KO$-Thom class of $\V$} is the image of the universal Thom class under
\beq\nonumber
\Th_n\in \KU^{\SW_{n}}(\R^{n}\sq \O_{n})\to \KU^{\SW_{n}}(\V)\ni \Th_\V,\qquad \Th_n\in \KO^{\SW_{n}}(\R^{n}\sq \O_{n})\to \KO^{\SW_{n}}(\V)\ni \Th_\V
\eeq
and a Spin$^c$, respectively, Spin structure provides an untwisted Thom class under the composition 
\beq
&&\KU^{\SW_{n}}(\R^{n}\sq \O_{n})\to \KU^{\SW_{n}}(\V)\simeq \KU^{n}(\V),\qquad\KO^{\SW_{n}}(\R^{n}\sq \O_{n})\to \KO^{\SW_{n}}(\V)\simeq \KO^{n}(\V)
\eeq
in the respective cases. When the context is clear, we drop the $\KR$, $\KU$, or $\KO$ descriptor, referring simply to a \emph{Thom class}. 

\end{defn}

%

\subsection{Comparison with previously constructed $\KO$- and $\KU$-Thom cocycles}\label{sec:Thomclassofvb}

In this subsection, we assume the topological stack~$\X$ is either not equipped with a Real structure (for the $\KU$ case) or has the trivial Real structure (for the~$\KO$ case). 
Furthermore, we specialize to vector bundles over global quotients, $\V\to X\sq G$. Such a bundle is equivalent to a $G$-equivariant vector bundle~$V\to X$, i.e., $\V=V\sq G\to X\sq G=\X$.

\begin{prop}\label{prop:standardThom}
Given a $G$-equivariant vector bundle $V\to X$, the twisted $\KU$-Thom class from Definition~\ref{defn:ThomofV} agrees with the twisted Thom class from \cite[\S3.6]{FHTI}. When~$V$ has a choice of spin structure, the untwisted $\KO$-Thom class from Definition~\ref{defn:ThomofV} agrees with the Thom class from \cite[Remark~3.2.22]{ST04}. These specialize to the $\KU$ Atiyah--Bott--Shapiro Thom class~\cite{ABS} $n\equiv 0$ mod 2 and the $\KO$ class in degrees $n\equiv 0$ mod~8. 
\end{prop}

\begin{proof}[Proof of Proposition~\ref{prop:standardThom}]
We begin by verifying the $\KU$-case, i.e., without any Real structures present. A map $X\sq G\to \pt\sq \O_n$ of stacks classifying a $G$-equivariant vector bundle $V\to X$ is equivalent to a zigzag of continuous functors between topological groupoids
\beq\label{eq:framebundle}
X\sq G\xleftarrow{\sim} \Fr(V)\sq (\O_n\times G)\to \pt \sq\O_n
\eeq
where $\Fr(V)$ is the (orthonormal) frame bundle of~$V$, and the second arrow is the projection. The pullback of the universal bundle $\R^n\sq \O_n\to \pt\sq \O_n$ along~\eqref{eq:framebundle} gives a zigzag
\beq\label{eq:framebundle2}
V\sq G\xleftarrow{\sim} (\Fr(V)\times \R^n)\sq (\O_n\times G)\to \R^n\sq \O_n.
\eeq
From~\eqref{eq:framebundle}, the pullback of the Thom twist gives the twisted Hilbert bundle 
\beq\label{eq:Eulerbundle}
(\Fr(V)\times \O_n\times \bS_n)\sq (\Pin^c_n\times G)\to \Fr(V)\sq (\O_n\times G)\simeq X\sq G
\eeq
over $X\sq G$. The Thom class is represented by the family of operators gotten by composing~\eqref{eq:preThom} with~\eqref{eq:framebundle2}. This is precisely the twisted Thom class constructed in~\cite[\S3.6]{FHTI}, see also~\cite[1.92 and 1.161]{Vienna}. Assuming that $V$ admits a Spin$^c$ structure, Lemma~\ref{lem:SWgerbe} provides 
\beq\label{eq:trivofSWquotient}
\Fr(V)\sq (\O_n\times G) \xleftarrow{\sim} \Spin^c(V)\sq (\Spin^c_n\times G)\to \Fr(V) \sq (\Pin^c_n\times G)
\eeq
where $\Spin^c(V)$ is the principal $\Spin^c$-bundle for $V$, and the zigzag~\eqref{eq:trivofSWquotient} determines a map of stacks that is a trivializing section of the Stiefel--Whitney gerbe. Pulling back the Thom twist along~\eqref{eq:trivofSWquotient} we obtain the Hilbert bundle~\eqref{eq:Eulerbundle} 
\beq\label{eq:STEuler}
&&(\Spin^c(V)\times \bS_n)\sq \Spin^c_n \simeq \Spin^c(V)\times_{\Spin^c_n} \bS_n\simeq \Spin^c(V)\times_{\Spin^c_n} \cCl_{n}
\eeq
which we identify with the complex version of Stolz and Teichner's Euler class in $\K$-theory. 

Repeating the above analysis for a trivial Real structure on $\X$ has the effect of replacing the  $\Spin^c$-groups by $\Spin$-groups and the complex Clifford algebras $\cCl_n$ by their real counterparts $\Cl_n$. In particular,~\eqref{eq:STEuler} determines a bundle of real, invertible $\Cl(V)$-$\Cl_{n}$ bimodules  over~$X$. The Thom cocycle pulled back along~\eqref{eq:trivofSWquotient} gives the $V$-family of operators that at a point $v\in V$ is (left) Clifford multiplication by $v\in V\in \Cl(V)$. This is precisely Stolz and Teichner's description of the Thom cocycle \cite[Remark~3.2.22]{ST04}. 

Finally, to compare with~\cite{ABS} in the complex case, one applies a fiberwise Morita equivalence to \eqref{eq:STEuler} between $\cCl_{n}$ and~$\C$ in degrees $0$ mod $2$. In the real case, one applies the Morita equivalence between $\Cl_{n}$ and~$\R$ in degrees $0$ mod $8$. These Morita equivalences afford isomorphisms on spaces of Fredholm operators (see Remark~\ref{rmk:RealMorita}), extracting the Atiyah--Bott--Shapiro Thom classes.
\ep

\subsection{The twisted Thom isomorphism}

\begin{prop}[Universal Thom isomorphism]\label{prop:Thomiso}
The external product with the universal twisted Thom class
\beq\label{eq:Thomisoeq}
\KR^\alpha(\pt\sq \O_{p,q})\xrightarrow{\smallsmile \Th_{p,q}} \KR^{\alpha\otimes\SW_{p,q}}(\R^{p,q}\sq \O_{p,q}),
\eeq
determines an isomorphism of abelian groups.
\end{prop}
\bp
This is a special case of \cite[Theorem~3.19]{Gomi}, following the argument in the case of twisted $\KU$-theory in \cite[pg 32-33]{FHTI}. Both for completeness and to emphasize the simplifications in this universal case we provide a detailed sketch below. 

Using the universal Thom-twisted Hilbert bundle from Remark~\ref{eq:universalThomHilb}, the $\O_{p,q}$-equivariant Atiyah--Singer map provides the weak equivalence on global sections
\beq
\KR^{\SW_{p,q}}(\R^{p,q}\sq \O_{p,q})&=&\pi_0\Gamma_c(\R^{p,q}\sq \Pin^c_{p,q},\Fred^{p,q}(L^2\O_{p,q}\otimes \ell^2\otimes \bS_{p,q}))^{C_2}\nonumber\\
&\xleftarrow{\sim}&\pi_0\Gamma(\pt\sq \Pin_{p,q}^c,\Fred^{p,q}_{\cCl(\R^{p,q})^\op}(L^2\O_{p,q}\otimes \ell^2\otimes \bS_{p,q}))^{C_2}\label{eq:ASmap}
\eeq
where $\cCl(\R^{p,q})$ denotes the Clifford algebra with its $\O_{p,q}$-action. In more detail, the source of the Atiyah--Singer map consists of Fredholm operators that are linear for the (right) $\cCl_{p,q}\otimes \cCl(\R^{p,q})^\op$-action on $\bS_{p,q}$. By its definition, the image of the Atiyah--Singer map are classes given by the external product of the Thom class~\eqref{eq:Thom} and Fredholm operators on $L^2\O_{p,q}\otimes \ell^2$. 
We refer to~\cite[\S3.6]{FHTI} and~\cite[Theorem 3.19 and Lemma B.12]{Gomi} for further discussion of the $\O_{p,q}$-equivariant Atiyah--Singer map.

Next, use the characterization from Remark~\ref{rmk:stackybS} of $\bS_{p,q}$ as an invertible bimodule between $\cCl(\R^{p,q})$ and $\cCl_{p,q}$, i.e., the Clifford algebra $\cCl(\R^{p,q})$ with its natural $O_{p,q}$-action and the Clifford algeba $\cCl_{p,q}$ with the trivial $\O_{p,q}$-action. This invertible bimodule is equivalently an invertible $\C$-$\cCl_{p,q}\otimes \cCl(\R^{p,q})^\op$ bimodule, and hence provides a Real $\Pin_{p,q}^c$-equivariant homeomorphism 
\beq
&&\Fred^{p,q}_{\cCl(\R^{p,q})^\op}(L^2\O_{p,q}\otimes \ell^2\otimes \bS_{p,q})
\simeq \Fred^0(L^2\O_{p,q}\otimes \ell^2).\label{eq:Mortiamap}
\eeq
The above determines an isomorphism on homotopy classes of Real sections over $\pt\sq \Pin^c_{p,q}$. To rephrase~\eqref{eq:Mortiamap} in the language of stacks (following Remark~\ref{rmk:stackybS}), the $\Pin^c_{p,q}$-action on $\bS_{p,q}= (\bS_{p,q})_{\cCl_{p,q}\otimes \cCl(\R^{p,q})}$ determines a bundle of invertible bimodules over $\pt\sq \Pin_{p,q}^c$. This invertible bimodule bundle then identifies the relevant spaces of Real sections over $\pt\sq \Pin_{p,q}^c$. 

Finally, we observe that in the target of~\eqref{eq:Mortiamap}, the $\Pin^c_{p,q}$-action factors through the $\O_{p,q}$-action. Hence, the bundle of Fredholm operators pulls back from $\pt\sq \O_{p,q}$, identifying an (a priori) twisted $\KR$-group with an untwisted one. From this we see
\beq
\KR^0(\pt\sq \O_{p,q})&=&\pi_0\Gamma(\pt\sq \O_{p,q},\Fred^0(L^2\O_{p,q}\otimes \ell^2))^{C_2}\nonumber\\
&\simeq&\pi_0\Gamma(\pt\sq \Pin^c_{p,q},\Fred^0(L^2\O_{p,q}\otimes \ell^2))^{C_2}.\nonumber
\eeq
This identification together with~\eqref{eq:ASmap} and~\eqref{eq:Mortiamap} provides the universal Thom isomorphism with $\alpha=0$. For general $\alpha$, one applies the same argument as above, but starting with a tensor product of a universal Thom-twisted Hilbert bundle and a universal $\alpha$-twisted Hilbert bundle. 
\ep

\begin{rmk}\label{rmk:cvs}
In the following, $\KR^{p^*\alpha\otimes \SW_\V}(\V)_{cv}$ denotes twisted $\KR$-theory with \emph{compact vertical support} along the fibers of the projection $\V\to \X$. If one assumes $\X$ is compact (meaning its coarse moduli space is compact) then $\KR^{p^*\alpha\otimes\tau_\V}(\V)_{cv}=\KR^{p^*\alpha\otimes\tau_\V}(\V)$, where we remind our convention for compactly supported K-theory from Remark~\ref{rmk:compactsupport}. 

\end{rmk}

\begin{cor}[The Thom isomorphism]\label{cor:Thomiso}
For a Real vector bundle $p\colon \V\to \X$ classified by a Real map $\X\to \pt\sq \O_{p,q}$, the pullback of the universal Real twisted Thom isomorphism~\eqref{eq:Thomisoeq} determines the twisted Thom isomorphism
\beq\label{eq:nonuniThomiso}
\KR^\alpha(\X)\xrightarrow{\sim} \KR^{p^*\alpha\otimes\tau_\V}(\V)_{cv}. 
\eeq
Analogous statements hold for the twisted $\KU$- and $\KO$-Thom isomorphisms. 
\end{cor}
\bp
The argument follows \cite[Theorem~3.19]{Gomi} and \cite[pg 32-33]{FHTI} even more closely than the previous result, so we will be brief. It suffices to prove the statement in the case of a Real global quotient $\X=X\sq G$ and $\V=V\sq G$ determined by a $G$-equivariant vector bundle $V\to X$. A map $X\sq G\to \pt\sq \O_{p,q}$ of Real stacks classifying a Real $G$-equivariant vector bundle $V\to X$ is equivalent to a $C_2$-equivariant zigzag of continuous functors between topological groupoids
\beq\label{eq:Realframebundle}
X\sq G\xleftarrow{\sim} \Fr(V)\sq (\O_{p,   q}\times G)\to \pt \sq\O_{p,q}
\eeq
where $\Fr(V)$ is the (orthonormal) frame bundle of~$V$, and the second arrow is the projection. The pullback of the universal bundle $\R^{p,q}\sq \O_{p,q}\to \pt\sq \O_{p,q}$ along~\eqref{eq:framebundle} gives a zigzag of $C_2$-groupoids
\beq\nonumber
&&V\sq G\xleftarrow{\sim} (\Fr(V)\times \R^{p,q})\sq (\O_n\times G)\to \R^{p,q}\sq \O_{p,q}\xrightarrow{\Th_{p,q}} \Fred_{\cCl_{p,q}}(\bS_{p,q})\sq \PU^\pm_{\cCl_{p,q}}(\bS_{p,q})
\eeq
that determines the Thom class of $V\sq G$. 
The pullback of the Atiyah--Singer map~\eqref{eq:ASmap} produces an isomorphism between $\KR$ classes on $V\sq G$, and $\KR$ classes over $X\sq G$ twisted by $\cCl(V)$, the fiberwise Clifford algebra of the vector bundle. These classes have compact vertical support. The pullback of the Morita equivalence replaces the $\cCl(V)$-twisting with the Thom twisting, and hence recovers the claimed twisted Thom isomorphism. 
\ep

\subsection{Twisted pushforwards along an embedding of smooth stacks}\label{sec:toppush}

For the remainder of this section, all stacks will be \emph{proper differentiable stacks}. By the linearization theorem, these are the same as local quotient stacks on the site of smooth manifolds, see Remark~\ref{rmk:linearization}. The bicategory of Real smooth stacks is defined in analogy to Definition~\ref{defn:Realstack}, and the basic results concerning Real topological stacks carry over to Real smooth stacks. We begin by spelling out the topological pushforward in twisted equivariant K-theory, e.g., see~\cite[\S3.6]{FHTI} and~\cite{Carey} for the non-equivariant case.

\begin{defn}[{\cite[\S3.2]{stringtopstacks}}]
A Real map of Real smooth stacks $f\colon \X\to \Y$ is an \emph{embedding} if for one (and hence, any) Real smooth atlas $Y\to \Y$, the pullback
\beq\label{eq:embeddingofstacks}
\begin{tikzpicture}[baseline=(basepoint)];
\node (A) at (0,0) {$X$};
\node (B) at (3.5,0) {$\X$};
\node (C) at (0,-1.25) {$Y$};
\node (D) at (3.5,-1.25) {$\Y$};
\node (P) at (.7,-.4) {\scalebox{1.5}{$\lrcorner$}};
\draw[->] (A) to  (B);
\draw[->] (B) to node [right] {$f$} (D);
\draw[->] (A) to (C);
\draw[->] (C) to  (D);
\path (0,-.75) coordinate (basepoint);
\end{tikzpicture}
\eeq
is equivalent to a smooth $C_2$-manifold $X$, and the map $X\to Y$ is a $C_2$-equivariant embedding of smooth manifolds. The $C_2$-equivariant normal bundle $TY/TX$ to the embedding $X\to Y$ descends to a vector bundle $\nu_f\to \X$,~\cite[Lemma/Definition 4.10]{Heinloth}, which is the \emph{normal bundle} to the embedding~$f\colon \X\to \Y$. 
\end{defn}

After choosing an invariant metric, the normal bundle $\nu_f\to \X$ is classified by a Real map $\X\to \pt\sq \O_{p,q}$. An invariant metric exists via the usual averaging trick due to our assumption that $\X$ and $\Y$ are proper.

\begin{defn}\label{defn:pushembedding}
The \emph{pushforward along an embedding} $f\colon \X\to \Y$ of smooth stacks is the composition 
\beq\label{eq:pushforwardalongembed}
&&\KR^{f^*\sigma\otimes\tau_{\nu_f}^{-1}}(\X)\xrightarrow{\Th_{\nu_f}} \KR^{\sigma|_{\nu_f}}(\nu_f)_{cv}\xrightarrow{{\rm Pontryagin-Thom}} \KR^{\sigma}(\Y)
\eeq
of the twisted Thom isomorphism for the normal bundle $\nu_f\to X$ of the embedding, followed by the Pontryagin--Thom collapse map.
\end{defn}


\begin{rmk}
The Pontryagin--Thom collapse map in~\eqref{eq:pushforwardalongembed} can equivalently be described as extending a class with support in $\nu_f\subset \Y$ by zero, i.e., to a class that restricts to zero on the complement $\Y\setminus \nu_f$ and restricts to the given (compactly supported) class on $\nu_f$. This uses the equivariant tubular neighborhood theorem~\cite{Kankaanrinta} to identify the normal bundle of the embedding with an open substack $\nu_f\subset \Y$. This open substack is characterized by the image of the normal bundle of the embedding of manifolds $X\to Y$ along the atlas~$Y\to \Y$ in~\eqref{eq:embeddingofstacks}.
\end{rmk}



\subsection{Pushforwards along submersions and the twisted equivariant Atiyah--Bott--Shapiro orientation}\label{sec:ABS}
Let $X\to Y$ be a Real $G$-equivariant proper submersion for a compact Real Lie group~$G$. By~\cite{Mostow,Palais}, there exists a $G$-equivariant Real embedding $X\hookrightarrow \R^{p,q}$ in a Real $G$-representation $\rho\colon G\to \O_{p,q}$. Taking the product of the embedding with the submersion gives a Real $G$-equivariant embedding $X\hookrightarrow \R^{p,q}\times Y$ with Real $G$-equivariant normal bundle $\nu\to X\sq G$. Consider the 2-commutative diagram in Real stacks
\beq\label{eq:normaltangent}
\begin{tikzpicture}[baseline=(basepoint)];
\node (A) at (0,0) {$X$};
\node (B) at (3,0) {$X\sq G$};
\node (C) at (0,-1.25) {$\R^{p,q}\times Y$};
\node (D) at (3,-1.25) {$(\R^{p,q}\times Y)\sq G$};
\node (E) at (7.5,-.625) {$\pt\sq (\O_{p,q}\times \O_{r,s})$};
\draw[->] (A) to  (B);
\draw[->,right hook-latex] (B) to (D);
\draw[->,right hook-latex] (A) to (C);
\draw[->] (C) to (D);
\draw[->] (B) to node [above] {$\phantom{A} TX\oplus \nu$} (E);
\draw[->] (D) to node [below] {$\rho\oplus TY$} (E);
\path (0,-.75) coordinate (basepoint);
\end{tikzpicture}
\eeq
gotten from restricting the $G$-equivariant Real tangent bundle of $\R^{p,q}\times Y$ along the embedded submanifold~$X$. The above setup determines various twistings: let $\tau_\nu,\tau_{TX}\in \Twist(X\sq G)$ and $\tau_{TY}\in \Twist(Y\sq G)$ denote the associated Thom twistings of the Real vector bundles, and $\tau_\rho\in \Twist(\pt\sq G)$ denote the pullback of the Thom twist $\SW_{p,q}$ on $\pt\sq \O_{p,q}$ to $\pt \sq G$ along the Real $G$-representation. Throughout the discussion below, we will use the same notation for a twist and its pullback along the obvious map, e.g.,~$\tau_\rho$ is regarded as a twisting for $X\sq G$, $Y\sq G$, and $\R^{p,q}\sq G$ using the maps of these stacks to $\pt\sq G$. For later use, we observe that the diagram~\eqref{eq:normaltangent} yields the isomorphism of twistings
\beq\label{Eq:isooftwistings}
\tau_{TX}\otimes \tau_\nu\simeq \tau_\rho\otimes \tau_{TY}\in \Twist(X\sq G). 
\eeq

For $\sigma$ a twisting on $Y\sq G$, consider the twisted topological pushforward~\eqref{eq:pushforwardalongembed} along the embedding $X\sq G\hookrightarrow (\R^{p,q}\times Y )\sq G$
 \beq\label{eq:ABS1}
 \KR^{\sigma\otimes \tau_{\nu}^{-1}}(X\sq G)\xrightarrow{\pi_!^{\rm top}} \KR^{\sigma}((\R^{p,q}\times Y)\sq G)_{cv}.
 \eeq
 Further consider the twisted Thom isomorphism
 \beq\label{eq:ABS2}
&&\KR^{\sigma}((\R^{p,q}\times Y)\sq G)_{cv}\xleftarrow{\smallsmile \Th_\rho} \KR^{\sigma\otimes\tau_\rho^{-1}}(Y\sq G)\simeq \KR^{\sigma\otimes\tau_{TX}^{-1}\otimes \tau_\nu^{-1}\otimes \tau_{TY}}(Y\sq G)
 \eeq
 for the Real vector bundle $(\R^{p,q}\times Y)\sq G\to Y\sq G$, where the isomorphism on the right uses~\eqref{Eq:isooftwistings}. 
 The twisted Atiyah--Bott--Shapiro map combines~\eqref{eq:ABS1} and~\eqref{eq:ABS2}.


\begin{defn} \label{defn:ABS}
For a Real $G$-equivariant proper submersion $\pi\colon X\to Y$, the \emph{topological pushforward} (also called the
\emph{twisted Atiyah--Bott--Shapiro orientation}) 
is the composition of the twisted topological pushforward and twisted Thom isomorphism, 
\beq\label{Eq:twistedABS}
&&\pi_!^{\rm top}\colon \KR^{\sigma\otimes \tau_\nu^{-1}}(X\sq G)\xrightarrow{\pi_!^{\rm top}}\KR^{\sigma}((\R^{p,q}\times Y)\sq G)_{cv}
\xleftarrow{\smallsmile \Th_\rho} \KR^{\sigma\otimes \tau_\rho^{-1}}(Y\sq G).
\eeq
Forgetting $C_2$-actions gives a $\KU$-variant of~\eqref{Eq:twistedABS}, and taking $X$ with trivial $C_2$-action gives a $\KO$-variant of~\eqref{Eq:twistedABS}. 
\end{defn}

\begin{ex}[Atiyah--Bott--Shapiro map for untwisted source]

Suppose that $\tau_{TX}^{-1}\otimes \tau_{TY}\in \Twist(X\sq G)$ is in the image of the pullback along $X\sq G\to Y\sq G$. From~\eqref{Eq:isooftwistings}, we conclude that $\tau_\nu$ also pulls back from $Y\sq G$ and there is a twisted topological pushforward 
 \beq\nonumber
 \KR^0(X\sq G)\simeq \KR^{\tau_\nu\otimes\tau_\nu^{-1}}(X\sq G)\xrightarrow{\pi_!^{\rm top}} \KR^{\tau_\nu}((\R^{p,q}\times Y)\sq G)
 \eeq
 that specializes~\eqref{eq:ABS1} with $\sigma=\tau_\nu$.  The isomorphism $\tau_\nu \otimes \tau_\rho^{-1}\simeq \tau_{TX}^{-1}\otimes  \tau_{TY}$ from~\eqref{Eq:isooftwistings} 
 further provides a twisted Thom isomorphism that specializes~\eqref{eq:ABS2}
 \beq\nonumber
\KR^{\tau_\nu}((\R^{p,q}\times Y)\sq G)\xleftarrow{\smallsmile \Th_\rho} \KR^{\tau_\nu\otimes\tau_\rho^{-1}}(Y\sq G)\simeq \KR^{\tau_{TX}^{-1}\otimes\tau_{TY}}(Y\sq G).
 \eeq
 Hence, when $\tau_{TX}^{-1}\otimes \tau_{TY}\in \Twist(X\sq G)$ pulls back from $Y\sq G$, we have the twisted Atiyah--Bott--Shapiro map 
\beq\label{eq:ABSuntwistedsource}
&&\KR^0(X\sq G)\xrightarrow{\pi_!^{\rm top}}\KR^{\tau_\nu}((\R^{p,q}\times Y)\sq G)
\xleftarrow{\smallsmile \Th_\rho } \KR^{\tau_{TX}^{-1}\otimes\tau_{TY}}(Y\sq G).
\eeq
In particular, the twisting on the far right does not depend on the embedding~$X\hookrightarrow \R^{p,q}$. 
\end{ex}

Note that a $G$-equivariant Real Spin$^c$ structure on the vertical tangent bundle of the fibration $X\to Y$ is equivalent to a trivialization of the twisting $\tau_{TX}^{-1}\otimes \tau_{TY}$. In the case that such a trivialization is fixed, one recovers the usual $G$-equivariant pushforward from~\eqref{eq:ABSuntwistedsource}. As a mild generalization, suppose that the vertical tangent bundle of the submersion $X\to Y$ is equipped with a Real Spin$^c$ structure in the sense of Definition~\ref{defn:spinstructure}, but the $G$-action is not compatible, i.e., the fibers of $X\to Y$ are compact spin manifold with $G$-action that need not preserve the fiberwise spin structure. Hence in this situation one also obtains a version of~\eqref{eq:ABSuntwistedsource} where the twisting in the target pulls back from~$\pt\sq G$. 

\begin{ex}[Twisted pushforwards to the $\pt\sq G$]\label{ex:ABSBG}
Specializing~\eqref{eq:ABSuntwistedsource} further to the case that $Y=\pt$, the twisting in the target of~\eqref{Eq:twistedABS} is classified by an element of the set $\HH^1(BG;\Z/2)\times \HH^3(BG;\underline{\Z})$. Specifically, such a twist is given by a grading $\epsilon \colon G\to \Z/2$ and Real $\U(1)$-central extension $G^\varphi$ of $G$. This gives the identification
$$
\KR^{\tau_\nu\otimes \tau_\rho^{-1}}(\pt\sq G)\simeq \KR^{\tau_{TX}^{-1}}(\pt\sq G)\simeq \Rep^{\tau_{TX}^{-1}}(G)
$$
between the twisted $\KR$-group in the target of~\eqref{Eq:twistedABS} and the Grothendieck group of Real, $\epsilon$-graded, $\varphi$-projective $G$-representations. This recovers the form~\eqref{Eq:twistedABS0} of the twisted Atiyah--Bott--Shapiro map stated in the introduction. 
\end{ex}

\section{Twisted equivariant Real power operations}

\subsection{Tensor powers of Fredholm operators}

Given a Real, $\Z/2$-graded Hilbert space $\cH$, the (graded) tensor powers $\cH^{\otimes k}$ for $k\in \mathbb{N}$ acquire Real graded  structures. Generalizing~\eqref{eq:tensoronU}, these tensor powers determine commuting squares of Real graded  group homomorphisms
\beq\label{eq:poweronU}
\begin{tikzpicture}[baseline=(basepoint)];
\node (A) at (0,0) {$\U(1)\wr\Sigma_k$};
\node (B) at (3,0) {$\U^\pm(\H)\wr\Sigma_k$};
\node (C) at (0,-1.25) {$\U(1)$};
\node (D) at (3,-1.25) {$\U^\pm(\H^{\otimes k})$};
\draw[->] (A) to  (B);
\draw[->] (B) to node [right] {$\otimes$} (D);
\draw[->] (A) to node [left] {$\times$} (C);
\draw[->] (C) to (D);
\path (0,-.75) coordinate (basepoint);
\end{tikzpicture}
\eeq
where the symmetric group action on $\H^{\otimes k}$ permutes the tensor factors. Combined with the action of $\U^\pm(\H)^{\times k}$ on each factor, we obtain an action of the wreath product $\U^\pm(\H)\wr \Sigma_k$ on $\H^{\otimes k}$. This action is unitary, and so gives the homomorphism denoted $\otimes$ in~\eqref{eq:poweronU}. Since the action is multilinear, it is compatible with the left vertical arrow that multiplies scalars. Passing to projective groups, we obtain Real graded  homomorphisms
\beq\label{eq:protopower}
&& \PU^\pm(\H)\wr \Sigma_k \to \PU^\pm(\H^{\otimes k}).
\eeq
The same construction generalizes to the Clifford linear case, 
\beq
&&\PU^\pm_{\cCl_{p,q}}(\cH)\wr \Sigma_k \to \PU^\pm_{\cCl_{p,q}^{\otimes k}}(\cH^{\otimes k}),\label{eq:protopowerCl}
\eeq
with the modification that the $k^{\rm th}$ power of a $\cCl_{p,q}$-module is a $\cCl_{p,q}^{\otimes k}\simeq \cCl_{kp,kq}$-module.

Generalizing~\eqref{eq:Fredtensor}, the $k^{\rm th}$ tensor power of a Fredholm operator determines a Real map
\beq\label{eq:Fredtensorpower}
  \Fred_{\cCl_{p,q}}(\H)^{\times k}&\to& \Fred_{\cCl_{p,q}^{\otimes k}}(\H^{\otimes k}),
\eeq
that is Real equivariant for the action of the wreath product. Hence, we obtain 2-commuting squares in stacks
\beq\nonumber
\begin{tikzpicture}[baseline=(basepoint)];
\node (A) at (0,0) {$ \Fred_{\cCl_{p,q}}(\H_n)^{\times k}\sq \PU^\pm_{\cCl_{p,q}}(\H)\wr\Sigma_k$};
\node (B) at (6,0) {$\Fred_{\cCl_{p,q}}(\H_n^{\otimes k})\sq \PU^\pm_{\cCl_{p,q}^{\otimes k}}(\H^{\otimes k})$};
\node (C) at (0,-1.25) {$\pt \sq \PU^\pm_{\cCl_{p,q}}(\H)\wr\Sigma_k$};
\node (D) at (6,-1.25) {$\pt \sq \PU^\pm_{\cCl_{p,q}^{\otimes k}}(\H^{\otimes k}).$};
\draw[->] (A) to  (B);
\draw[->] (B) to (D);
\draw[->] (A) to (C);
\draw[->] (C) to (D);
\path (0,-.75) coordinate (basepoint);
\end{tikzpicture}
\eeq
Forgetting the $C_2$-actions gives a diagram in topological stacks, whereas taking $C_2$-fixed points gives an analogous diagram for real Hilbert spaces and projective orthogonal groups. 

\subsection{Constructing the total power operations}\label{sec:powerconstruction}

\begin{defn} The \emph{$k^{\rm th}$ total power operation} of a twist~\eqref{eq:twistingdefn} is the composition
\beq
&&\Pow_k\tau\colon \X^{\times k}\sq \Sigma_k\xrightarrow{\Sym^k\tau} \pt\sq  \PU^\pm_{\cCl_{p,q}}(\H^{\otimes k})\wr \Sigma_k \to \pt\sq \PU^\pm_{\cCl_{p,q}^{\otimes k}}(\H^{\otimes k})\label{eq:kpowertwist2}
\eeq
of the $k^{\rm th}$ symmetric power functor on Real stacks (see~\S\ref{sec:sympowers}) with the map of Real stacks determined by the corresponding homomorphism~\eqref{eq:protopower}. 
\end{defn}

\begin{lem}\label{lem:powmulttwist}
The power operation~\eqref{eq:kpowertwist2} on twists assignment extends to a symmetric monoidal functor
$$
\Pow_k\colon \Twist(\X)\to \Twist(\X^{\times k}\sq \Sigma_k). 
$$
\end{lem}
\bp
The monoidal structure on $\Twist(\X)$ is defined using the tensor product of unitary maps (see~\eqref{eq:tensorproductoftwists} and~\eqref{eq:tensoroftwists}), whereas the power operation is defined using tensor powers. The associativity and symmetry of the tensor product then provides canonical isomorphisms $\Pow_k(\tau\otimes \tau')\simeq \Pow_k(\tau)\otimes \Pow_k(\tau')$ in the data of a symmetric monoidal functor. 
\ep
\begin{defn}\label{defn:powerop}
The \emph{$k^{\rm th}$ total power operation} is the map of abelian groups
\beq\label{eq:poweropdefine}
 \KR^{\tau}(\X)\xrightarrow{\Pow_k}\KR^{\Pow_k\tau}(\X^{\times k}\sq \Sigma_k) 
\eeq
determined by the Real maps of stacks
\beq\label{eq:cocyclepowerop}
&&\X^{\times k}\sq \Sigma_k\xrightarrow{\Sym^k(F)} \Fred_{\cCl_{p,q}}(\cH)^{\times k}\sq \PU^\pm_{\cCl_{p,q}}(\cH)\wr \Sigma_k\to \Fred_{\cCl_{p,q}^{\otimes k}}(\cH^{\otimes k})\sq \PU^\pm_{\cCl_{p,q}^{\otimes k}}(\cH^{\otimes k})
\eeq
where the first arrow uses the description of twisted cocycles from~\eqref{eq:twistKasmapofstacks} and functoriality of~$\Sym^k$, and the second arrow uses that the tensor product map~\eqref{eq:Fredtensorpower} is equivariant relative to the homomorphism~\eqref{eq:protopower}. Forgetting the $C_2$-actions in~\eqref{eq:cocyclepowerop} defines the $\KU$-variant of~\eqref{eq:poweropdefine}. When $\X$ has the trivial $C_2$-action,~\eqref{eq:cocyclepowerop} determines the $\KO$-variant of~\eqref{eq:poweropdefine}. By construction, total power operations are multiplicative but not additive. 
\end{defn}

\begin{rmk}
In good cases, one can also define the $k^{\rm th}$ Adams operation as the composition 
\beq\nonumber
\KR^{\tau}(\X)&\xrightarrow{\Delta^*\circ \Pow_k}& 
\KR^{\tau^{\otimes k}}(\X \times \pt\sq \Sigma_k)\stackrel{?}{\dashrightarrow} \KR^{\tau^{\otimes k}}(\X)\otimes \Rep(\Sigma_k)\to \KR^{\tau^{\otimes k}}(\X) 
\eeq
gotten from the composition of the total power operation, followed by restriction to the diagonal (with trivial $\Sigma_k$-action), 
$$
\X\times \pt\sq \Sigma_k \xhookrightarrow{\Delta} \Sym^k(\X)
$$
and the restriction to the representation of $\Sigma_k$ whose character is the indicator function on the long cycle $(12\cdots k)\in \Sigma_k$. However, in the twisted setting the arrow labeled by a question mark need not exist as the $\Sigma_k$-equivariant structure on $\tau^{\otimes k}$ can be nontrivial. There are twisted generalizations of the Adams operations whose existence depends upon properties of various tensor powers of the twisting, see~\cite[\S10]{AtiyahSegaltwistco} and \cite[\S8]{DonovanKaroubi}. 
\end{rmk}

\subsection{Properties of the power operation}

\begin{prop} \label{prop:Atiyah}
Definition~\ref{defn:powerop} recovers Atiyah's total power operation in (untwisted) Real equivariant K-theory.
\end{prop}
\bp
Atiyah's $k^{\rm th}$ total power operation is induced by sending a vector bundle $E\to X$ to the external tensor power $E^{\boxtimes k}$ as a $\Sigma_k$-equivariant vector bundle over $X^{\times k}$; this original construction \cite[Proposition 2.2]{AtiyahPow} applies to complex vector bundles over a space, i.e., $\KU$-theory in degree zero. The identical formula applies to Real equivariant vector bundles, defining maps
\beq\nonumber
&&\KR(X\sq G)=\KR_G(X)\xrightarrow{\rm Atiyah} \KR_{G\wr \Sigma_k}(X^{\times k})=\KR(X^{\times k}\sq G\wr \Sigma_k)
\eeq
that equal $\Pow_k$ from Definition~\ref{defn:powerop}. For example, choose a $G$-invariant metric on~$E$, identify $E$ with a (finite rank) Hilbert bundle $E\to X$, and take the family of zero operators as a family of odd self-adjoint Fredholm operators. Alternatively, choose an embedding $E\subset X\times \H$ that realizes $E$ as a summand in the trivial Hilbert bundle, and take the family of Fredholm operators given by the projection with kernel $E$, see Remark~\ref{rmk:vectorbundle}. 
\ep

Next we verify that the twisted power operations constructed above satisfy properties analogous to the classical power operations on an $H_\infty$-ring spectrum as defined in~\cite[Chapter VIII,\S1]{bmms}. In our setting, these properties are inherited from natural transformations among the symmetric power functors on topological stacks, see Lemmas~\ref{lem:sym1} and~\ref{lem:sym2}. 

\begin{prop}\label{prop:powersproperties}
 For a twist $\tau\colon \X\to \pt\sq \PU^\pm_{\cCl_{p,q}}(\cH)$ and $k,j\in \mathbb{N}$, there is a natural isomorphism of twists 
\beq\label{eq:alpha}
 &&\alpha_{k,j}^*\Pow_{k+j}(\tau)\simeq \Pow_j(\tau)\boxtimes\Pow_k(\tau)\colon \X^{\times j}\sq \Sigma_j\times\X^{\times k}\sq \Sigma_k \to \pt \sq \PU^\pm_{\cCl_{p,q}^{\otimes (j+k)}}(\cH^{\otimes (j+k)})
\eeq
 for the natural transformation $\alpha_{k,j}$ from~\eqref{eq:defalpha}, and a natural isomorphism of twists
\beq\label{eq:beta}
 &&\beta^*\Pow_{jk}\simeq \Pow_j\circ \Pow_k\colon (\X^{\times k}\sq \Sigma_k)^{\times j}\sq \Sigma_j\to \pt \sq \PU^\pm_{\cCl_{p,q}^{\otimes (kj)}}((\cH^{\otimes k})^{\otimes j}))
\eeq
 for the natural transformation $\beta_{k,j}$ from~\eqref{eq:defbeta}.
Relative to these isomorphisms of twists, for $[F]\in \KR^{\tau}(\X)$ we have the equality of classes,
\begin{enumerate}
    \item[(i)] $\alpha^*\Pow_{k+j}([F])=\Pow_j([F])\boxtimes\Pow_k([F]) \in \KR^{\Pow_j(\tau)\boxtimes\Pow_k(\tau)}(\X^j\sq\Sigma_j\times \X^k\sq\Sigma_k)$, and
    \item[(ii)] $\beta^*\Pow_{jk}([F])=\Pow_j(\Pow_k([F])) \in \KR^{\Pow_j(\Pow_k(\tau))}(\X^{jk}\sq(\Sigma_k\wr\Sigma_j)$,
\end{enumerate}
in the respective cases. 
\end{prop}
\bp Since $ \alpha_{k,j}^*\Pow_{k+j}(\tau)$ is defined as the composition, 
\beq\nonumber
&&\X^{\times j}\sq \Sigma_j\times\X^{\times k}\sq \Sigma_k\xrightarrow{\alpha_{k,j}}\X^{\times {k+j}}\sq \Sigma_{k+j}\xrightarrow{\Sym^{k+j}\tau} \pt\sq \PU^\pm(\H)\wr \Sigma_{k+j}\to \pt \sq \PU^\pm(\H^{\otimes {k+j}})
\eeq
and $\Pow_j(\tau)\boxtimes\Pow_k(\tau)$ is the composition 
\beq\nonumber
\X^{\times j}\sq \Sigma_j\times\X^{\times k}\sq \Sigma_k  &\xrightarrow{\Sym^k\times \Sym^j\tau}& \pt\sq \PU^\pm(\H)\wr \Sigma_{j}\times \pt\sq \PU^\pm(\H)\wr \Sigma_{k}\nonumber\\
&\to&  \pt\sq \PU^\pm(\H^{\otimes j}) \times \pt\sq \PU^\pm(\H^{\otimes k})\nonumber\\
& \xrightarrow{\boxtimes}& \pt \sq \PU^\pm(\H^{\otimes {k+j}})\nonumber
\eeq
the assertion~\eqref{eq:alpha} follows from associativity and symmetry of the graded tensor product. 
The proof for \eqref{eq:beta} is similar. The equalities of classes (1) and (2) repeat the arguments for maps into the stacks $\Fred_{\cCl_{p,q}}(\cH)\sq \PU^\pm_{\cCl_{p,q}}(\cH)$. 
\ep

\subsection{Power operations applied to universal Thom classes}

\begin{proof}[Proof of Theorem~\ref{thm:A}]
First we verify that the Thom twists are compatible. The $k^{\rm th}$ total power of the $(p,q)^{\rm th}$ Thom twist is determined by the Real graded  homomorphisms
\beq\nonumber
\O_{p,q}\wr\Sigma_k\to \PU^\pm(\bS_{p,q})\wr\Sigma_k\to  \PU^\pm(\bS_{p,q}^{\otimes k})\simeq \PU^\pm(\bS_{kp,kq})
\eeq
where the first two arrows come from applying the power operation \eqref{eq:kpowertwist2} to the Thom twist Definition~\ref{defn:Thomtwist}, and the last isomorphism is a feature specific to the spinor representation, coming from the fact that $(\cCl_{p,q})^{\otimes k}=\cCl_{kp,qk}$ are isomorphisms of projective representations. This shows that the $k^{\rm th}$ power of the $(p,q)^{\rm th}$ Thom twist restricts from the $(kp,kq)^{\rm th}$ Thom twist.

Applying~\eqref{eq:cocyclepowerop} to the Thom cocycle~\eqref{eq:preThom} we obtain a map
$$
(\R^{p,q})^{\times k}\sq\O_{p,q}\wr \Sigma_k\to \End_{\Cl_{p,q}}(\bS_{p,q})^{\times k} \sq \PU^\pm(\bS_{p,q})\wr\Sigma_k\to \End_{\Cl_{p,q}}(\bS_{p,q}^{\otimes k}) \sq \PU^\pm(\bS_{p,q}^{\otimes k})
$$
for the family of Real operators (using~\eqref{eq:Thomfamofops}) given by
$$
(\R^{p,q})^{\times k}\subset \cCl_{p,q}^{\otimes k}\to \End_{\Cl_{p,q}^{\otimes k}}(\bS_{p,q}^{\otimes k})\simeq \End_{\Cl_{kq,kq}}(\bS_{kp,kq}). 
$$
This family of operators determines the $k^{\rm th}$ total power of the $(p,q)^{\rm th}$ Thom class, and restricts from the family of operators associated to the $(kp,kq)^{\rm th}$ Thom class. The theorem is proved. 
\ep

\begin{proof}[Proof of Corollary~\ref{cor:A}]
The same argument as in the proof of Theorem~\ref{thm:A} can be applied to the complex Clifford algebras (without Real structure) to obtain the KU statement, and the real Clifford algebras to obtain the KO statement. 
\ep

\begin{prop}\label{prop:ThomisocommuteswithPow}
For any twisting $\sigma\in \Twist(\X)$, the $k^{\rm th}$ total power operation commutes with the Thom isomorphism~\eqref{eq:nonuniThomiso}, 
\beq\nonumber
\begin{tikzpicture}[baseline=(basepoint)];
\node (A) at (0,0) {$\KR^\sigma(\X)$};
\node (B) at (5,0) {$\KR^{\sigma\otimes\SW_\V}(\V)$};
\node (C) at (0,-1.25) {$\KR^{\Pow_k\sigma}(\X^{\times k}\sq \Sigma_k)$};
\node (D) at (5,-1.25) {$\KR^{\Pow_k\SW(\sigma\otimes \V)}(\V^{\times k}\sq \Sigma_k)$};
\draw[->] (A) to node [above] {$\smallsmile \Th_\V$}  (B);
\draw[->] (B) to node [right] {$\Pow_k$} (D);
\draw[->] (A) to node [left] {$\Pow_k$} (C);
\draw[->] (C) to node [below] {$\smallsmile \Th_{\V^{\times k}\sq \Sigma_k}$} (D);
\path (0,-.75) coordinate (basepoint);
\end{tikzpicture}
\eeq
with analogous commutative squares for the $\KU$- and $\KO$-cases. 
\end{prop}

\bp
It suffices to consider the universal case, 
\beq\nonumber
\begin{tikzpicture}[baseline=(basepoint)];
\node (A) at (0,0) {$\KR^\sigma(\pt\sq O_{p,q})$};
\node (B) at (7,0) {$\KR^{\sigma\otimes \SW_{p,q}}(\R^{p,q}\sq \O_{p,q})$};
\node (C) at (0,-1.25) {$\KR^{\Pow_k\sigma}(\pt \sq \O_{p,q} \wr \Sigma_k)$};
\node (D) at (7,-1.25) {$\KR^{\Pow_k(\sigma\otimes \SW_{p,q})}((\R^{p,q})^{\times k}\sq \Sigma_k).$};
\draw[->] (A) to node [above] {$\Th_{p,q}$} (B);
\draw[->] (B) to node [right] {$\Pow_k$} (D);
\draw[->] (A) to node [left] {$\Pow_k$} (C);
\draw[->] (C) to node [above] {$\Th_{\R^{pk,qk}\sq \O_{p,q}\wr \Sigma_k}$} (D);
\path (0,-.75) coordinate (basepoint);
\end{tikzpicture}
\eeq
Lemma~\ref{lem:powmulttwist} provides the isomorphism of twistings $\Pow_k(\sigma\otimes \tau_{p,q})\simeq \Pow_k(\sigma)\otimes \Pow_k(\tau_{p,q})$. Then
after identifying the Thom class $\Th_{\R^{pk,qk}\sq \O_{p,q}\wr \Sigma_k}$ with the restriction of the universal Thom class $\Th_{pk,qk}$, the commutativity of the above square follows from multiplicativity of the power operation and Theorem~\ref{thm:A}: the maps from the upper left to the lower right are both given by the external product with $\Pow_k\Th_{p,q}=\res(\Th_{pk,qk})$. \ep

\subsection{The twisted Atiyah--Bott--Shapiro map commutes with the power operation}

Let $X$ and $Y$ be smooth manifolds with the smooth action of a compact Lie group~$G$.

\begin{thm}\label{thm:ABS}
For all $n\in \N$, the twisted Atiyah--Bott--Shapiro map~\eqref{Eq:twistedABS} commutes with the~$n^{\rm th}$ total power operation, i.e., for any Real $G$-equivariant proper submersion $\pi_!\colon X\to Y$, the diagram commutes
\beq\label{eq:bigABS}
\begin{tikzpicture}[baseline=(basepoint)];
\node (A) at (0,0) {$\KR^{\sigma \otimes \tau_\nu^{-1}}(X\sq G)$};
\node (B) at (6,0) {$\KR^{\sigma\otimes \tau_\rho^{-1}}(Y\sq G)$};
\node (C) at (0,-1.25) {$\KR^{\Pow_n(\sigma \otimes \tau_\nu^{-1})}(X^{\times n}\sq G\wr \Sigma_n)$};
\node (D) at (6,-1.25) {$\KR^{\Pow_n(\sigma\otimes\tau_\rho^{-1})}(Y^{\times n}\sq G\wr \Sigma_n).$};
\draw[->] (A) to node [above] {$\pi_!^{\rm top}$} (B);
\draw[->] (B) to node [right] {$\Pow_n$} (D);
\draw[->] (A) to node [left] {$\Pow_n$} (C);
\draw[->] (C) to node [below] {$(\Pow_n\pi)_!^{\rm top}$} (D);
\path (0,-.75) coordinate (basepoint);
\end{tikzpicture}
\eeq
The analogous statement also holds in KU- and KO-theory. 
\end{thm}

\bp
Extend~\eqref{eq:bigABS} to the diagram
\beq\nonumber
\resizebox{1\textwidth}{!}{$\begin{tikzpicture}[baseline=(basepoint)];
\node (A) at (0,0) {$\KR^{\sigma\otimes\tau_\nu^{-1}}(X\sq G)$};
\node (B) at (4,0) {$\KR^{\sigma}(\nu)$};
\node (C) at (8.5,0) {$\KR^{\sigma}((\R^{p,q}\times Y)\sq G)$};
\node (D) at (13.5,0) {$\KR^{\sigma\otimes\tau_\rho^{-1}}(Y\sq G)$};
\node (AA) at (0,-1.5) {$\KR^{\Pow_n(\sigma\otimes\tau_\nu^{-1})}(X^{\times n}\sq G\wr\Sigma_n)$};
\node (BB) at (4,-1.5) {$\KR^{\Pow_n\sigma}(\nu^{\times n}\sq \Sigma_n)$};
\node (CC) at (8.5,-1.5) {$\KR^{\Pow_n\sigma}((\R^{np,nq}\times Y^{\times n})\sq G\wr \Sigma_n)$};
\node (DD) at (13.5,-1.5) {$\KR^{\Pow_n(\sigma\otimes\tau_\rho^{-1})}(Y^{\times n}\sq G\wr \Sigma_n )$};
\draw[->] (A) to node [above] {$\smallsmile\Th_{\nu}$}  (B);
\draw[->] (B) to node [above] {extend by 0} (C);
\draw[->] (D) to node [above] {$\smallsmile\Th_\rho$}(C);
\draw[->] (AA) to (BB);
\draw[->] (BB) to (CC);
\draw[->] (DD) to (CC);
\draw[->] (A) to node [left] {$\Pow_n$} (AA);
\draw[->] (B) to node [left] {$\Pow_n$} (BB);
\draw[->] (C) to node [left] {$\Pow_n$} (CC);
\draw[->] (D) to node [left] {$\Pow_n$} (DD);
\path (0,-.75) coordinate (basepoint);
\end{tikzpicture}$}
\eeq
where the top row determines the twisted Atiyah--Bott--Shapiro map for fibration $X\sq G\to Y\sq G$, the lower row determines the twisted Atiyah--Bott--Shapro map for the fibration $X^{\times n}\sq G\wr\Sigma_n\to Y^{\times n}\sq G\wr\Sigma_n$, and the vertical arrows are the $n^{\rm th}$ total power operation. The canonical isomorphisms of twists from 
$$
\Pow_n(\sigma\otimes \tau_\nu^{-1})\simeq \Pow_n(\sigma)\otimes \tau_{\nu^{\oplus n}}^{-1},\qquad \Pow_n(\sigma\otimes \tau_\rho^{-1})\simeq \Pow_n(\sigma)\otimes \tau_{\rho^{\oplus n}}^{-1}
$$
follow from Lemma~\ref{lem:powmulttwist} together with the fact that the (external) direct sum of bundles is compatible with the (external) tensor product of Thom classes. Commutativity of the squares on the far left and the far right in the diagram above then follows from  Proposition~\ref{prop:ThomisocommuteswithPow}. The middle square commutes by inspection: applying the power and extending by zero yields the same cocycle as first extending by zero and then applying the power operation. The $\KU$- and $\KO$-variants follow from the same structures. The theorem is proved. 
\ep

\begin{proof}[Proof of Theorem~\ref{thm:3}]
Using the pushforward to $\pt\sq G$ from Example~\ref{ex:ABSBG}, the result follows from setting $Y=\pt$ in Theorem~\ref{thm:ABS}. 
\ep

\appendix

\section{Topological groupoids, stacks, and gerbes}\label{appen:stacks}

Following~\cite{FHTI}, all topological spaces are assumed to be locally contractible, paracompact, and completely regular. This implies the existence of partitions of unity~\cite{Doldpartition} and locally contractible slices for actions by compact Lie groups~\cite{Mostow,Palais}. Let $\Top$ denote the category of such topological spaces and continuous maps. 

 A \emph{stack} is a presheaf of groupoids on $\Top$ satisfying descent for all \'etale covers. A \emph{topological stack} is a stack $\X$ for which there exists a representable epimorphism $X\to \X$ where $X$ is a topological space viewed as a stack, e.g., see~\cite[Definition 2.3]{Heinloth}. A choice of $X\to \X$ is an \emph{atlas}, which we emphasize is not part of the data of a topological stack. 
 
\subsection{Topological groupoids and local quotient stacks}\label{sec:locquostacks}


\begin{defn}
A \emph{topological groupoid} is a groupoid object in topological spaces.
\end{defn}

We use notation $\bX=\{\bX_1\rightrightarrows \bX_0\}$ to denote the space of objects $\bX_0$ and morphisms $\bX_1$ of a topological groupoid $\bX$. Source, target, unit, and compositions maps are denoted $s,t,u,$ and~$c$. Topological groupoids, continuous functors and continuous natural transformations form a (strict) 2-category. There is a functor from the 2-category of topological groupoids to the bicategory of topological stacks, 
\beq\label{eq:grpdtostack}
&&\Grpd\to \Stack,
\eeq
sending a groupoid $\bX$ to the stack of locally trivial $\bX$-bundles on $\Top$, e.g., see~\cite[\S3]{Heinloth}.


\begin{ex}
Given a topological group $G$ acting continuously on a space $X$, the \emph{quotient groupoid} has objects the space $X$ and morphisms the space $G\times X$. The source map is the projection, the target map is the action map, and composition is determined by multiplication in $G$. The stack $X\sq G$ assigns to a space $Y$ the groupoid whose objects are principal $G$-bundles $P\to Y$ together with a $G$-equivariant map $P\to X$. 
\end{ex}

Conversely, given a topological stack with a choice of \emph{atlas} $\bX_0\to \X$ (i.e., representable epimorphism), there is a topological groupoid with objects $\bX_0$ and morphisms $\bX_1=\bX_0\times_\X \bX_0$. In this case, the stack of principal $\bX$ bundles is equivalent to the stack $\X$, and 
we say \emph{the groupoid $\bX$ presents the stack $\X$}, or that \emph{$\bX$ is a presentation of $\X$}. 

Given a continuous functor $F\colon \bX\to \bY$ between topological groupoids, let $F_0\colon \bX_0\to \bY_0$ and $F_1\colon \bX_1\to \bY_1$ denote the maps on objects and morphisms, respectively. Consider the diagrams of topological spaces
\beq\label{eq:essentialequiv}
\begin{tikzpicture}[baseline=(basepoint)];
\node (A) at (0,0) {$\bY_1\times_{\bY_0}\bX_0$};
\node (B) at (3,0) {$\bX_0$};
\node (C) at (0,-1.5) {$\bY_1$};
\node (D) at (3,-1.5) {$\bY_0$};
\draw[->] (A) to node [above] {$\pr_2$} (B);
\draw[->] (B) to node [right] {$\pr_1$} (D);
\draw[->] (A) to node [left] {$s$} (C);
\draw[->] (C) to node [below] {$F_0$} (D);
\path (0,-.75) coordinate (basepoint);
\end{tikzpicture}\qquad 
\begin{tikzpicture}[baseline=(basepoint)];
\node (A) at (0,0) {$\bX_1$};
\node (B) at (3,0) {$\bY_1$};
\node (C) at (0,-1.5) {$\bX_0\times \bX_0$};
\node (D) at (3,-1.5) {$\bY_0\times \bY_0.$};
\draw[->] (A) to node [above] {$F_1$} (B);
\draw[->] (B) to node [right] {$s\times t$} (D);
\draw[->] (A) to node [left] {$s\times t$} (C);
\draw[->] (C) to node [below] {$F_0\times F_0$} (D);
\path (0,-.75) coordinate (basepoint);
\end{tikzpicture}
\eeq

\begin{defn}
A continuous functor $F\colon \bX\to \bY$ is a \emph{local equivalence} if $\bY_1\times_{\bY_0}\bX_0\to \bY_0$ on the left in~\eqref{eq:essentialequiv} admits local sections, and the diagram on the right is a pullback. 
\end{defn}

\begin{rmk} \label{defn:localequiv}
The two conditions in Definition~\ref{defn:localequiv} are continuous versions of essential surjectivity and fully faithfulness. 
\end{rmk}

Under~\eqref{eq:grpdtostack}, local equivalences of topological groupoids are sent to 1-isomorphisms of stacks. Pronk shows that \emph{all} 1-isomorphisms of stacks come from local equivalences. 

\begin{thm}[{\cite{Pronk}}]\label{thm:Pronk}
After localizing the 2-category of topological groupoids at the local equivalences, \eqref{eq:grpdtostack} induces an equivalence with the bicategory of topological stacks.
\end{thm}

For stacks $\X$ and $\Y$ presented by groupoids $\bX$ and $\bY$ respectively, a consequence of the above is that a map of stacks $\X\to \Y$ can be presented by a span
\beq\label{eq:anafunctor}
\bX\xleftarrow{\sim} \mathbb{Z}\rightarrow \bY
\eeq
where the left arrow in~\eqref{eq:anafunctor} is a local equivalence. 


\begin{ex}
Given a groupoid $\bX$ and a continuous map $\phi\colon S\to \bX_0$, define the \emph{pullback groupoid} $\phi^*\bX$ with objects $(\phi^*\bX)_0=S$ and morphisms the pullback
\beq\nonumber
\begin{tikzpicture}[baseline=(basepoint)];
\node (A) at (0,0) {$(\phi^*\bX)_1$};
\node (B) at (3,0) {$\bX_1$};
\node (C) at (0,-1.5) {$S\times S$};
\node (D) at (3,-1.5) {$\bX_0\times \bX_0$};
\draw[->] (A) to  (B);
\draw[->] (B) to node [right] {$s\times t $} (D);
\draw[->] (A) to (C);
\draw[->] (C) to node [below] {$\phi\times \phi$} (D);
\path (0,-.75) coordinate (basepoint);
\end{tikzpicture}\qquad 
\eeq
i.e., for any pair $x,y\in S$ we take all the morphisms between $\phi(x)$ and $\phi(y)$ in $\bX_1$. There is an evident functor $\phi^*\bX\to \bX$. If $\phi\colon S=\coprod U_i\to \bX_0$ is an open cover of $\bX_0$, then $\phi^*\bX\to \bX$ is a local equivalence. In particular, an open cover $\{U_i\}$ of a topological space $X$ determines the \emph{\v{C}ech groupoid} $\check{C}(U_i)$, which is equivalent to~$X$.
\end{ex}

\begin{defn}\label{sec:coarse}
The \emph{coarse moduli space} of a groupoid $\bX$, denoted $|\bX|$, is the space of isomorphism classes of objects endowed with the quotient topology for $q\colon \bX_0\to |\bX|$.
\end{defn}

 If $\bX\to \bY$ is a local equivalence, the induced map of coarse moduli spaces $|\bX|\to |\bY|$ is a homeomorphism~\cite{Noohiclassifying}. For a subspace $S\subset |\bX|$ let $\bX_S:=\phi^*\bX$ denote the pullback groupoid for the map $\phi\colon q^{-1}(S)\to \bX_0$. When $S$ is open (respectively, closed) we refer to~$\bX_S$ as an \emph{open} (respectively, \emph{closed}) subgroupoid of $\bX$. An \emph{open covering} of a topological groupoid is an open cover $\{U_i\}$ of $|\bX|$, determining the collection $\{\bX_{U_i}\}$ of open subgroupoids. A \emph{local quotient groupoid} \cite[\S{A.2.2}]{FHTI} is a topological groupoid that is locally of the form $X\sq G$ for $X$ a topological space and $G$ a compact Lie group, i.e., a groupoid $\bX$ with an open cover $\{U_i\}_{i\in I}$ of $|\bX|$ for which there exist local equivalences~$\bX_{U_i}\simeq X_i\sq G_i$ for each $i\in I$. Define the following subcategory of topological stacks. 

\begin{defn}[{\cite[\S6]{Heinloth}}]\label{defn:locqutstack}
The bicategory of \emph{local quotient stacks} is the full sub-bicategory of topological stacks in the image of local quotient groupoids under~\eqref{eq:grpdtostack}.
\end{defn}

\begin{rmk}\label{rmk:linearization}
For a local quotient stack, the diagonal $\X\to \X\times \X$ is a proper map; this implies that the automorphism group of any object is compact. For smooth stacks, the converse holds: any proper Lie groupoid represents a local quotient stack~\cite[Theorem~2.3]{Zung}. A smooth stack determines a topological one, e.g., by taking the topological groupoid underlying a Lie groupoid presentation. Hence, all proper Lie groupoids determine objects in the bicategory of Definition~\ref{defn:locqutstack}. 
\end{rmk}

\subsection{$G$-stacks and Real stacks}\label{sec:Realstacks}

Let $G$ be a Lie group. A \emph{$G$-stack} is a topological stack~$\X$ with action map~$\mu\colon G\times \X\to \X$ and 2-morphisms in the diagrams
\beq
\begin{tikzpicture}[baseline=(basepoint)];
\node (A) at (0,0) {$G\times G\times \X $};
\node (B) at (3,0) {$G\times \X$};
\node (C) at (0,-1.5) {$G\times \X$};
\node (D) at (3,-1.5) {$\X$};
\node (E) at (1.5,-.75) {$\twocommute$};
\draw[->] (A) to node [above] {$\id_G\times \mu$} (B);
\draw[->] (B) to node [right] {$m\times \id_\X$} (D);
\draw[->] (A) to node [left] {$\mu$} (C);
\draw[->] (C) to node [below] {$\mu$} (D);
\path (0,-.75) coordinate (basepoint);
\end{tikzpicture}\qquad 
\begin{tikzpicture}[baseline=(basepoint)];
\node (A) at (0,0) {$G\times \X $};
\node (B) at (3,0) {$\X$};
\node (C) at (1.5,-1.5) {$\X$};
\node (E) at (1.5,-.75) {$\twocommute$};
\draw[->] (A) to node [above] {$\mu$} (B);
\draw[->] (C) to node [left] {$1\times \id_\X$} (A);
\draw[->] (C) to node [right=10pt] {$\id_\X$} (B);
\path (0,-.75) coordinate (basepoint);
\end{tikzpicture}\label{Eq:Gstack}
\eeq
satisfying further coherence conditions. A \emph{$G$-equivariant map} of $G$-stacks is a map of topological stacks $f\colon \X\to \Y$ together with 2-commuting data 
\beq
\begin{tikzpicture}[baseline=(basepoint)];
\node (A) at (0,0) {$G\times \X $};
\node (B) at (3,0) {$G\times \Y$};
\node (C) at (0,-1.5) {$\X$};
\node (D) at (3,-1.5) {$\Y$};
\node (E) at (1.5,-.75) {$\twocommute$};
\draw[->] (A) to node [above] {$\id_G\times f$} (B);
\draw[->] (B) to node [right] {$\mu_Y$} (D);
\draw[->] (A) to node [left] {$\mu_X$} (C);
\draw[->] (C) to node [below] {$f$} (D);
\path (0,-.75) coordinate (basepoint);
\end{tikzpicture}\label{eq:equivariantstack}
\eeq
satisfying compatibility conditions for the 2-morphisms in~\eqref{Eq:Gstack}. An \emph{isomorphism} between $G$-equivariant maps is an isomorphism of maps of stacks $f\Rightarrow f'$ compatible with the equivariance data~\eqref{eq:equivariantstack}. All together, $G$-stacks form a bicategory with an evident forgetful map to topological stacks, e.g., see \cite[\S3]{SP11}. 

The following is well-known.

\begin{lem}
There is an equivalence of bicategories between stacks with $G$-action and stacks over $\pt\sq G$. 
\end{lem}
\begin{proof}[Proof sketch]
First we sketch the construction of a functor from $G$-stacks to stacks over $\pt\sq G$. Given a $G$-stack $\X$, the quotient stack $\X\sq G$ 
has an evident map to $\pt\sq G$, and a map of $G$-stacks $\X\to \X'$ determines a map on quotients $\X\sq G\to \X'\sq G$ over $\pt\sq G$. 

Given a stack $\Y$ over $\pt\sq G$, we obtain a $G$-stack $\X$ by pullback 
\beq\label{eq:Gpullback}
&&\begin{tikzpicture}[baseline=(basepoint)];
\node (A) at (0,0) {$\X$};
\node (B) at (3,0) {$*$};
\node (C) at (0,-1.5) {$\Y$};
\node (D) at (3,-1.5) {$\pt\sq G$};
\node (P) at (.7,-.4) {\scalebox{1.5}{$\lrcorner$}};
\draw[->] (B) to  (D);
\draw[->] (C) to (D);
\draw[->] (A) to (B);
\draw[->] (A) to (C);
\path (0,-.75) coordinate (basepoint);
\end{tikzpicture}
\eeq
i.e., $\X$ is the total space of the principal $G$-bundle classified by $\Y\to\pt\sq G$. A map of stacks $\Y\to \Y'$ over $\pt\sq G$ determines a map of principal $G$-bundles, and hence a morphism of $G$-stacks. This completes the sketch of a functor from stacks over $\pt\sq G$ to $G$-stacks. 

To see that the above functors are inverse to each other, observe that the description of $\X$ as a principal $G$-bundle over $\Y$~\eqref{eq:Gpullback} produces an equivalence $\Y\simeq \X\sq G$. Furthermore, this equivalence is natural in the sense that a map $\Y\to\Y'$ is compatible with the induced morphism $\X\sq G\to \X'\sq G$. 
Finally, the equivalence on 2-morphisms and compatibility with compositions comes from universal properties of $G$-quotients and pullbacks~\eqref{eq:Gpullback}. 
\ep


Using~\eqref{eq:grpdtostack}, a groupoid object in $G$-spaces determines a stack with $G$-action. There are two valid definitions of a $G$-action on a groupoid: one for weak actions and one for strict actions, in the sense that diagrams analogous to~\eqref{Eq:Gstack} commute strictly. By \cite[Proposition 1.5]{Romagny} these two notions are equivalent: any groupoid with weak $G$-action is equivalent to a groupoid with strict $G$-action. Hence, in certain constructions it can be useful to choose a presentation of a given $G$-stack in which the action is strict. However, strict $G$-actions need not pass to strict $G$-actions under an equivalence of topological groupoids. 


\begin{defn}\label{defn:Realstack}
The bicategory of \emph{Real stacks} has objects stacks with $C_2$-action, 1-morphisms $C_2$-equivariant map of stacks, and 2-morphisms isomorphisms between $C_2$-equivariant maps. 
\end{defn}

\begin{ex}
We recall from~\cite{Atiyahreal} that a \emph{Real space} is a space with $C_2$-action and a \emph{Real map} of Real spaces is a $C_2$-equivariant map. There is an evident faithful functor from Atiyah's category of Real spaces to the bicategory of Real stacks. 
\end{ex}

\begin{ex}
We recall from \cite{Atiyahbottp} that a \emph{Real group} is a group $G$ with involution $(\overline{\phantom{g}})\colon G\to G$, see~\S\ref{sec:Realgroups} above. This involution determines a $C_2$-action on the groupoid $\pt\sq G$, and hence a Real structure on the underlying stack $\pt\sq G$. 
\end{ex}

\begin{defn}
    A \emph{Real presentation} of a Real stack is a groupoid presentation with a strict $C_2$-action (which exists by \cite[Proposition 1.5]{Romagny}). A \emph{Real atlas} of a Real stack $\X$ is a $C_2$-equivariant atlas $\bX_0\to \X$ whose associated groupoid presentation is Real. 
\end{defn}

\subsection{Symmetric powers of $G$-stacks}\label{sec:sympowers}
\begin{defn}
The \emph{$k^{\rm th}$ symmetric power} of a topological stack $\X$ is 
$$
\Sym^k(\X):=\X^{\times k}\sq \Sigma_k.
$$ 
This extends to a functor $\Sym^k\colon \Stack\to \Stack/(\pt\sq \Sigma_k)$ from stacks to over~$\pt\sq \Sigma_k$. 
\end{defn}

\begin{lem} \label{lem:sym1}
For $j,k\in \N$, the map of stacks
$$
\X^{\times j}\sq \Sigma_j\times \X^{\times k}\sq \Sigma_k\simeq (\X^{\times (j+k)})\sq (\Sigma_j\times \Sigma_k)\to (\X^{\times (j+k)})\sq \Sigma_{j+k}
$$
induced by the inclusion of groups $\Sigma_j\times \Sigma_k\hookrightarrow \Sigma_{j+k}$ determines a natural transformation 
\beq\label{eq:defalpha}
\alpha_{j,k}\colon \Sym^j(-)\times \Sym^k(-)\to \Sym^{j+k}(-)
\eeq
between functors from stacks to stacks over $\pt\sq \Sigma_{j+k}$. 
\end{lem}

\begin{lem}\label{lem:sym2}
For $j,k\in \N$, the map of stacks
$$
(\X^{\times k}\sq \Sigma_k)^{\times j}\sq \Sigma_j\simeq \X^{\times jk}\sq \Sigma_k\wr\Sigma_j\to \X^{\times jk}\sq \Sigma_{jk}
$$
induced by the inclusion of groups $\Sigma_k\wr\Sigma_j=(\Sigma_k^{\times j})\rtimes \Sigma_j\subseteq\Sigma_{jk}$ determines a natural transformation 
\beq\label{eq:defbeta}
\beta_{j,k}\colon \Sym^j(\Sym^k(-))\to \Sym^{jk}(-)
\eeq
between functors from stacks to stacks over $\pt\sq \Sigma_{jk}$. 
\end{lem}

We omit the (straightforward) proofs of the above. 

\begin{ex}
For a global quotient stack $X\sq G$, there is an equivalence of topological stacks $\Sym^k(X\sq G)\simeq X^{\times k}\sq G\wr \Sigma_k$ for the action of the wreath product $G\wr \Sigma_k$ on $X^{\times k}$. We recall the homomorphism of groups
$$
G\times \Sigma_k\xhookrightarrow{\Delta\times \id_{\Sigma_k}} G^{\times k}\rtimes \Sigma_k\simeq G\wr\Sigma_k
$$
induced by the diagonal $G\hookrightarrow G^{\times k}$. This gives a map 
$$
(X^{\times k} \sq \Sigma_k)\sq G \simeq X^{\times k}\sq (\Sigma_k\times G)\to \Sym^k(X\sq G)
$$

\end{ex}

\begin{ex} More generally, for a $G$-stack $\X$, there is a map
$$
(\X^{\times k}\sq \Sigma_k)\sq G\simeq \X^{\times k}\sq (\Sigma_k\times G)\to \X^{\times k}\sq G\wr \Sigma_k\simeq \Sym^k(\X\sq G). 
$$
The quotient on the left gives $\Sym^k(\X)$ the structure of a $G$-stack, which maps to the symmetric power $\Sym^k(\X\sq G)$ of the $G$-quotient. 
\end{ex}

\begin{defn}
The \emph{$k^{\rm th}$ symmetric power of a $G$-stack $\X$} is the stack $\Sym(\X^{\times k})\simeq \X^{\times k}\sq \Sigma_k$ with the diagonal $G$-action on $\X^{\times k}$. For $G=C_2$, this defines the symmetric powers of a Real stack as a Real stack. 
\end{defn}

\subsection{Real graded  gerbes}\label{sec:appenRealgerbe}

There are several equivalent approaches to gerbes on topological stacks, all tracing back to~\cite{Giraud}. Perhaps the most practical point of view treats gerbes as central extensions of topological groupoids, e.g., see~\cite[\S2.2]{FHTI}. On the other hand, the global geometry of gerbes and the larger bicategory in which they reside is most transparent in the language of higher categorical principal bundles, e.g., see~\cite[\S{A}]{Gajer} and~\cite{Moerdijk2002}. The literature on this subject is vast, and we do not give exhaustive references below. 
However, for the reader's convenience we indicate some of the proofs. 

We begin with the basic notions. Multiplication on $\U(1)$ makes $\pt\sq \U(1)$ an abelian group object in stacks, yielding a bicategory of stacks with $\pt\sq \U(1)$-action~\cite[\S3.3]{SP11}. 

\begin{defn}\label{defn:gerbe}
Let $\X$ be a topological stack. A (bundle) \emph{gerbe} over $\X$ is a $\pt\sq \U(1)$-principal bundle $\widehat{\X} \to \X$, i.e., $\widehat{\X}$ is equipped with a $\pt\sq \U(1)$-action and there exists an atlas $u\colon X\to \X$ with a $\pt\sq \U(1)$-equivariant equivalence $u^*\widehat{\X}\simeq X\times \pt\sq \U(1)$ over $X$. 
\end{defn}

For a local quotient stack $\X$ and a gerbe $\widehat{\X}\to \X$, the stack $\widehat{\X}$ is again a local quotient stack~\cite[Corollary 6.3]{Heinloth}. Gerbes are the objects of a symmetric monoidal 2-stack~\cite[\S4.1]{NikolausSchweigert}; in particular for each stack $\X$ there is a symmetric monoidal 2-groupoid $\Gerbe_\X$ and a map of stacks $\X\to \Y$ determines a (weak) pullback 2-functor $\Gerbe_\Y\to \Gerbe_\X$. 


Next we observe that $\pt\sq \U(1)$ is an abelian group object in Real stacks for the Real structure on $\pt\sq \U(1)$ generated by complex conjugation on $\U(1)$. 

\begin{defn} 
Let $\X$ be a Real stack. A \emph{Real gerbe} over $\X$ is a Real stack $\widehat{\X}$ and a Real map $\widehat{\X}\to \X$ that (after forgetting the Real structure) is a gerbe over $\X$ and with the property that on a Real atlas $X\to \X$, the $\pt\sq \U(1)$-equivariant isomorphism $u^*\widehat{\X}\simeq X\times \pt\sq \U(1)$ is furthermore a Real isomorphism for the Real structure on $\widehat{\X}$, the Real structure on $\X$, and the Real structure on $\pt\sq \U(1)$ described above. 
\end{defn}

\begin{defn}
A \emph{grading} for a topological stack $\X$ is a map of stacks  $\X\to \pt\sq (\Z/2)$. 
\end{defn}
We observe that a grading is the same data as a metrized real line bundle on $\X$. 

\begin{defn} \label{defn:gradedgerbe}
A \emph{graded gerbe} over a stack $\X$ is a grading on $\X$ and a $\U(1)$-bundle gerbe 
\beq\label{eq:gerbedefn}
&&\pt\sq \U(1) \to \widehat{\X}\to  \X,\qquad \X\to \pt\sq(\Z/2)
\eeq
A \emph{Real graded  gerbe} over a $C_2$-stack $\X$ is a graded gerbe~\eqref{eq:gerbedefn} in $C_2$-stacks: $\pt\sq \U(1)$ carries a $C_2$-action from complex conjugation, $\widehat{\X}$ is endowed with a $C_2$-action, the maps in~\eqref{eq:gerbedefn} are $C_2$-equviariant, and the grading $\X\to \pt\sq(\Z/2)$ is $C_2$-equivariant for the trivial $C_2$-action on~$\pt\sq(\Z/2)$. 
\end{defn}

Real graded gerbes are the objects of 2-stack with a forgetful functor to stacks. The 1-morphisms are Real $\pt\sq \U(1)$-equivariant maps $\widehat{\X}\to \widehat{\Y}$ covering a Real map $\X\to \Y$ of stacks and compatible with the gradings.

\begin{defn}\label{defn:trivializationgerbe}
The \emph{trivial} Real graded gerbe over a Real stack $\X$ is the Real stack $\pt\sq \U(1)\times \X$ with diagonal Real structure and grading given by the trivial map to~$*$ followed by the unit map $*\to \pt\sq(\Z/2)$. A \emph{trivialization} of a Real graded gerbe is an equivalence with the trivial Real graded gerbe. 
\end{defn}


\begin{rmk}\label{rmk:orientifold2}
Given a Real gerbe $\widehat{\X}\to \X$, we obtain a gerbe $\widehat{\X}\sq C_2\to \X\sq C_2$ over $\pt\sq C_2$. This is the \emph{orientifold} associated to a Real gerbe. Conversely, given a gerbe $\widehat{\Y}\to \Y$ over $\pt\sq C_2$, the pullbacks~\eqref{eq:Gpullback} construct a $C_2$-equivariant gerbe, compare~\cite{Vienna}. 
\end{rmk}

\begin{ex}\label{ex:centralext}
For a compact Lie group $G$ with involution (i.e., Real structure), the data of a Real graded  central extension determines a Real graded  gerbe over $\pt\sq G$,
\beq\label{eq:gradedcentralextofgroups}
\begin{array}{c} 1\to \U(1)  \to \widehat{G}\to  G\to 1, \\ G\xrightarrow{\epsilon} \Z/2 \end{array} \quad \implies \quad \begin{array}{c} \pt\sq \U(1)\to \widehat{\X}={[\pt\sq \hat{G}]} \to {[\pt\sq G]}=\X \\ 
{[\pt\sq G]}  \xrightarrow{\epsilon} [\pt\sq(\Z/2)] \end{array}
\eeq
where the $C_2$-action is generated by the involution on $G$ and complex conjugation on $\U(1)$. 
\end{ex}

Next we unpack Definition~\ref{defn:gradedgerbe} in a groupoid presentation of $\X$. 

\begin{defn}[{\cite[Definition 2.6]{FHTI}}]\label{defn:gradedcentralext}
A \emph{$\pt\sq \U(1)$-central extension} of a topological groupoid~$\bX$ is a principal $\U(1)$-bundle $P\to \bX_1$ and an isomorphism of $\U(1)$-bundles over $\bX_2=\bX_1\times_{\bX_0} \bX_1$ 
\beq\label{eq:lambda}
\lambda\colon p_2^*P\otimes_{\U(1)} p_1^*P\to c^*P,\qquad p_i(f_2,f_1)=f_i, \ c(f_2,f_1)=f_2\circ f_1
\eeq
satisfying the condition that the diagram over $\bX_3=\bX_1\times_{\bX_0} \bX_1\times_{\bX_0} \bX_1$ commute
\beq\nonumber
\begin{tikzpicture}[baseline=(basepoint)];
\node (A) at (0,0) {$p_3^*P\otimes_{\U(1)} p_2^*P\otimes_{\U(1)} p_1^*P$};
\node (B) at (5,0) {$c_{32}^*P\otimes p_3^*P$};
\node (C) at (0,-1.5) {$p_1^*P\otimes c_{21}^*P$};
\node (D) at (5,-1.5) {$c_{321}^*P$};
\draw[->] (A) to node [above] {$p_{32}^*\lambda\otimes \id$}  (B);
\draw[->] (B) to node [right] {$(c_{32}\times p_1)^*\lambda $} (D);
\draw[->] (A) to node [left] {$\id\otimes p_{21}^*\lambda$} (C);
\draw[->] (C) to node [below] {$(p_3\times c_{21})^*\lambda$} (D);
\path (0,-.75) coordinate (basepoint);
\end{tikzpicture}\qquad \begin{array}{l} p_i(f_3,f_2,f_1)=f_i\\ p_{ij}(f_3,f_2,f_1)=(f_i,f_j) \\ c_{ij}(f_3,f_2,f_1)=f_i\circ f_j\\ c_{321}(f_3,f_2,f_1)=f_3\circ f_2\circ f_1\end{array}
\eeq
for the indicated maps out of $\bX_3$. 
A \emph{graded central extension} is a central extension together with a functor $\bX\to \pt\sq(\Z/2)$, i.e., a map $\bX_1\to \Z/2$ satisfying a cocycle condition. 
\end{defn}

A graded central extension determines a groupoid $\widehat{\bX}=\{P\rightrightarrows \bX_0\}$ over $\bX$ and a continuous functor $\epsilon\colon \bX\to \pt\sq(\Z/2)$: the source and target in $\widehat{\bX}$ are $P\to\bX_1\rightrightarrows \bX_0$, and composition is determined by~\eqref{eq:lambda}. If $\bX$ is a local quotient groupoid, then so is $\widehat{\bX}$~\cite[Corollary 2.18]{FHTI}.

The grading $\epsilon$ allows us to view the $\U(1)$-bundle $P\to \bX_1$ as having even or odd fibers $P_x$, corresponding to the value of $\epsilon(x)\in \{0,1\}$. The cocycle condition for $\bX_1\to \Z/2$ implies that~\eqref{eq:lambda} is a map of graded $\U(1)$-bundles for the graded tensor product. 
This permits the following (graded) tensor product of central extensions. 

\begin{defn}[{\cite[page~15]{FHTI}}]\label{defn:gradedtensorproduct}
The \emph{external tensor product} of graded $\pt\sq \U(1)$-central extensions $(P,\lambda,\epsilon)$ of $\bX$ and $(Q,\mu,\delta)$ of $\bY$ is determined by the external product of principal $\U(1)$-bundles $P\boxtimes_{\U(1)} Q\to \bX_1\times \bY_1$, the sum of the gradings $(\epsilon+ \delta)\colon\bX_1\times \bY_1\to \Z/2$, and the map of principal bundles
\beq
p_2^*(P\boxtimes_{\U(1)} Q)\boxtimes_{\U(1)} p_1^*(P\boxtimes_{\U(1)} Q) &\xrightarrow{\sigma}& (p_2^*P\boxtimes_{\U(1)} p_1^*P)\boxtimes_{\U(1)} (p_2^*Q\boxtimes_{\U(1)} \boxtimes_{\U(1)} p_1^*Q) \nonumber\\&\xrightarrow{\lambda\boxtimes \mu}& c^*(P\boxtimes_{\U(1)} Q)\label{eq:gradedtensorproduct}
\eeq
where $\sigma$ is the braiding isomorphism $p_2^*Q\boxtimes_{\U(1)} p_1^*P\xrightarrow{\sim} p_1^*P\boxtimes_{\U(1)} p_2^*Q$ which carries a sign if the fibers of $P$ and $Q$ are both odd. 

\end{defn}

\begin{lem} \label{lem:BX}
Let $\bX_0\to \X$ be a Real atlas determining a groupoid presentation $\bX$ of a Real topological stack $\X$. There is a natural equivalence of the following 2-groupoids:
\begin{enumerate}
\item  Real graded  central extensions of $\bX$; and
\item Real graded  gerbes that trivialize on $\bX_0\to \X$.
\end{enumerate}
Furthermore, the monoidal structure on Real graded  gerbes is compatible with the tensor product of Real graded  central extensions. 
\end{lem}
\begin{proof}[Proof sketch]


The statement without gradings is well-known, e.g., \cite[Remark 5.5]{Heinloth} or~\cite[Proposition 4.12]{BX}; adding Real and graded structures presents no new difficulties. 
%
%
\ep

%

\begin{rmk}
We note that complex conjugation on~$\U(1)$ sends a graded central extension to its dual, in the sense that the $\U(1)$-principal bundles are dual to each other. Hence, a Real structure on a central extension (or gerbe) can be described in terms of an involutive isomorphism with its dual. This gives the description of Real gerbes as \emph{Jandl gerbes}, e.g., see~\cite[\S2.3]{orientifold}, \cite[\S3.3]{Waldorf2007} and~\cite[\S4.2]{NikolausSchweigert}. 
\end{rmk}

\begin{rmk}
Definitions of Real (or Jandl, or orientifold) gerbes in the literature vary in two essentially superficial ways. First, the 2-category of $C_2$-objects in topological stacks (where actions are weak) is equivalent a localization of the 2-category of topological groupoids with strict $C_2$-actions \cite[Proposition 1.5]{Romagny}. This allows one to specify either weak or strict equivariance data in the definition of a Real gerbe over a topological groupoid. Unless stated otherwise, we assume that $C_2$-actions are weak. The second variation in the definitions of Real gerbes uses the equivalence of 2-categories between $C_2$-stacks and stacks over $\pt\sq C_2$, where the latter leads to the definition of an orientifold, see~\S\ref{sec:Realstacks} and Remarks~\ref{rmk:orientifold} and~\ref{rmk:orientifold2}. Jandl gerbes and orientifolds were shown to be equivalent as 2-stacks in~\cite[Corollary A.2.4]{NikolausSchweigert}. \end{rmk}

\bibliographystyle{amsalpha}
\bibliography{references}

\newcommand{\etalchar}[1]{$^{#1}$}
\providecommand{\bysame}{\leavevmode\hbox to3em{\hrulefill}\thinspace}
\providecommand{\MR}{\relax\ifhmode\unskip\space\fi MR }
\providecommand{\MRhref}[2]{%
  \href{http://www.ams.org/mathscinet-getitem?mr=#1}{#2}
}
\providecommand{\href}[2]{#2}
\begin{thebibliography}{BMMS86}

\bibitem[ABG{\etalchar{+}}14]{ABGHR}
M.~Ando, A.~Blumberg, D.~Gepner, M.~Hopkins, and C.~Rezk, \emph{Units of ring
  spectra, orientations and {T}hom spectra via rigid infinite loop space
  theory}, J. Topol. \textbf{7} (2014), no.~4, 1077--1117.

\bibitem[ABS64]{ABS}
M.~Atiyah, R.~Bott, and A.~Shapiro, \emph{Clifford modules}, Topology
  \textbf{3} (1964), no.~suppl. 1, 3--38.

\bibitem[AS69]{Atiyah_Segal_complete}
M.~Atiyah and G.~Segal, \emph{Equivariant {K-theory} and completion}, J.
  Differential Geom. \textbf{3} (1969), no.~1-2, 1--18.

\bibitem[AS04]{AtiyahSegaltwists}
\bysame, \emph{{Twisted K-theory}}, Ukrainian Math. Bull. \textbf{1} (2004).

\bibitem[AS05]{AtiyahSegaltwistco}
\bysame, \emph{{Twisted K-theory and cohomology}}, {ArXiv Mathematics e-prints}
  (2005).

\bibitem[Ati66a]{Atiyahreal}
M.~Atiyah, \emph{K-theory and reality}, The Quarterly Journal of Mathematics
  \textbf{17} (1966).

\bibitem[Ati66b]{AtiyahPow}
M.~F. Atiyah, \emph{Power operations in {$K$}-theory}, Quart. J. Math. Oxford
  Ser. (2) \textbf{17} (1966), 165--193. \MR{202130}

\bibitem[Ati68]{Atiyahbottp}
M.~Atiyah, \emph{Bott periodicity and the index of elliptic operators},
  Quarterly Journal of Mathematics \textbf{19} (1968), 113--140.

\bibitem[BGNX07]{stringtopstacks}
K.~Behrend, G.~Ginot, B.~Noohi, and P.~Xu, \emph{String topology for loop
  stacks}, Comptes Rendus. Math\'ematique \textbf{344} (2007), no.~4, 247--252.

\bibitem[BM00]{Bouwknegt}
P.~Bouwknegt and V.~Mathai, \emph{{D-branes, B fields and twisted K theory}},
  JHEP \textbf{03} (2000), 007.

\bibitem[BMMS86]{bmms}
R.~R. Bruner, J.~P. May, J.~E. McClure, and M.~Steinberger, \emph{{$H_\infty $}
  ring spectra and their applications}, Lecture Notes in Mathematics, vol.
  1176, Springer-Verlag, Berlin, 1986. \MR{836132}

\bibitem[BMRS06]{BMRS}
J.~Brodzki, V.~Mathai, J.~Rosenberg, and R.~Szabo, \emph{D-branes, {RR}-fields
  and duality on noncommutative manifolds}, Communications in Mathematical
  Physics \textbf{277} (2006).

\bibitem[BX06]{BX}
K.~Behrend and P.~Xu, \emph{Differentiable stacks and gerbes}, preprint (2006).

\bibitem[CW08]{Carey}
A.~Carey and B.L. Wang, \emph{{Thom isomorphism and push-forward map in twisted
  K-theory}}, Journal of K-Theory \textbf{1} (2008).

\bibitem[DEF{\etalchar{+}}99]{strings1}
P.~Deligne, P.~Etingof, D.~Freed, L.~Jeffrey, D.~Kazhdan, J.~Morgan,
  D.~Morrison, and E.~Witten, \emph{{Quantum Fields and Strings: {A} Course for
  Mathematicians, Volume 1}}, American Mathematical Society, 1999.

\bibitem[DK70]{DonovanKaroubi}
P.~Donovan and M.~Karoubi, \emph{Graded brauer groups and {K}-theory with local
  coefficients}, Inst. Hautes Etudes Sci. Publ. Math. (1970), no.~38.

\bibitem[Dol63]{Doldpartition}
A.~Dold, \emph{Partitions of unity in the theory of fibrations}, Ann. of Math.
  \textbf{78} (1963).

\bibitem[EU14]{EspinozaUribe}
J.~Espinoza and B.~Uribe, \emph{Topological properties of the unitary group},
  JP Journal of Geometry and Topology \textbf{16} (2014).

\bibitem[FH21]{FreedHopkins}
D.~Freed and M.~Hopkins, \emph{Reflection positivity and invertible topological
  phases}, Geom. Topol. \textbf{25} (2021).

\bibitem[FHT11a]{FHTI}
D.~Freed, J.~Hopkins, and C.~Teleman, \emph{Loop groups and twisted
  {$K$}-theory {I}}, J. Topol. \textbf{4} (2011), no.~4, 737--798.

\bibitem[FHT11b]{FHTIII}
\bysame, \emph{Loop groups and twisted {$K$}-theory {III}}, Ann. of Math. (2)
  \textbf{174} (2011), no.~2, 947--1007.

\bibitem[FM13]{FreedMoore}
D.~Freed and G.~Moore, \emph{Twisted equivariant matter}, Ann. Henri Poincare
  \textbf{14} (2013).

\bibitem[Fre12]{Vienna}
D.~Freed, \emph{Lectures on twisted {K}-theory and orientifolds}, lectures at
  {ESI} {Vienna} (2012).

\bibitem[Gaj97]{Gajer}
P.~Gajer, \emph{Geometry of {Deligne} cohomology}, Inventiones mathematicae
  \textbf{127} (1997).

\bibitem[Gir71]{Giraud}
J.~Giraud, \emph{Cohomologie non ab\'elienne}, Springer-Verlag, 1971.

\bibitem[Gom23]{Gomi}
K.~Gomi, \emph{{Freed--Moore} {K}-theory}, To appear in: Communications in
  Analysis and Geometry (2023).

\bibitem[GSW11]{orientifold}
K.~Gawedzki, R.~Suszek, and K.~Waldorf, \emph{Bundle gerbes for orientifold
  sigma models}, Adv. Theor. Math. Phys. \textbf{15} (2011), no.~3.

\bibitem[Guk00]{GukovDbrane}
S.~Gukov, \emph{{K theory, reality, and orientifolds}}, Commun. Math. Phys.
  \textbf{210} (2000), 621--639.

\bibitem[Hei05]{Heinloth}
J.~Heinloth, \emph{Notes on differentiable stacks}, In: Mathematisches
  Institut, Georg-August-Universit\"at G\"ottingen: Seminars Winter Term
  2004/2005 (2005).

\bibitem[HJ19]{HebestreitJoachim}
F.~Hebestreit and M.~Joachim, \emph{Twisted spin cobordism and positive scalar
  curvature}, Journal of Topology \textbf{13} (2019), no.~1, 1--58.

\bibitem[HK24]{ZachYigal}
Z.~Halladay and Y.~Kamel, \emph{Real spin bordism and orientations of
  topological {K}-theory}, preprint (2024).

\bibitem[HMSV16]{HMSV}
P.~Hekmati, M.~Murray, R.~Szabo, and R.~Vozzo, \emph{Real bundle gerbes,
  orientifolds and twisted kr-homology}, Advances in Theoretical and
  Mathematical Physics \textbf{23} (2016).

\bibitem[Kan07]{Kankaanrinta}
M.~Kankaanrinta, \emph{{Equivariant collaring, tubular neighbourhood and gluing
  theorems for proper Lie group actions}}, Algebraic \& Geometric Topology
  \textbf{7} (2007), no.~1, 1 -- 27.

\bibitem[Kap00]{KapustinDbrane}
A.~Kapustin, \emph{{D-branes in a topologically nontrivial B field}}, Adv.
  Theor. Math. Phys. \textbf{4} (2000), 127--154.

\bibitem[Kub16]{Kubota}
Y.~Kubota, \emph{Notes on twisted equivariant {K}-theory for {$C^*$-algebras}},
  International Journal of Mathematics \textbf{27} (2016).

\bibitem[LM89]{LM}
H.B. Lawson and M.L. Michelsohn, \emph{Spin geometry (pms-38)}, Princeton
  University Press, 1989.

\bibitem[Mat71]{Matumoto}
T.~Matumoto, \emph{Equivariant {K}-theory and {Fredholm} operators}, J. Fac.
  Sci. Univ. Tokyo Sect. I A Math. \textbf{18} (1971).

\bibitem[MM97]{Minasian}
R.~Minasian and G.~Moore, \emph{{{K}-theory and Ramond-Ramond charge}}, JHEP
  \textbf{11} (1997), 002.

\bibitem[MMS03]{MMS}
V.~Mathai, M.~Murray, and D.~Stevenson, \emph{Type {I D}-branes in an {H}-flux
  and twisted {KO}-theory}, JHEP \textbf{7} (2003).

\bibitem[Moe02]{Moerdijk2002}
Ieke Moerdijk, \emph{Introduction to the language of stacks and gerbes},
  preprint (2002).

\bibitem[Mos57]{Mostow}
G.~D. Mostow, \emph{Equivariant embeddings in {Euclidean} space}, Ann. of Math.
  \textbf{65} (1957).

\bibitem[Mou11]{Moutuou}
E.M. Moutuou, \emph{Twistings of {KR} for real groupoids}, preprint (2011).

\bibitem[Mou12]{MoutuouThesis}
\bysame, \emph{Twisted groupoid {KR}-theory}, Universit\'e de Lorraine (2012).

\bibitem[Noo12]{Noohiclassifying}
B.~Noohi, \emph{Homotopy types of topological stacks}, Advances in Mathematics
  \textbf{230} (2012).

\bibitem[NS11]{NikolausSchweigert}
T.~Nikolaus and C.~Schweigert, \emph{Equivariance in higher geometry}, Advances
  in Mathematics \textbf{226} (2011), no.~4, 3367--3408.

\bibitem[Pal61]{Palais}
R.~Palais, \emph{On the existence of slices for actions of non-compact {Lie}
  groups}, Ann. of Math. \textbf{73} (1961).

\bibitem[Pro96]{Pronk}
D.~Pronk, \emph{Etendues and stacks as bicategories of fractions}, Compositio
  Math. \textbf{102} (1996).

\bibitem[Rom05]{Romagny}
M.~Romagny, \emph{{Group actions on stacks and applications}}, Michigan
  Mathematical Journal \textbf{53} (2005), no.~1.

\bibitem[Sch18]{Schottenloher}
M.~Schottenloher, \emph{The unitary group in its strong topology}, Advances in
  Pure Mathematics \textbf{8} (2018).

\bibitem[SP11]{SP11}
C.~J. Schommer-Pries, \emph{Central extensions of smooth 2-groups and a
  finite-dimensional string 2-group}, Geometry \& Topology \textbf{15} (2011).

\bibitem[SS21]{Sati2021EquivariantPI}
H.~Sati and U.~Schreiber, \emph{Equivariant principal infinity-bundles},
  preprint (2021).

\bibitem[ST04]{ST04}
S.~Stolz and P.~Teichner, \emph{What is an elliptic object?}, Topology,
  geometry and quantum field theory, London Math. Soc. LNS 308, Cambridge Univ.
  Press (2004), 247--343.

\bibitem[TXLG04]{TuStacks}
J.L. Tu, P.~Xu, and C.~Laurent-Gengoux, \emph{Twisted {K}-theory of
  differentiable stacks}, {Annales scientifiques de l'\'Ecole Normale
  Sup\'erieure}, Serie 4 \textbf{37} (2004), no.~6.

\bibitem[Wal07]{Waldorf2007}
K.~Waldorf, \emph{More morphisms between bundle gerbes}, Theory and
  Applications of Categories [electronic only] \textbf{18} (2007), 240--273.

\bibitem[Wit98]{WittenDbrane}
E.~Witten, \emph{{D-branes and K-theory}}, JHEP \textbf{12} (1998), 019.

\bibitem[Zun06]{Zung}
N.~Zung, \emph{Proper groupoids and momentum maps: linearization, affinity, and
  convexity}, Ann. Sci. \'Ecole Norm. Sup. \textbf{39} (2006).

\end{thebibliography}

\end{document}